\documentclass[12pt,reqno]{amsart}
\usepackage{amsmath}
\usepackage{amssymb}
\usepackage{amsthm}
\usepackage{mathrsfs}
\usepackage{latexsym} 
\usepackage{xcolor}
\usepackage{comment}
\usepackage{graphicx}
\usepackage{caption}

\newtheorem{theorem}{Theorem}[section]
\newtheorem{lemma}{Lemma}[section]
\newtheorem{corollary}{Corollary}[section]
\newtheorem{remark}{Remark}[section]
\newtheorem{definition}{Definition}[section]
\newtheorem{proposition}{Proposition}[section]
\newtheorem{example}{Example}[section]
\newtheorem{assumption}{Assumption}[section]
\numberwithin{equation}{section}
\newcommand{\bth}{\begin{theorem}}
\newcommand{\ethe}{\end{theorem}}

\newcommand{\bre}{\begin{remark}}
\newcommand{\ere}{\end{remark}}

\newcommand{\ble}{\begin{lemma}}
\newcommand{\ele}{\end{lemma}}
\newcommand{\bde}{\begin{definition}}
\newcommand{\ede}{\end{definition}}

\newcommand{\bco}{\begin{corollary}}
\newcommand{\eco}{\end{corollary}}

\newcommand{\bpr}{\begin{proposition}}
\newcommand{\epr}{\end{proposition}}

\newcommand{\bexer}{\begin{exercise}}
\newcommand{\eexer}{\end{exercise}}
\newcommand{\breh}{\begin{hint}}
\newcommand{\ereh}{\end{hint}}

\newcommand{\halmos}{\hfill \qed}

\newcommand{\bexam}{\begin{example}}
\newcommand{\eexam}{\end{example}}

\newcommand{\pr} {{\bf Proof.}}

\newcommand{\bfi}{\begin{fig}}
\newcommand{\efi}{\end{fig}}

\newcommand{\beao}{\begin{eqnarray*}}
\newcommand{\eeao}{\end{eqnarray*}\noindent}

\newcommand{\beam}{\begin{eqnarray}}
\newcommand{\eeam}{\end{eqnarray}\noindent}
\newcommand{\E}{\mathbf{E}}
\newcommand{\PP}{\mathbf{P}}

\newcommand{\xto}{x\to\infty}

\newcommand{\bF}{\overline{F}}
\newcommand{\bB}{\overline{B}}
\newcommand{\bG}{\overline{G}}
\newcommand{\bV}{\overline{V}}

\newcommand{\bbr}{{\mathbb R}}

\newcommand{\bbb}{{\mathbb B}}

\newcommand{\bbn}{{\mathbb N}}

\newcommand{\vep}{\varepsilon}

\allowdisplaybreaks[1]
\begin{document}
\title[Uniform asymptotics for a multidimensional renewal risk model]{Uniform asymptotics 
for a multidimensional renewal risk model with random number of delayed claims and multivariate 
subexponentiality}

\author[D.G. Konstantinides, C.D. Passalidis, M. Yuan]{Dimitrios G. Konstantinides, Charalampos  D. Passalidis, Meng Yuan} 

\address{Dept. of Statistics and Actuarial-Financial Mathematics,
University of the Aegean, Karlovassi, GR-83 200 Samos, Greece}
\address{School of Data Sience and Artificial Inteligence,
Dongbei University of Finance and Economics, Dalian, Liaoning, China}
\email{konstant@aegean.gr} \email{sasd24009@sas.aegean.gr} \email{mengyuan\_prob@126.com}

\date{{\small \today}}

\begin{abstract}
In this paper we examine a multivariate risk model, with common renewal counting process, constant interest rate, and each 
claim vector is accompanied by a random number of delayed claim vectors. The interest is focused on the 
asymptotic behavior of the entrance probability of the discounted aggregate claims into some rare-sets, over a 
finite and an infinite time horizon. Our results study the case where the main claims and the delayed claims 
have in some sense, asymptotic equivalent tails, but also the case where the delayed claims are negligible with
comparisons with the main claims. More precisely, our estimations over finite time horizon are equipped with local uniformity, and are valid under the assumption of class $\mathcal{S}_A$, of multivariate subexponential distributions for the  
claim distributions. On the case of infinite time horizon we need a mild restriction on the distribution class $\mathcal{A}^*_A$ of multivariate subexponential 
distributions with positive lower Karamata index. The asymptotic relations reflect completely as all the sources of randomness, under the concrete 
rare-sets $A$, and the different dependence structures as well, without loosing elegance in spite of their generality.
Further, we provide some more explicit formulas, together with relaxations of some assumptions, for the claim distributions 
from the multivariate regular variation. For the proof of the main results on infinite time case and for the construction of 
examples of multivariate distributions we need some closure properties of distribution classes $\mathcal{A}^*$ and $\mathcal{A}^*_A$. Especially, we present some 
necessary and sufficient conditions for the closure property with respect to convolution and some sufficient conditions for the closure property with respect to product convolution. Finally, we carry out some numerical studies  to show the accuracy of our asymptotic estimations. The examples contain also cases with moderate heavy tails. 
\end{abstract}

\maketitle
%\vspace{3mm}
\vspace{3mm}


%\textit{Journal of Economic Literature}: Primary C46 ;\quad Secondary G32.

\section{Introduction} \label{sec.KPY.1}

\subsection{Model description} \label{subs.KPY.1.1}

In modern insurance companies the claims are split into two types, the main and the delayed claims. The first type contains claims that are produced 
directly by an event, for example in case of natural disaster, such type of claims describe the compensation of the property elements, for example 
immobility insurance, or car insurance,  or of the health complications, for example injure expenses, or instant death coverage. The second type of 
claims includes the delayed claims, that were produced after some time delay, by the event, that causes the main claim. In case of natural disaster 
such type of claims can include medical expenses, or death coverage after some curing period.

This way, we find an increasing number of papers focused on risk models with delayed claims, since any ignorance of them, eventually leads to 
underestimation of the insurance risk, see relative discussions and literature overview on risk models with delayed claims in 
subsection 1.2. 

Further, the claim distributions with heavy tails and dependence among the portfolios, attracts the attention of an increasing number of researchers 
in the area of multivariate risk models with heavy-tailed claims, see for example \cite{jiang:wang:chen:xu:2015}, \cite{gao:yang:2014}, \cite{li:2016}, 
\cite{yang:su:2023}, \cite{cheng:konstantinides:wang:2024}, \cite{chen:konstantinides:passalidis:2025}, \cite{konstantinides:passalidis:2024j}, 
\cite{konstantinides:liu:passalidis:2025}, \cite{konstantinides:passalidis:2025o}, among others.

In combination of the previous attempts, we consider an insurer that operates $d$-lines of business, with $d \in \bbn$, and each main claim vector, 
generates random number of delayed claim vectors. Concretely, the main claim vectors are depicted as sequence $\{{\bf X}^{(i)}\,,\; i \in \bbn\}$ of 
$d$-dimensional random vectors, that takes values in the non-negative orthant and their arrival times are denoted by the sequence $\{\tau_i\,,\; i \in \bbn\}$, 
with $\tau_0=0$ conventionally. The sequence of $\{\tau_i\,,\; i \in \bbn\}$ constitutes a counting process 
$N(t):=\sup \{i \in \bbn\;:\; \tau_i \leq t\}$, for any $t\geq 0$, where we put $\sup \emptyset =0$ conventionally. 
Let us suppose that the $\{N(t)\,,\; t\geq 0\}$ has finite renewal function of the form
\beao
\lambda(t) :=\E\left[ N(t)\right] = \sum_{i=1}^{\infty} \PP[\tau_i \leq t]\,,
\eeao
for any fixed $t \geq 0$.

To avoid trivial cases, we study time horizons from the interval $\Lambda = \{t\;:\;\lambda(t) >0\}$, and we assume that each 
vector ${\bf X}^{(i)}=(X_{1}^{(i)},\,\ldots,\,X_{d}^{(i)})^{\top}$ (where ${\bf z}^{\top}$ is the transposed vector 
of ${\bf z}$) can have zero components, but not all of them equal to zero, in order to get a renewal epoch at moment $\tau_i$.

Each main claim vector ${\bf X}^{(i)}$, generates $M_i$ number of delayed claim vectors, where $M_i$ represents a random variable with values from 
$\bbn_0:=\bbn \cup \{0\}$. These delayed claim vectors are denoted by ${\bf Y}^{(i,j)} =(Y_1^{(i,j)},\,\ldots,\,Y_d^{(i,j)})^{\top}$, 
with $1\leq j \leq M_i$, and with support of their distribution, contained in non-negative orthant, while their arrival times are $\tau_i + D_{ij}$, 
where $D_{ij}$ represent a non-negative random variable, that denotes the delay time of the $j$-th delayed claim vector, generated from the $i$-th 
main claim vector. Again, to avoid trivialities, for any $i=1,\,\ldots,\,N(t)$, and any  $1\leq j \leq M_i$, the ${\bf Y}^{(i,j)}$ can have 
zero components, but not all of them equal to zero. However, we notice that a zero component of ${\bf X}^{(i)}$, let say  ${\bf X}_1^{(i)}=0$, 
does not imply necessarily that the first component of ${\bf Y}^{(i,j)}$ is also equal to zero, which is very common in insurance context. 
For example, a traffic accident can generate main claim in the insurance portfolio of the car and zero  claim to the corresponding life insurance, 
at the moment $\tau_i$, but after some time interval $D_{i1}$ can lead to claim in the life insurance portfolio.

Thus, if $r\geq 0$ represents the constant interest rate, then the insurer's discounted aggregate claims up to moment $t \geq 0$, are obtained by the formula
\beam \label{eq.KPY.1.1} \notag
{\bf D}_r(t) &:=& \sum_{i=1}^{N(t)} {\bf X}^{(i)}\,e^{-r\,\tau_i} +  \sum_{i=1}^{N(t)}\sum_{j=1}^{M_i} {\bf Y}^{(i,j)}\,e^{-r\,(\tau_i + D_{ij})}\,{\bf 1}_{\{\tau_i + D_{ij} \leq t\}} \\[2mm]
&=&\begin{pmatrix}   
\sum_{i=1}^{N(t)} X_{1}^{(i)}\,e^{-r\,\tau_i} +  \sum_{i=1}^{N(t)}\sum_{j=1}^{M_i} Y_{1}^{(i,j)}\,e^{-r\,(\tau_i + D_{ij})}\,{\bf 1}_{\{\tau_i + D_{ij} \leq t\}}   \\ \vdots \\ \sum_{i=1}^{N(t)} X_d^{(i)}\,e^{-r\,\tau_i} +  \sum_{i=1}^{N(t)} \,\sum_{j=1}^{M_i}  Y_d^{(i,j)}\,e^{-r\,(\tau_i + D_{ij})}\,{\bf 1}_{\{\tau_i + D_{ij} \leq t\}} 
\end{pmatrix}\,,
\eeam
where ${\bf 1}_{\{E\}}$, denotes the indicator function of some event $E$. We note that the discounted aggregate claims on infinite time 
horizon are given as in relation \eqref{eq.KPY.1.1} after replacement of $t$ by $\infty$, keeping in mind that $N(\infty)= \infty$. Our 
goal is to find the asymptotic behavior of the following probabilities
\beam \label{eq.KPY.1.2} 
\PP\left( {\bf D}_r(t) \in x\,A \right)\,,\qquad \qquad \PP\left( {\bf D}_r(\infty) \in x\,A \right)\,,
\eeam
as $\xto$, that represent the entrance probability of discounted aggregate claims, over finite or infinite time horizon, respectively, 
into some rare set $x\,A$, where $A$ belongs to a general family of sets, that contains many interesting events for the actuarial practice, 
see $\mathscr{R}$ in Section 2. As we will discuss later, see Remark \ref{rem.KPY.3.1}, our asymptotic expressions for 
probabilities in \eqref{eq.KPY.1.2}, are directly connected with finite-time and infinite time ruin probabilities, respectively.

We establish asymptotic estimations for the entrance probability over finite time horizon that are locally uniform with 
respect to time, when the claim distribution belongs to class of multivariate subexponential distributions $\mathcal{S}_A$, while we obtain 
infinite time asymptotics (as also some globally uniform results) under a weak restriction to distribution class $\mathcal{A}_A^*$, in the 
sense that this class contains most of the important distributions of class $\mathcal{S}_A$. We should mention that even in the one-dimensional 
case with $A=(1,\,\infty)$, our results remain new, see also the discussions on remarks after the main results. 

Let us present two assumptions, that we use in the rest of the paper.

\begin{assumption} \label{ass.KPY.1.1}
The sequences $\{{\bf X}^{(i)}\,,\; i \in \bbn\}$, $\{{\bf Y}^{(i,j)}\,,\; i \in \bbn\,,\; j \in \bbn\}$, $\{M_{i}\,,\; i \in \bbn\}$, 
$\{D_{ij}\,,\; i \in \bbn\,,\; j \in \bbn\}$ and $\{N(t)\,,\; t\geq 0\}$ are mutually independent.  
\end{assumption}

\begin{assumption} \label{ass.KPY.1.2}
The sequence $\{\theta_{i}:=\tau_i - \tau_{i-1}\,,\; i \in \bbn\}$ are independent and identically distributed (i.i.d.) random variables. 
The $\{{\bf X}^{(i)}\,,\; i \in \bbn\}$, $\{{\bf Y}^{(i,j)}\,,\; i \in \bbn\,,\; j \in \bbn\}$ are sequences of i.i.d. random vectors with 
common distribution $F$, $G$, respectively and $\{M_{i}\,,\; i \in \bbn\}$, $\{D_{ij}\,,\; i \in \bbn\,,\; j \in \bbn\}$ are sequences of i.i.d. 
random variables with common distribution $K$, $H$, respectively.  
\end{assumption}

\bre \label{rem.KPY.1.1}
Assumptions \ref{ass.KPY.1.1} and \ref{ass.KPY.1.2}, look as the simplest possible for the risk model \eqref{eq.KPY.1.1}, that seems already 
complex enough. For example, somebody could relax the independence assumption between $\{{\bf X}^{(i)}\,,\; i \in \bbn\}$, 
$\{{\bf Y}^{(i,j)}\,,\; i \in \bbn\,,\; j \in \bbn\}$, or even to relax the assumption that $\{N(t)\,,\; t\geq 0\}$ is renewal process. 
However, although such generalizations could be feasible in cases when the distributions $F$ and $G$ belong to the class of multivariate 
regular variation, symbolically $MRV$, or even in cases when there are appropriate conditions on Matuszewska indexes, the fact that we study 
the whole class of multivariate subexponential distributions $\mathcal{S}_A$, and the slightly smaller class $\mathcal{A}_A^*$, such kind of 
generalizations are not easy, as it happens to be in smaller classes. According to our knowledge, for a random number of delayed claims in 
multidimensional setup, there is only one paper \cite{yuan:lu:fu:2025}, with arbitrarily dependent counting processes, not necessarily renewal, 
and with L\'{e}vy processes for logarithmic returns, but with $F$ and $G$ from class $MRV$, with asymptotic dependent components and non-uniform 
asymptotic expressions with respect to time.  

In this paper although we examine a more simple model, with a common renewal counting process and deterministic interest rate, we proceed further to uniform asymptotic estimates on class $\mathcal{S}_A$, without assumption that the vectors ${\bf X}$ and ${\bf Y}$ have necessarily asymptotic dependent components. Further, we consider also the case, when the delayed claims have 'negligible' distribution tails with respect to that of the main claims.
\ere

The rest of the paper is organized as follows. In subsection $1.2$ we present a short overview of the delayed risk models in one- or multi-dimensions, 
in order to make clear the practicability of our results. In Section $2$, we demonstrate the necessary preliminary concepts, related with heavy-tailed 
distribution classes and related indexes. In Section $3$, we present the first main result, which is referred to local uniform asymptotic expressions 
for the first probability in \eqref{eq.KPY.1.2}, when the distributions $F$ and $G$ either belong to class $\mathcal{S}_A$ and are weakly equivalent, 
or $F$ belongs to $\mathcal{S}_A$ and $G$ has negligible tail with respect to the tail of $F$. The proof is accommodated in the same section. In Section 
$4$, we give the second main result, where we provide the estimation of the entrance probability for infinite time horizon, under the restriction of 
the claim distributions into class $\mathcal{A}_A^*$. The proof and the preliminary lemmas are included in this section. The lemmas have their own merit, 
since they imply some closure properties of the classes $\mathcal{A}_A^*$ and $\mathcal{A}^*$, with respect to convolution and product convolution. 
In the $MRV$ subcase we take more explicit expressions for both results and some other conditions can be relaxed.  
We also provide more explicit forms for the asymptotic estimations of the probabilities on \eqref{eq.KPY.1.2}, (in non-$MRV$ case) 
for some concrete sets $A$ which are interesting in insurance practice, under some weak and general dependence structure, see Remarks \ref{rem.KPY.3.3}, 
\ref{rem.KPY.3.4}, \ref{rem.KPY.4.2} below. 
In Section 5, we provide some numerical studies in order to demonstrate the accuracy of our asymptotic results. 
In comparison to many similar studies in risk theory, we provide also some examples with moderate heavy-tailed marginals, like lognormal distribution. 
Such approximations often had been avoided, because moderate heavy tailed distributions gives smaller accuracy.

All the main results highlights the multivariate linear single big jump principle in the asymptotic behavior of the entrance probabilities from \eqref{eq.KPY.1.2}.

\subsection{Brief overview of delayed risk model} \label{subs.KPY.1.2}

The delayed risk model is a more realistic form of insurance risk, that fits better to the most of actuarial applications. 
In \cite{waters:papatriandafylou:1985} was introduced such a model for the discrete time case, while for continuous time the first 
papers are \cite{yuen:guo:2001}, \cite{yuen:guo:ng:2005}. However, these papers focused in the case of claim distributions with light tails. 

In \cite{li:2013} we find an application of the randomly weighted sums on the ruin probability over infinite time horizon for 
models with heavy-tailed distributions of main and delayed claims. Since then, we meet a spectrum of papers that concentrate on this topic, 
relaxing either dependence conditions or distribution classes, see \cite{gao:zhuang:huang:2019}, \cite{yang:li:2019}, \cite{yang:wang:zhang:2021}, 
\cite{yuan:lu:2023}, \cite{lu:qin:yuan:2025} among others. These papers above, consider the case where $M_i=1$ holds almost surely, 
and also consider one or two dimensional risk models. However, the generation of exactly one delayed claim, 
seems too restrictive in modern insurance industry, where the number of delayed claims should be random. Hence, in \cite{li:2023b} were 
examined ruin probabilities, and discounted aggregate claims, over finite and infinite time horizon, in a renewal risk model (one-dimensional), with Assumptions 
\ref{ass.KPY.1.1} and \ref{ass.KPY.1.2} to remain intact. However, the estimations over the finite time horizon were not uniform. In two-dimensional set up, in 
\cite{jia:chen:cheng:2025} we have extension to the study of finite time, under some weak dependence structures. In \cite{liu:gao:chen:2024} we find a substantial 
extension of the result in \cite{li:2023b}, in one-dimensional set up, through a quasi-renewal process and with uniformity, however it was restricted to a smaller 
distribution class than the subexponential one.

In multidimensional set up, we know only \cite{yuan:lu:fu:2025}, and as was mentioned these results do not overlap completely ours, due to the more general model 
there, with $d$ different counting processes, arbitrarily dependent and not necessarily renewal ones, as also the presence of stochastic interest rate. However, that 
paper is focused in class $MRV$ for the $F,\,G$, with the ${\bf X}$ and ${\bf Y}$ to have asymptotically dependent components, while the $F,\,G$ have a relation of 
asymptotic equivalence between the tails. 

In the present paper, we examine the case, when $F,\,G$ belong to class $\mathcal{S}_A$ (or to class $\mathcal{A}_A^*$), which contains the $MRV$ and which are more 
wide, since they include also moderate heavy-tailed marginals, while we do not necessarily assume that the components of the vectors ${\bf X}$ and 
${\bf Y}$ are asymptotically dependent. Even more, we examine also the cases where the tail of $G$ is asymptotic negligible with respect to the tail of $F$, and 
these results show that the influence of the delayed claims to entrance probabilities \eqref{eq.KPY.1.1} and \eqref{eq.KPY.1.2} is asymptotically negligible.

Hence, the present paper is inspired by \cite{li:2023b} and \cite{yuan:lu:fu:2025}. Our results remain new even with the respect to one-dimensional 
subcase, with $A=(1,\,\infty)$, see in Remarks \ref{rem.KPY.3.1}, and \ref{rem.KPY.4.1}.

\section{Preliminaries} \label{sec.KPY.2}

In this section we present some necessary concepts, related with the heavy-tailed distribution classes on multivariate set up, as well as some 
examples of multivariate distribution. 

In what follows, all the vectors are of dimension $d$, with $d \in \bbn$, which are denoted in bold script, and the operations are defined 
component-wisely, namely ${\bf z} \pm {\bf y} =(z_1\pm y_1,\,\ldots,\,z_d\pm y_d)^{\top}$ for any vectors ${\bf z}$ and ${\bf y}$, and 
the scalar product of ${\bf z}$ with some finite scalar $k$, is defined as usual by $k\,{\bf z} =(k\,z_1,\,\ldots,\,k\,z_d)^{\top}$. 
For any real numbers $a_1,\,\ldots,\,a_n$ we denote their maximum as $\bigvee_{i=1}^n a_i:= \max \{a_1,\,\ldots,\,a_n\}$ and  their minimum 
as $\bigwedge_{i=1}^n a_i:= \min \{a_1,\,\ldots,\,a_n\}$. We note that the sum over an empty index set is equal to zero. 
For a set $\bbb \subset \bbr^d=(-\infty,\,\infty)^d$, we denote by $\overline{\bbb}$ its closed hull, by $\bbb^c$ its complement set, 
by $\partial \bbb$ its border. A set $\bbb \subset \bbr^d$ is called increasing if for any ${\bf z} \in \bbb$ and any 
${\bf y} \in \bbr_+^d:=[0,\,\infty)^d$, it holds ${\bf z} + {\bf y} \in \bbb$.

Hereafter, all the the limit relations by default hold as $\xto$, except otherwise stated. For two $d$-variate positive functions ${\bf f}$, ${\bf g}$, we 
write ${\bf f}(x\,\bbb) \sim c\,{\bf g}(x\,\bbb)$, with some $c \in (0,\,\infty)$, ${\bf f}(x\,\bbb) =o\left[{\bf g}(x\,\bbb)\right]$, and 
${\bf f}(x\,\bbb) =O\left[{\bf g}(x\,\bbb)\right]$, with $\bbb \in \bbr_+^d \setminus \{{\bf 0}\}$, if it holds
\beao
\lim\dfrac{{\bf f}(x\,\bbb)}{{\bf g}(x\,\bbb)}=c\,,\quad \lim\dfrac{{\bf f}(x\,\bbb)}{{\bf g}(x\,\bbb)}=0\,,\quad \limsup \dfrac{{\bf f}(x\,\bbb)}{{\bf g}(x\,\bbb)}< \infty\,, 
\eeao
respectively, and we denote by ${\bf f}(x\,\bbb) \asymp {\bf g}(x\,\bbb)$, if both relations ${\bf f}(x\,\bbb) =O\left[{\bf g}(x\,\bbb)\right]$ and ${\bf g}(x\,\bbb) =O\left[{\bf f}(x\,\bbb)\right]$ are true. For two $(d+1)$-variate functions ${\bf f}^*$ and ${\bf g}^*$, we write ${\bf f}^*(x\,\bbb;\,y) \sim c\,{\bf g}^*(x\,\bbb;\,y)$, uniformly for $y \in \Delta$, with $\Delta$ a non-empty set and $c \in (0,\,\infty)$, if it holds
\beao
\lim \sup_{y \in \Delta} \left|\dfrac{{\bf f}^*(x\,\bbb;\,y)}{{\bf g}^*(x\,\bbb;\,y)}- c \right| =0\,.
\eeao
Finally, if a random variable (or, vector) $\Theta_1$ follows distribution $V_1$, symbolically we write $\Theta_1 \stackrel{d}{\sim} V_1$, 
while by $V_1*V_2$ we denote the convolution of distribution $V_1$ and $V_2$. By $V_1^{n*}$ we denote the $n$-th order convolution power of $V_1$. 
If $V_1$ is one-dimensional we denote by $\bV_1(x) =1-V_1(x)$, for $x \in \bbr$, its distribution tail.

\subsection{Heavy-tailed distribution classes} \label{subs.KPY.2.1}

The multivariate heavy-tailed classes play important role in actuarial and financial mathematics, mostly in modeling the dependence effect, that appears 
in real world applications. Special attention was gotten to the class of multivariate regularly varying distributions ($MRV$) that has a variety of applications 
in mathematical statistics as well as in applied probability, see \cite{resnick:2007}, \cite{buraczewski:damek:mikosch:2016}, \cite{samorodnitsky:2016} among others.

During the last years we observed an attempt to study properties and applications of bigger (than $MRV$) classes, mostly motivated by seminal paper 
\cite{samorodnitsky:sun:2016}, where was introduced the following family of sets
\beam \label{eq.KPY.2.1} 
\mathscr{R}:=\left\{ A \subsetneq \bbr^d\;:\; A\; \text{open,\,increasing},\;A^c \;\text{convex},\;{\bf 0} \notin \overline{A} \right\}\,.
\eeam
The family $\mathscr{R}$ is closed with respect to positive scalar multiplication.

\bre \label{rem.KPY.2.1}
The set family \eqref{eq.KPY.2.1} contains interesting sets, for the actuarial practice, as the following
\beam \label{eq.KPY.2.2} 
A_1=\left\{ {\bf y}\;:\; \sum_{j=1}^d l_j\,y_j >c \right\}\,,
\eeam
with $l_1,\,\ldots,\,l_d \geq 0$, $\sum_{j=1}^d l_j =1$ and $c>0$, and
\beam \label{eq.KPY.2.3} 
A_2=\left\{ {\bf y}\;:\; y_{j} >c_j\,,\; \exists \;j=1,\,\ldots,\,d  \right\}\,,
\eeam
with $c_j>0$, for all $j=1,\,\ldots,\,d$.
So, the probabilities in \eqref{eq.KPY.1.2} with $A=A_i$, for $i=1,\,2$, play crucial role in estimation of the insurer's solvency, in the sense of estimation 
of the initial capital adequacy in dynamic framework. We also observe that for $d=1$, set $A_2$ with $c_1=1$, is reduced to the set $(1,\,\infty)$ and hence 
the probabilities \eqref{eq.KPY.1.2} provide the distribution tail of the discounted aggregate claims (either over finite or over infinite time horizon, respectively). 
This probability is well-studied in risk theory and further in the applied probability area. For more examples of sets from family $\mathscr{R}$, see 
\cite[Sec. 4]{samorodnitsky:sun:2016} and \cite[Sec. 2]{shen:yuan:lu:2022}. 
\ere

Hereafter, the distributions $F$ have support from the non-negative orthant, and all the one-dimensional distributions have infinite 
right endpoint of their support.

Let $A \in \mathscr{R}$ and ${\bf Z} \stackrel{d}{\sim} V$, with $V$ supported in $\bbr_+^d$. From \cite[Lem. 4.5]{samorodnitsky:sun:2016} 
we obtain the random variable
\beam \label{eq.KPY.2.4} 
Z_A:=\sup \left\{ u\;:\; {\bf Z} \in u\,A \right\}\,,
\eeam
has proper distribution $V_A$, whose tail is given by 
\beam \label{eq.KPY.2.5} 
\bV_A(x)=\PP\left(\sup_{{\bf p} \in I_A} {\bf p}^{\top}\,{\bf Z} >x \right)=\PP\left({\bf Z} \in x\,A \right)\,,
\eeam
for any $x>0$, where $I_A \subsetneq \bbr^d$ represents an index set (see \cite[Lem. 4.3(c), 4.5]{samorodnitsky:sun:2016} for the 
existence of $I_A$ for any $A \in \mathscr{R}$, and of the proper distribution of \eqref{eq.KPY.2.4}, respectively).

Thus, we have the definition of the multivariate subexponentiality as follows. We say that $V$ belongs to the class of multivariate 
subexponential distributions on $A$, symbolically $V \in \mathcal{S}_A$, if $V_A$ is subexponential, symbolically $V_A \in \mathcal{S}$, that means it holds
\beao
\lim \dfrac{\overline{V_A^{n*}}(x)}{\bV_A(x)} =n\,,
\eeao
for all (or, equivalently, for some) integer $n\geq 2$. The class of multivariate subexponential distributions (on the whole $\mathscr{R}$) 
is defined as $\mathcal{S}_{\mathscr{R}}:=\bigcap_{A \in \mathscr{R}} \mathcal{S}_A$.

Through this way, \cite{konstantinides:passalidis:2024g} introduced and studied the class of multivariate long tailed distributions on $A$, 
symbolically $V \in \mathcal{L}_A$, when $V_A \in \mathcal{L}$, namely it holds
\beao
\lim \dfrac{\bV_A(x-y)}{\bV_A(x)} = 1\,, 
\eeao
for all (or, equivalently, for some) $y>0$, and $\mathcal{L}_{\mathscr{R}}:=\bigcap_{A \in \mathscr{R}} \mathcal{L}_A$, respectively. 

For an one-dimensional distribution $B$, with infinite right endpoint, its lower Karamata index is defined as
\beam \label{eq.KPY.2.6} 
K_B^{-}=- \lim_{v \downarrow 1} \dfrac {\log \overline{B^*}(v)}{\log v}\,,
\eeam
and its lower Matuszewska index is defined as
\beao
J_B^- = -\lim_{v\to \infty} \dfrac{\log \overline{B^*}(v)}{\log v}\,,
\eeao
where
\beao
\overline{B^*}(v)= \limsup \dfrac {\bB(v\,x)}{\bB(x)}\,,
\eeao
see in \cite[Subsec. 2.1]{bingham:goldie:teugels:1987} for more discussions. We note that if $\Theta_1 \stackrel{d}{\sim} B_1$ and there 
exists some constant $K>0$, with $K\,\Theta_1 \stackrel{d}{\sim} B_1'$, then we obtain $K_{B_1}^{-} = K_{B_1'}^{-}$. Further, if $\bB_1(x) \sim c\,\bB_2(x)$, 
with $c \in (0,\,\infty)$, then we find $K_{B_1}^{-} =K_{B_2}^{-}$.

For an one-dimensional distribution $B$, we say that $B \in \mathcal{A}^*$, if $B \in \mathcal{S}$ and $K_B^{-}> 0$. This definition appeared 
in \cite{tang:yuan:2016}, where was written that the most important distributions from class $\mathcal{S}$ have positive lower Karamata index, 
as for example the regularly varying and the rapidly varying distributions. The difference between the classes $\mathcal{A}^*$ and $\mathcal{S}$ 
is very small, however the restriction in $\mathcal{A}^*$ permits the study of infinite randomly weighted sums, that is impossible for the class 
$\mathcal{S}$, see \cite[Th. 3.1, 3.2]{tang:yuan:2016}. In \cite{konstantinides:passalidis:2025o} was introduced the class $\mathcal{A}_A^*$ as 
follows: $V \in \mathcal{A}_A^*$ if $V_A \in \mathcal{A}^*$, and further $\mathcal{A}_{\mathscr{R}}^* = \bigcap_{A \in \mathscr{R}} \mathcal{A}_A^*$.

We also denote by $\mathcal{T}^*$ the distribution class $\mathcal{L}$, with positive lower Karamata index, namely $B_1 \in \mathcal{T}^*$, if $B_1 \in 
\mathcal{L}$ and $K_{B_1}^- >0$. 

Here, we assume that the distributions $F$ and $G$ belong to class $\mathcal{A}_A^*$ to find the asymptotic expressions on infinite 
time horizon as well as some uniform asymptotic expressions for all time horizons. For this purpose we need some necessary preliminary lemmas, 
that are related with the closure properties with respect to convolution for distributions from classes $\mathcal{A}^*$ and $\mathcal{A}_A^*$ 
and with respect to product convolution, see in Section 4.2. This kind of properties are also helpful in construction of examples in the classes 
$\mathcal{A}_A^*$ and $\mathcal{A}_{\mathscr{R}}^*$. We notice, that for all distributions $B$, it holds $K_B^- \leq J_B^-$, and that $J_B^- >0$ 
if and only if $B \in \mathcal{P_D}$, namely if it holds
\beao
\limsup \dfrac{\bB(v\,x)}{\bB(x)} < 1\,,
\eeao
for all (or, equivalently, for some) $v>1$. Thus the class $\mathcal{A}:=\mathcal{S} \cap \mathcal{P_D}$, introduced in \cite{konstantinides:tang:tsitsiashvili:2002}, contains $\mathcal{A}^*$, although the difference is negligible. As before, we say $V \in \mathcal{A}_A$, if $V_A \in \mathcal{A}$, hence it holds $\mathcal{A}_A^* \subsetneq \mathcal{A}_A$. 

An one-dimensional distribution $B$, is called regularly varying with index $\alpha \in (0,\,\infty)$, symbolically $B \in \mathcal{R}_{-\alpha}$, 
if it holds
\beao
\lim \dfrac{\bB(b\,x)}{\bB(x)} =b^{-\alpha}\,,
\eeao 
for any $b >0$. From \cite[Theorem 1.5.2]{bingham:goldie:teugels:1987}, we know that the above relation holds 
uniformly in the following sense: for any fixed $\epsilon>0$ it holds
\beam \label{eq.KPY.2.6a} 
\lim \sup_{b \in [\epsilon,\infty) } \left|\dfrac{\bB(b\,x)}{b^{-\alpha}\bB(x)}- 1 \right| =0\,.
\eeam
Hence, we say that distribution $V$ belongs to the class $MRV$, if there exists one-dimensional distribution $B \in \mathcal{R}_{-\alpha}$, with $\alpha \in (0,\,\infty)$ and some non-degenerate to zero Radon measure $\mu$ such that 
\beam \label{eq.KPY.2.7} 
\lim \dfrac {\PP ({\bf Z} \in x\,\bbb )}{\bB(x)} = \mu(\bbb)\,,
\eeam
for any Borel set $\bbb \in \overline{\bbr_+^d}$, bounded away from the ${\bf 0}$, which satisfies $\mu(\partial \bbb) =0$. 
Then we denote $V \in MRV(\alpha,\,B,\,\mu)$.

Let us notice that in the proof of \cite[Prop. 4.14]{samorodnitsky:sun:2016} was found that for any $A \in \mathscr{R}$ we 
have $\mu(\partial A) = 0$, $\mu(A)< \infty$. Therefore if $V \in MRV(\alpha,\,B,\,\mu)$, with $\mu(A)>0$, then 
$V_A\in\mathcal{R}_{-\alpha}$ (the opposite does not holds in general) and hence by the one-dimensional inclusions
(see \cite[Ch. 2]{leipus:siaulys:konstantinides:2023}) we find that if $\mu(A)>0$ for all $A \in \mathscr{R}$, then it holds
\beam \label{eq.KPY.2.8} 
MRV \subsetneq \mathcal{A}_{\mathscr{R}}^* \subsetneq \mathcal{S}_{\mathscr{R}}  \subsetneq \mathcal{L}_{\mathscr{R}} \,,
\eeam
where $MRV$ represents the union of all $MRV(\alpha,\,B,\,\mu)$ distributions, see in 
\cite[Prop. 2.1, Rem. 2.2]{konstantinides:passalidis:2025o}. Relation \eqref{eq.KPY.2.8} is valid also for the 
classes $\mathcal{A}_A^*$, $\mathcal{S}_A$, $\mathcal{L}_A$, for any $A \in \mathscr{R}$ such that the assumption $\mu(A)>0$ is valid. We refer the reader to \cite[Exam. 2.1]{chen:passalidis:2028}, for the cases where $V \in MRV(\alpha,\,B,\,\mu)$ with $\mu(A)=0$, for some set $A \in \mathscr{R}$. However, under the $MRV$, $\mu(A_i)>0$, for $i=1,2$ is true, except in the cases where neither of the marginals is not regularly varying. We observe that the first 
inclusion of \eqref{eq.KPY.2.8} is not trivial, but contains a rich spectrum of distributions that have moderately heavy tails, 
which does not contain class $MRV$. The second inclusion of $\mathcal{A}_{\mathscr{R}}^* \subsetneq \mathcal{S}_{\mathscr{R}}$ 
contains trivial examples for practical applications. 

In next Section we construct some examples to show this type of inclusions. It is worth to mention that up to class 
$\mathcal{S}_{\mathscr{R}}$ in relation \eqref{eq.KPY.2.8}, the distributions are consistent with respect to multivariate 
linear single big jump principle, due to the property of insensitivity with respect to dimension (see also 
\cite{konstantinides:passalidis:2025c} for nonlinear single big jump principle).

\subsection{Examples in classes $\mathcal{S}_A$ and $\mathcal{A}_A^*$} \label{subs.KPY.2.2}

Here we present two examples for the distribution classes $\mathcal{S}_A$ and $\mathcal{A}_A^*$, for the sets $A_1$ and $A_2$ 
from relations \eqref{eq.KPY.2.2} and \eqref{eq.KPY.2.3} and an example for the class $\mathcal{A}_{\mathscr{R}}^*$. The examples 
for classes $\mathcal{S}_A$ and $\mathcal{A}_A^*$ are based only on conditions for the marginal distributions and the dependence 
structure among the components, that facilitates their control and construction. We see that these examples are not restricted to the $MRV$ case.

The first example is focused on class $\mathcal{A}_{A_1}^*$, with $A_1$ from relation \eqref{eq.KPY.2.2}, see 
in \cite[Exam. 4.1]{konstantinides:liu:passalidis:2025} for the case of class $\mathcal{S}_{A_1}$. We recall that 
for set $A_1$ a possible choice of $I_{A_1}$ is in the form 
$\left\{\left(l_1/c,\,\ldots,\,l_d/c \right)^{\top} \right\}$ and hence we obtain
\beam \label{eq.KPY.2.9} 
\bV_{A_1}(x) = \PP\left({\bf Z} \in x\,A_1 \right) = \PP\left(\sum_{i=1}^d \dfrac{l_i}c\,Z_i > x \right) \,,
\eeam 
for any $x>0$. The regression dependence, see below \eqref{eq.KPY.2.10}, that depicts the components dependence of vector ${\bf Z}$ in 
the following examples, contains several commonly used copulas and was introduced by \cite{lehmann:1966}, see also in \cite[p. 876]{geluk:tang:2009} 
for more information. 

\bexam \label{exam.KPY.2.1}
Let ${\bf Z}=(Z_1,\,\ldots,\,Z_d)^{\top}$ be a non-negative random vector with distribution $V$. Let the marginals $V_i \in \mathcal{A}^*$, with 
$i=1,\,\ldots,\,d$ and $V_i * V_j \in \mathcal{S}$, for any $1\leq i \neq j \leq d$. We assume that the $Z_1,\,\ldots,\,Z_d$ are regression 
dependent random variables, namely there exist constants $x_0 >0$, $C>0$, such that it holds
\beam \label{eq.KPY.2.10} 
\PP\left(Z_i > x_i\;|\; Z_j =x_j\,,\;j \in J_i \right) \leq C\,\PP\left(Z_i > x_i \right) \,,
\eeam
for any $x_i \wedge x_j \geq x_0$, and $\emptyset \neq J_i \subseteq \{1,\,\ldots,\,d\} \setminus \{i\}$, with $i=1,\,\ldots,\,d$. We also 
assume that $l_1,\,\ldots,\,l_d \in [a,\,b]^d$, with $0< a \leq b \leq 1$. Then, since $c>0$, from relation \eqref{eq.KPY.2.9}, under the 
previous assumptions, we can apply \cite[Th. 2]{wang:2011} to obtain the asymptotic formula
\beam \label{eq.KPY.2.11} 
\bV_{A_1}(x) \sim \sum_{i=1}^d \PP\left(\dfrac{l_i}c\,Z_i > x \right) =: \sum_{i=1}^d
\overline{V'_i}(x) \,,
\eeam
where, 
\beao
\dfrac{l_i}c\,Z_i \stackrel{d}{\sim} V_i'\,,
\eeao
for any $i=1,\,\ldots,\,d$. From \eqref{eq.KPY.2.11} and the fact that independence represents a special sub-case of \eqref{eq.KPY.2.10}, 
again through \cite[Th. 2]{wang:2011} we find 
\beam \label{eq.KPY.2.12} 
\overline{V_1'*V_2'}(x) \sim \overline{V'_1}(x) + \overline{V'_2}(x) \,,
\eeam
from which via the Lemma \ref{lem.KPY.4.1}, since $V_i' \in \mathcal{A}^*$, is implied the $V_1'*V_2' \in \mathcal{A}^*$. 

Furthermore, for the same reason, 
\beam \label{eq.KPY.2.13} 
\overline{V_1'*V_2'*V_3'}(x) \sim \overline{V_1'*V_2'}(x)+ \overline{V'_3}(x) \,,
\eeam
and hence via the Lemma \ref{lem.KPY.4.1}, since $V_1'*V_2' \in \mathcal{A}^*$ and $V_3' \in \mathcal{A}^*$ is implied that 
$V_1'*V_2'*V_3' \in \mathcal{A}^*$. Further, by induction we obtain
\beam \label{eq.KPY.2.14} 
\overline{V_1'* \cdots *V_d'}(x) \sim \sum_{i=1}^d \overline{V'_i}(x) \,,
\eeam
with $V_1'* \cdots *V_d' \in \mathcal{A}^*$. From this last relation, together with \eqref{eq.KPY.2.11} and \eqref{eq.KPY.2.14}, because of the 
closure property of $\mathcal{A}^*$ with respect to the strong tail equivalence, we get $V_{A_1}  \in \mathcal{A}^*$, to conclude that $V \in \mathcal{A}_{A_1}^*$.
\eexam

\bre \label{rem.KPY.2.2}
Example \ref{exam.KPY.2.1}, under the assumption $V_1,\,\ldots,\,V_d \in \mathcal{S}$, instead of the inclusion into $\mathcal{A}^*$, implies, 
following a similar argument, that $V \in \mathcal{S}_{A_1}$, using \cite[Cor. 1.1]{leipus:siaulys:2020}, instead of Lemma \ref{lem.KPY.4.1}. In 
fact, for $d=2$ and $c=1$, this is \cite[Exam. 4.1]{konstantinides:liu:passalidis:2025}. As was mentioned in \cite[Rem. 2.3]{konstantinides:passalidis:2025o}, 
the \cite[Exam. 4.2 - 4.4]{konstantinides:liu:passalidis:2025}, with the additional requirement of positive lower Karamata index of the marginal distributions, 
belong to $\mathcal{A}_{A_1'}^*$
with
\beao
A_1' = \left\{{\bf y}\;:\;\sum_{j=1}^2 l_j\,y_j >1 \right\}\,.
\eeao
However, via almost the same approach, we can check that they also belong to $\mathcal{A}_{A_1}^*$. Further, their Example 4.5 
also belongs to $\mathcal{A}_{A_1''}^*$, with $A_1''$ to be like the $A_1'$ with $l_1=l_2=1/2$.  
\ere

Next, we find sufficient conditions for distribution membership to class $\mathcal{S}_{A_2}$ and $\mathcal{A}_{A_2}^*$. 
A possible choice for set $I_{A_2}$, is the $\left\{ \left( 1/c_1,\,\ldots,\, 1/c_d\right)^{\top} \right\}$, and thus we obtain
\beam \label{eq.KPY.2.15} 
\overline{V_{A_2}}(x) = \PP\left({\bf Z} \in x\,A_2 \right) = \PP\left(\bigvee_{i=1}^d \dfrac{Z_i}{c_i} > x \right)\,,
\eeam

\bexam \label{exam.KPY.2.2}
Let ${\bf Z}=(Z_1,\,\ldots,\,Z_d)^{\top}$ be a non-negative random vector with distribution $V$. We assume that the assumptions of 
Example \ref{exam.KPY.2.1} are valid, with the only difference that $V_i \in \mathcal{B} \in \left\{ \mathcal{A}^*\,,\;\mathcal{S} \right\}$, 
for any $i=1,\,\ldots,\,d$. Then, since $0< 1/c_i < \infty$, applying \cite[Th. 2]{wang:2011} on \eqref{eq.KPY.2.15} we find
\beam \label{eq.KPY.2.16} 
\overline{V_{A_2}}(x) \sim \sum_{i=1}^d \PP\left(\dfrac{Z_i}{c_i} > x \right)=:\sum_{i=1}^d \overline{V'_i}(x) \,,
\eeam
with 
\beao
\dfrac{Z_i}{c_i} \stackrel{d}{\sim} V'_i \in \mathcal{B}\,.
\eeao 
As before, we obtain relation \eqref{eq.KPY.2.12} and it holds $V_1'* V_2' \in \mathcal{B}$, see Lemma \ref{lem.KPY.4.1} for 
$\mathcal{B} =\mathcal{A}^* $, and \cite[Cor. 1.1]{leipus:siaulys:2020} for $\mathcal{B} =\mathcal{S}$. Employing again induction 
we can show that \eqref{eq.KPY.2.14} is true and further that $V_1'* \cdots *V_d' \in \mathcal{B}$. From this last, together 
with relations \eqref{eq.KPY.2.14} and \eqref{eq.KPY.2.16} we have that $V_{A_2} \in \mathcal{B}$ (by the closure property 
of $\mathcal{B}$ with respect to strong tail equivalence), and therefore it holds
$V \in \mathcal{B}_A$.
\eexam 

Finally, we consider an example for distribution from class $\mathcal{A}_{\mathscr{R}}^*$, which is based on \cite[Exam. 4.17]{samorodnitsky:sun:2016}.

\bexam \label{exam.KPY.2.3}
Let ${\bf Z}=(Z_1,\,\ldots,\,Z_d)^{\top}$ be a non-negative random vector with $\PP({\bf Z} = {\bf 0})=0$. 
We define a random variable $\Theta \stackrel{d}{\sim} B$, whose distribution tail is given by
\beao
\bB(x) := \PP \left(\sum_{j=1}^d Z_j > x \;\bigg|\;\dfrac{{\bf Z}}{\sum_{j=1}^d Z_j} \in (\eta_1,\,\ldots,\,\eta_d) \right)\,,
\eeao
for any $(\eta_1,\,\ldots,\,\eta_d) \in \Delta_d$, with $\Delta_d$ representing the unit ball. We define the random variable
\beam \label{eq.KPY.2.17} 
W:= \inf \left\{u>0\;:\;u\,\dfrac{{\bf Z}}{\sum_{j=1}^d Z_j} \in A \right\}\,,
\eeam
for some fixed $A \in \mathscr{R}$. We assume that $\Theta$ and $W$ are independent, and we define the rotational invariant 
random variable as follows 
\beam \label{eq.KPY.2.18} 
Z_A \stackrel{d}{=}\dfrac{\Theta}{W} \stackrel{d}{\sim} V_A \,,
\eeam
with $\stackrel{d}{=}$ to indicate equality in distribution. Let $B \in \mathcal{A}^*$. Then by relation \eqref{eq.KPY.2.17} 
we find that the $W$ is bounded away from the zero almost surely, hence the $1/W$ is bounded from above almost surely. Therefore, 
from this last and relation \eqref{eq.KPY.2.18} we obtain that $V_A \in \mathcal{A}^*$ as it follows by Lemma \ref{lem.KPY.4.3}. 
So we have $V \in  \mathcal{A}_A^*$. Due to arbitrary choice of $A \in \mathscr{R}$, we conclude $V \in  \mathcal{A}_{\mathscr{R}}^*$.
\eexam

\section{Local uniform asymptotic results} \label{sec.KPY.3}

Now, we present the first main result and its proof, that provides local uniform, with respect to time, asymptotic expressions for 
the first probability in \eqref{eq.KPY.1.2}, when the claim vectors follow a distribution from the class of multivariate subexponential 
distributions. By the term local uniformity we mean that the results are valid uniformly for $t \in \Lambda_T$, where 
$\Lambda_T : = [0,\,T]\cap \Lambda$, for all fixed $T \in \Lambda$. Let us notice that hereafter the distribution of 
the number of delayed claims $M$, has light tail, that means there exists some $\delta>0$, such that $\E\left[e^{\delta\,M} \right] < \infty$.

\subsection{Main result} \label{sec.KPY.3.1}

The following theorem represents the first main result of the paper.

\bth \label{th.KPY.3.1}
Let $A \in \mathscr{R}$ be some fixed set, let consider the discounted aggregate claims of relation \eqref{eq.KPY.1.1}. We suppose 
that Assumptions \ref{ass.KPY.1.1} and \ref{ass.KPY.1.2} are satisfied and there exists some $\delta > 0$, such that 
$\E\left[e^{\delta\,M} \right]< \infty$. Let some fixed $T \in \Lambda$.
\begin{enumerate}
\item[(i)]
If $F,\,G \in \mathcal{S}_A$ and $G(x\,A) \asymp F(x\,A)$, then it holds
\beam \label{eq.KPY.3.1} 
&&\PP \left({\bf D}_r(t) \in x\,A \right) \\[2mm] \notag
&&\sim \int_0^t  \PP\left( {\bf X} \in x\,e^{r\,s}\,A \right)\,\lambda(ds) + \E[M] \int_0^t \int_0^{t-s} \PP\left( {\bf Y} \in x\,e^{r\,(s+y)}\,A \right)\,H(dy)\,\lambda(ds) \,,
\eeam
uniformly for $t \in \Lambda_T$.
\item[(ii)]
If $F \in \mathcal{S}_A$ and $G(x\,A) =o\left[F(x\,A)\right]$, then it holds
\beam \label{eq.KPY.3.2} 
\PP \left({\bf D}_r(t) \in x\,A \right) \sim \int_0^t  \PP\left( {\bf X} \in x\,e^{r\,s}\,A \right)\,\lambda(ds) \,,
\eeam
uniformly for $t \in \Lambda_T$.
\end{enumerate}
\ethe 

\bre \label{rem.KPY.3.1}
If in Theorem \ref{th.KPY.3.1}(i) [or (ii)] we consider that $F,\,G \in \mathcal{S}_{\mathscr{R}}$ and $G(x\,A) \asymp F(x\,A)$, for 
any $A \in \mathscr{R}$ (or, $F \in \mathcal{S}_{\mathscr{R}}$ and  $G(x\,A) =o\left[F(x\,A)\right]$ for any $A \in \mathscr{R}$, respectively), 
then relation \eqref{eq.KPY.3.1} (or, relation \eqref{eq.KPY.3.2}, respectively), holds uniformly for $t \in \Lambda_T$, for any $A \in \mathscr{R}$. 
This can be easily verified, by the arbitrariness in the choice of $A \in \mathscr{R}$, in Lemmas \ref{lem.KPY.3.1} - \ref{lem.KPY.3.4} below. 
Theorem \ref{th.KPY.3.1} remains new result, even in the one-dimensional subcase, with $A=(1,\,\infty)$, since contributes the local uniformity to 
\cite[Lem. 3.5]{li:2023b}. In the case where the premiums have bounded densities, then the results of Theorem \ref{th.KPY.3.1} gives easily the finite-time ruin probability (over a ruin set connected with $A$) via some classical technics, see for example proof of \cite[Cor. 3.1]{konstantinides:passalidis:2024j} for more details. Hence in the one-dimensional subcase we can equip local uniformity in the 
\cite[Th. 2.1]{li:2023b}, as also we can relax the condition of constant premiums density into bounded density. We omit such a section in order to keep the paper compact. Furthermore, Theorem \ref{th.KPY.3.1}, when $M=0$ almost surely, is reduced to \cite[Th. 3.1]{konstantinides:liu:passalidis:2025}.     
\ere

\bre \label{rem.KPY.3.2}
Relations  \eqref{eq.KPY.3.1} and  \eqref{eq.KPY.3.2} establish the multivariate linear single big jump principle in the behavior of 
discounted aggregate claims. In \eqref{eq.KPY.3.2}, where the delayed claims distribution has negligible tail, in comparison with the 
tail distribution of the main claims, this principle becomes more obvious for the sequence of the main claim vectors, and this result, 
taking into account the \cite[Th. 3.1]{konstantinides:liu:passalidis:2025}, shows that when $G(x\,A)=o\left[F(x\,A)\right]$ the 
underestimation of the delayed claim does not bring any important asymptotic effect on the insurer. However, in \eqref{eq.KPY.3.1}, 
when $G(x\,A)\asymp F(x\,A)$, this relation shows that not only the multivariate linear single big jump principle of discounted aggregate claims comes from the main claims but 
also from the delayed ones, which depicts that both of them play crucial role for the insurer's solvency. Further, relation \eqref{eq.KPY.3.1} 
contains the $\E[M]$, that refers to the number of the delayed claims.
\ere

In the following remarks, we get more direct expressions for the estimations  \eqref{eq.KPY.3.1} and  \eqref{eq.KPY.3.2}, when the set $A$ is 
chosen between $A_1$ and $A_2$.

\bre \label{rem.KPY.3.3}
Let hold again the assumption of Example \ref{exam.KPY.2.1} for the vectors ${\bf X}$ and ${\bf Y}$, with marginal distributions from 
class $\mathcal{S}$, see Remark \ref{rem.KPY.2.2}. Then it holds
\beam \label{eq.KPY.3.3} 
\PP \left({\bf X} \in x\,A_1 \right) \sim \sum_{j=1}^d \PP\left( \dfrac{l_j}c\,X_j > x \right)\,, \qquad \PP \left({\bf Y} \in x\,A_1 \right) \sim \sum_{j=1}^d \PP\left( \dfrac{l_j}c\,Y_j > x \right)\,.
\eeam 
From \eqref{eq.KPY.3.1} and  \eqref{eq.KPY.3.3} is implied that 
\beam \label{eq.KPY.3.4} 
&&\PP \left({\bf D}_r(t) \in x\,A_1 \right) \sim  \\[2mm] \notag
&&\sum_{j=1}^d\left[ \int_0^t \PP\left( \dfrac{l_j}c\,X_j > x e^{r\,s} \right) \lambda(ds) +
\E[M] \int_0^t \int_0^{t-s} \PP\left(\dfrac{l_j}c Y_j > x e^{r\,(s+y)} \right) H(dy) \lambda(ds) \right] .
\eeam 
uniformly for $t \in \Lambda_T$.

Let us note that the uniformity remains through the application of \cite[Lem. 8]{wang:2011}. We recall by 
Example \ref{exam.KPY.2.1}, Remark \ref{rem.KPY.2.2} that $V_i' * V_j' \in \mathcal{S}$ and 
$e^{-r\,T} \leq e^{-r\,(s+y)} \leq e^{-r\,s} \leq 1$.

Similarly, by relation \eqref{eq.KPY.3.2} we have uniformly for $t \in \Lambda_T$:
\beam \label{eq.KPY.3.5} 
\PP \left({\bf D}_r(t) \in x\,A_1 \right) \sim \sum_{j=1}^d \int_0^t \PP\left( \dfrac{l_j}c\,X_j > x\,e^{r\,s} \right)\,\lambda(ds) \,.
\eeam

A set of the form $A_1$, and particularly with $c=1$, leads to discounted aggregate 
claims and ruin probabilities, that usually are denoted as $\psi_{sum}$. In risk models without 
delayed claims, we find many papers with an asymptotic expression similar to formula \eqref{eq.KPY.3.5}, 
with the only difference that in the integral of \eqref{eq.KPY.3.5} the probabilities are of the form
\beao
\PP\left(X_j > (l_1+\cdots +l_d)^{-1}\,x\,e^{r\,s} \right)=\PP\left(X_j > x\,e^{r\,s} \right)\,.
\eeao
Among the first papers in this direction we mention \cite{chen:wang:wang:2013}, and follow many other papers, see for example 
\cite{lu:qin:yuan:2025} for a model with delayed claims. The reason of this difference in the formula, is that, along that line, 
the claim vectors essentially satisfy a kind of multivariate non-linear single big jump principle (more precisely, satisfy the condition for 
$(\mathcal{D} \cap \mathcal{L})^{(2)}$, see in \cite{konstantinides:passalidis:2025c}), while here satisfy the multivariate linear single big jump principle.  

We should also mention that the first estimations for sets of the form $A_1$, in case of multivariate linear single big jump principle, 
appears in \cite{samorodnitsky:sun:2016}, that shows another reason of importance of the approach through the family $\mathscr{R}$, 
since such results do not exist even in $MRV$, the smallest class that satisfy the multivariate linear single big jump principle.
\ere
 
\bre \label{rem.KPY.3.4}
Let hold the assumptions of Example \ref{exam.KPY.2.2} for the vectors ${\bf X}$ and ${\bf Y}$, with marginal distributions from 
class $\mathcal{S}$. Then it holds
\beam \label{eq.KPY.3.6} 
\PP \left({\bf X} \in x\,A_2 \right) \sim \sum_{j=1}^d \PP\left( \dfrac{1}{c_j}\,X_j > x \right)\,, \qquad \PP \left({\bf Y} \in x\,A_2 \right) \sim \sum_{j=1}^d \PP\left( \dfrac{1}{c_j}\,Y_j > x \right)\,.
\eeam 
From \eqref{eq.KPY.3.1} and  \eqref{eq.KPY.3.6} is implied that 
\beam \label{eq.KPY.3.7} 
&&\PP \left({\bf D}_r(t) \in x\,A_2 \right) \sim \\[2mm] \notag
&&\sum_{j=1}^d \left[ \int_0^t \PP\left( \dfrac 1{c_j}\,X_j > x e^{r\,s} \right) \lambda(ds) +
\E[M] \int_0^t \int_0^{t-s} \PP\left(\dfrac 1{c_j}\,Y_j > x e^{r\,(s+y)} \right) H(dy) \lambda(ds) \right] .
\eeam 
uniformly for $t \in \Lambda_T$. The uniformity remains due to the relation
\beao
\PP \left(\bigvee_{j=1}^d w_j\,X_j > x \right) \sim \sum_{j=1}^d \PP(w_j\,X_j > x)\,,
\eeao
uniformly for $w_1,\,\ldots,\,w_d \in [a,\,b]^d$, with $0<a \leq b < \infty$, which follows by Bonferroni inequalities. The same argument for $Y_j$.

Similarly, by relation \eqref{eq.KPY.3.2} and \eqref{eq.KPY.3.6}  we have uniformly for $t \in \Lambda_T$ the asymptotic relation
\beam \label{eq.KPY.3.8} 
\PP \left({\bf D}_r(t) \in x\,A_2 \right) \sim \sum_{j=1}^d \int_0^t \PP\left( \dfrac 1{c_j}\,X_j > x\,e^{r\,s} \right)\,\lambda(ds) \,.
\eeam

In practical applications, a usual choice of $c_j$ is the $l_j$, because the 
$l_j\,x$ denotes the initial capital at the $j$-th line of business, with $l_1,\,\ldots,\,l_d >0$ and $l_1+ \cdots + l_d =1$. 
In this case, the form of \eqref{eq.KPY.3.8} seems to coincide with similar relations for ruin probabilities in both linear 
and non-linear big jump, usually denoted by $\psi_{min}$, see for example  \cite{chen:wang:wang:2013}.

An intuitive interpretation for the validity of the linear and non-linear case, is the nature  of the set $A_2$, in 
combination with the weak dependence structure of the components of ${\bf X}$, assumed in Example \ref{exam.KPY.2.2}. 
\ere

The following corollary provides more explicit forms of relations \eqref{eq.KPY.3.1} and \eqref{eq.KPY.3.2} in the case 
where the distribution class is restricted to $MRV$. The proof of this corollary follows easily by 
Theorem \ref{th.KPY.3.1} and the property \eqref{eq.KPY.2.6a} of $\mathcal{R}_{-\alpha}$ (recall that if 
$F \in MRV(\alpha,\,B,\,\mu)$, with $\mu(A)>0$, then $F_A \in \mathcal{R}_{-\alpha}$), therefore it is omitted.

\bco \label{cor.KPY.3.1}
Let $A \in \mathscr{R}$ be some fixed set. We consider the discounted aggregate claims from relation 
\eqref{eq.KPY.1.1} and some fixed $T \in \Lambda$.
\begin{enumerate}
\item[(i)]
Under the conditions of Theorem \ref{th.KPY.3.1} (i), with restrictions $F \in MRV(\alpha,\,B,\,\mu)$ and 
$G \in MRV(\alpha,\,Q,\,\nu)$, with $\mu(A)\wedge\nu(A)>0$, it holds
\beam \label{eq.KPY.3.b} \notag
&&\PP\left( {\bf D}_r(t) \in x\,A \right) \sim \mu(A)\,\bB(x)\,\int_0^t e^{-\alpha\,r\,s}\,\lambda(ds) \\[2mm]
&&+ \nu(A)\,\overline{Q}(x)\,\E[M]\,\int_0^t \int_0^{t-s} e^{-\alpha\,r\,(s+y)}\,H(dy)\,\lambda(ds)\,,
\eeam
uniformly for $t \in \Lambda_T$.
\item[(ii)]
Under the conditions of Theorem \ref{th.KPY.3.1} (ii), with restriction $F \in MRV(\alpha,\,B,\,\mu)$, with  $\mu(A)>0$, it holds
\beam \label{eq.KPY.3.c} 
\PP\left( {\bf D}_r(t) \in x\,A \right) \sim \mu(A)\,\bB(x)\,\int_0^t e^{-\alpha\,r\,s}\,\lambda(ds)\,,
\eeam
uniformly for $t \in \Lambda_T$.
\end{enumerate}
\eco
\bre \label{rem.KPY.3.a'}
We note that in Corollary \ref{cor.KPY.3.1} (i), for conservation of assumption $F(x\,A) \asymp G(x\,A)$, it is 
enough to hold $\bB(x) \asymp \overline{Q}(x)$, due to the fact that for any $A \in \mathscr{R}$ it holds 
$\mu(A),\, \nu(A) \in (0,\,\infty)$.
Correspondingly, in order to find if $G(x\,A)=o[F(x\,A)]$, it is enough to check if it holds $G(x\,A)=o[\bB(x)]$.
\ere
\subsection{Argumentation of Theorem \ref{th.KPY.3.1}} \label{sec.KPY.3.2}

Before presentation of the argumentation of the previous result, we need some preliminary lemmas. The first lemma is inspired by 
\cite[Lem. 3.2]{yang:li:2019}, with randomly weighted one-dimensional sums. In order to achieve the uniformity of Theorem \ref{th.KPY.3.1} 
we need its non-randomly weighted multivariate version. The following lemma generalizes the \cite[Lem. 5.1]{konstantinides:liu:passalidis:2025}. 
Hereafter, we denote
\beao
Z_A^{(i)} := \sup\left\{u\;:\; {\bf Z}^{(i)} \in u\,A \right\}\,,
\eeao
for any $i \in \bbn$, and for two random variables $Z_A^{(i)}$, $Z_A^{(j)}$ with $i\neq j$, we write $Z_A^{(i)} \leq_{st} Z_A^{(j)}$ if for any $x \in \bbr$ it holds
\beao
\PP(Z_A^{(i)} >x) \leq\PP(Z_A^{(j)} >x)\,.  
\eeao 

\ble \label{lem.KPY.3.1}
Let $A \in \mathscr{R}$ be some fixed set. We consider $n+m$ non-negative, independent random vectors ${\bf Z}^{(1)},\,\ldots,\,{\bf Z}^{(n+m)}$, with 
distributions $V_1,\,\ldots,\,V_{n+m}$, respectively. Let $0 <a \leq b < \infty$ be some constants and ${\bf c}_n:=(c_1,\,\ldots,\,c_n)^{\top} \in [a,\,b]^n$, 
while the $c_{n+1},\,\ldots,\,c_{n+m}$ take any values from the $[0,\,b]^m$, with the assumption that for any $n+1 \leq j \leq n+m$, there exists some 
$i_j \in \{1,\,\ldots,\,n\}$, such that $c_j \leq c_{i_j}$. If $V_1,\,\ldots,\,V_n \in \mathcal{L}_A$ and there exists $V \in \mathcal{S}_A$, with 
$V_i(x\,A) \asymp V(x\,A)$, for any $i=1,\,\ldots,\,n$, and $V_j(x\,A)=o\left[V(x\,A)\right]$, for any $j=n+1,\,\ldots,\,n+m$, then it holds
\beam \label{eq.KPY.3.9} 
\PP \left(\sum_{i=1}^{n+m} c_i\,{\bf Z}^{(i)} \in x\,A \right) \sim \sum_{i=1}^n \PP\left( c_i\,{\bf Z}^{(i)} \in x\,A \right) \,,
\eeam
uniformly for ${\bf c}_n \in [a,\,b]^n$, in the sense
\beao
\lim \sup_{{\bf c}_n \in [a,\,b]^n} \left| \dfrac{\PP \left(\sum_{i=1}^{n+m} c_i\,{\bf Z}^{(i)} \in x\,A \right)}{\sum_{i=1}^n \PP\left( c_i\,{\bf Z}^{(i)} \in x\,A \right)} -1\right|=0\,.
\eeao 
\ele

\pr~
For the lower bound for \eqref{eq.KPY.3.9} we obtain
\beao
\PP \left(\sum_{i=1}^{n+m} c_i\,{\bf Z}^{(i)} \in x\,A \right) \geq \PP \left(\sum_{i=1}^{n} c_i\,{\bf Z}^{(i)} \in x\,A \right) \sim \sum_{i=1}^{n} \PP \left(c_i\,{\bf Z}^{(i)} \in x\,A \right)\,,
\eeao
uniformly for ${\bf c}_n \in [a,\,b]^n$, where at the first step we used the increasing property of the set $A$, and that the summands 
are non-negative, while at the second step we take into account \cite[Lem. 5.1]{konstantinides:liu:passalidis:2025}.

It remains to show the upper bound for \eqref{eq.KPY.3.9}. At first, we observe that from $V_1,\,\ldots,\,V_n \in \mathcal{L}_A$, $V \in \mathcal{S}_A$ 
and $V_i(x\,A) \asymp V(x\,A)$, for any $i=1,\,\ldots,\,n$, is implied by \cite[Prop. 3.3]{konstantinides:passalidis:2025h} that $V_1,\,\ldots,\,V_n \in \mathcal{S}_A$.

Further, since $V_j(x\,A) = o\left[V(x\,A)\right]$, for $j=n+1,\,\ldots,\,n+m$, hence for any $\vep \in (0,\,1)$ we can find some $x_0 =x_0(\vep) >0$, 
such that it holds
\beam \label{eq.KPY.3.10} 
V_j(x\,A) \leq \vep\,V_{i_j}(x\,A)\,,
\eeam
for any $x \geq x_0$, recall that $V_{i_j}(x\,A) \asymp V(x\,A)$ for all $i_j \in \{1,\,\ldots,\,n\}$. We now construct the new random vectors ${\bf Z}^{(n+1)*},\,\ldots,\,{\bf Z}^{(n+m)*}$ 
which are independent each other and independent of the rest sources of randomness, with distributions
\beam \label{eq.KPY.3.11} 
V_j^*(x\,A) =\max\{\vep\,V_{i_j}(x\,A),\,V_j(x\,A)\} \geq V_j(x\,A)\,,
\eeam
respectively, for any $x\geq 0$. We see from relations \eqref{eq.KPY.3.10}, \eqref{eq.KPY.3.11} that the ${\bf Z}^{(n+1)*},\,\ldots,\,{\bf Z}^{(n+m)*}$ 
are also non-negative. Let $Z_A^{(i)*}:=\sup\{u\;:\;{\bf Z}^{(i)*} \in u\,A\}$, 
for any $j=n+1,\,\ldots,\,n+m$. From relation \eqref{eq.KPY.3.11} we obtain that 
$Z_A^{(i)} \leq_{st} Z_A^{(j)*}$. From relation \eqref{eq.KPY.3.11} is implied
\beam \label{eq.KPY.3.12} 
V_j^*(x\,A) \sim \vep\,V_{i_j}(x\,A)\,,
\eeam 
and consequently $V_j^* \in \mathcal{S}_A$, by the closure property of $\mathcal{S}_A$ with respect to strong tail equivalence, see 
\cite[Prop. 4.12(a)]{samorodnitsky:sun:2016}. Hence, from \eqref{eq.KPY.3.12} we find $V_j^*(x\,A) \asymp V(x\,A)$. Therefore, we find
\beam \label{eq.KPY.3.13} \notag
&&\PP \left(\sum_{i=1}^{n+m} c_i\,{\bf Z}^{(i)} \in x\,A \right) \leq \PP \left(\sum_{i=1}^{n+m} c_i\, Z_A^{(i)} > x\right) \\[2mm] \notag
&&\leq \PP \left(\sum_{i=1}^{n} c_i\, Z_A^{(i)} + \sum_{j=n+1}^{n+m} c_{i_j}\, Z_A^{(j)*} >x\right)\\[2mm] \notag
&&\sim \sum_{i=1}^{n}\PP \left( c_i\, Z_A^{(i)} >x \right) + \sum_{j=n+1}^{n+m}\PP \left( c_{i_j}\, Z_A^{(j)*} >x\right)\\[2mm] 
&&= \sum_{i=1}^{n}\PP \left( c_i\, {\bf Z}^{(i)} \in x\,A\right) + \sum_{j=n+1}^{n+m}\PP \left( c_{i_j}\, {\bf Z}^{(j)*} \in x\,A\right)\,,
\eeam
uniformly for ${\bf c}_n \in [a,\,b]^n$, where at the first step we used \cite[Prop. 2.4]{konstantinides:passalidis:2024g}, at the 
second step we apply the inequalities $Z_A^{(j)} \leq_{st} Z_A^{(j)*}$, $c_j \leq c_{i_j}$ and at the third step we use \cite[Lem. 1]{tang:yuan:2014} 
since the $Z_A^{(i)}$ are subexponential with weak equivalent tails and all the coefficients belong to the interval $[a,\,b]$. By \eqref{eq.KPY.3.12} we obtain
\beam \label{eq.KPY.3.14}
\sum_{j=n+1}^{n+m}\PP \left( c_{i_j}\,{\bf Z}^{(j)*} \in x\,A\right) \sim \vep\,\sum_{j=n+1}^{n+m}\PP \left( c_{i_j}\,{\bf Z}^{(i_j)} \in x\,A\right)\,.
\eeam 
From relations \eqref{eq.KPY.3.13} and \eqref{eq.KPY.3.14} and the arbitrary choice of $\vep >0$, we have the upper bound. 
~\halmos

\bre \label{rem.KPY.3.5}
Due to uniformity of relation \eqref{eq.KPY.3.9}, if in Lemma \ref{lem.KPY.3.1} hold the same conditions, but the weighted coefficients 
become random variables $\Theta_i$, with $i=1,\,\ldots,\,n+m$, with the $(\Theta_1,\,\ldots,\,\Theta_n)^{\top} \in [a,\,b]^n$ almost surely, 
the $(\Theta_{n+1},\,\ldots,\,\Theta_{n+m})^{\top} \in [0,\,b]^m$ almost surely, and for any $j=n+1,\,\ldots,\,n+m$ there exists some 
$i_j=1,\,\ldots,\,n$, such that $\Theta_j \leq_{st} \Theta_{i_j}$, then by dominated convergence theorem we obtain
\beam \label{eq.KPY.3.15}
\PP \left( \sum_{i=1}^{n+m} \Theta_{i}\,{\bf Z}^{(i)} \in x\,A\right) \sim \sum_{i=1}^{n}\PP \left(\Theta_{i}\,{\bf Z}^{(i)} \in x\,A\right)\,.
\eeam 
Relation \eqref{eq.KPY.3.15} demonstrates the multivariate linear single  big jump principle for the randomly weighted sum $\sum_{i=1}^{n} \Theta_{i}\,{\bf Z}^{(i)}$.
\ere

The following lemma represent a multivariate version of the Kesten inequality, and in fact is a reformulation of \cite[Prop. 4.12(c)]{samorodnitsky:sun:2016}.

\ble \label{lem.KPY.3.2}
Let $A \in \mathscr{R}$ be some fixed set, and ${\bf Z}^{(1)},\,\ldots,\,{\bf Z}^{(n)}$ be i.i.d. copies of non-negative random vectors, 
with ${\bf Z} \stackrel{d}{\sim} V \in \mathcal{S}_A$. Then for any $\vep > 0$, there exists some constant $C_{\vep}>0$, such that it holds
\beam \label{eq.KPY.3.16}
\PP \left( {\bf Z}^{(1)} +\cdots +  {\bf Z}^{(n)} \in x\,A\right) \leq C_{\vep}\,(1+\vep)^n \,\PP \left({\bf Z}\in x\,A\right)\,,
\eeam 
for any $x>0$ and any $n \in \bbn$.
\ele

The next auxiliary assertion follows immediately by proof of \cite[Th. 3.1]{konstantinides:liu:passalidis:2025}, see there relation $(33)$ 
and \cite[Eq. (3.3)]{hao:tang:2008} for the one-dimensional case.

\ble \label{lem.KPY.3.3}
Let $A \in \mathscr{R}$ be some fixed set, and $\{{\bf Z}^{(i)}\,,\;i \in \bbn\}$ be i.i.d. copies of non-negative random vectors, 
with ${\bf Z} \stackrel{d}{\sim} V \in \mathcal{S}_A$. We assume that the arrival times $\{\tau_i\,,\;i \in \bbn\}$ constitutes a renewal 
counting process $\{N(t)\,,\;t\geq 0\}$ with finite renewal function $\lambda(t)$.  We suppose that $\{{\bf Z}^{(i)}\,,\;i \in \bbn\}$ 
and  $\{N(t)\,,\;t\geq 0\}$ are independent, $r\geq 0$ is the constant interest rate, and we consider some fixed $T \in \Lambda$. Then it holds
\beam \label{eq.KPY.3.17}
\lim_{N \to \infty} \sup_{t \in \Lambda_T} \sup_{x\geq 0} \dfrac{\sum_{n=N+1}^{\infty} \sum_{i=1}^{n} \PP \left( {\bf Z}^{(i)}\,e^{-r\,\tau_i} \in x\,A\,,\;N(t)=n\right)}{\int_0^t \PP \left( {\bf Z} \in x\,e^{r\,s}\,A\right)\,\lambda(ds)}=0\,.
\eeam
\ele

The last auxiliary result, is related with the closure properties of distribution class $\mathcal{S}_A$ with respect to convolution.

\ble \label{lem.KPY.3.4}
Let $A \in \mathscr{R}$ be some fixed set, and $V_1,\,V_2$ be distributions with support in the non-negative orthant. We suppose that 
$V_1 \in \mathcal{S}_A$ and $V_2(x\,A)=o\left[V_1(x\,A)\right]$. Then it holds
\beam \label{eq.KPY.3.18}
V_1*V_2(x\,A) \sim V_1(x\,A)\,,
\eeam 
and further $V_1*V_2 \in \mathcal{S}_A$.
\ele

\pr~
The proof could be an easy consequence of \cite[Cor. 3.1]{konstantinides:passalidis:2024g}, and some one-dimensional closure properties 
with respect to convolution of subexponential distributions. However, for sake of completeness, we provide an alternative proof.

Let ${\bf Z}^{(1)}\,,\;{\bf Z}^{(2)}$ be non-negative random vectors with distributions $V_1$ and $V_2$, respectively. Then
\beao
V_1*V_2(x\,A)=\PP \left( {\bf Z}^{(1)}+ {\bf Z}^{(2)} \in x\,A\right) \leq\PP \left( Z_A^{(1)}+ Z_A^{(2)}> x\right) \sim\PP \left( Z_A^{(1)}> x\right)=V_1(x\,A)\,,
\eeao
where at the second step we used \cite[Prop. 2.4]{konstantinides:passalidis:2024g}, and at the third step we take advantage from 
\cite[Cor. 3.18]{foss:korshunov:zachary:2013}, since $\bV_A^{(2)}(x)= o\left[\bV_A^{(1)}(x)\right]$, see also 
\cite[Prop. 3.13(v)]{leipus:siaulys:konstantinides:2023}. The last relation provides the upper bound for \eqref{eq.KPY.3.18}.

For the lower bound, since the set $xA$ is increasing and the ${\bf Z}^{(1)}\,,\;{\bf Z}^{(2)}$ are non-negative, via Bonferroni inequality we obtain
\beao 
&&V_1*V_2(x\,A) \geq \PP \left(\bigcup_{i=1}^2 \left\{ {\bf Z}^{(i)} \in x\,A \right\}\right) \\[2mm]
&&\geq \PP \left( {\bf Z}^{(1)} \in x\,A\right) + \PP \left( {\bf Z}^{(2)} \in x\,A\right) - \PP \left( {\bf Z}^{(1)} \in x\,A\right)\,\PP \left( {\bf Z}^{(2)} \in x\,A\right) \\[2mm]
&&\sim V_1(x\,A) + V_2(x\,A) = [1+o(1)]\,V_1(x\,A)\,.
\eeao
Which completes the proof of relation \eqref{eq.KPY.3.18}.

From  \eqref{eq.KPY.3.18}, keeping in mind that $V_1 \in \mathcal{S}_A$, from the closure property of $\mathcal{S}_A$ with respect to 
strong tail equivalence, see \cite[Prop. 4.12(a)]{samorodnitsky:sun:2016}, we obtain  $V_1*V_2 \in \mathcal{S}_A$.    
~\halmos

\noindent{\bf Proof of Theorem \ref{th.KPY.3.1}.}~
We show explicitly the assertion $(i)$, since the assertion $(ii)$ is implied by similar way, and we indicate only the points where 
the assertion $(ii)$ uses different arguments.

For sake of simplicity, we denote by $h(x;\,t)$ the right member of relation \eqref{eq.KPY.3.1} and by $h^*(x;\,t)$ the right 
member of relation \eqref{eq.KPY.3.2}.

\begin{enumerate}
\item[(i)]
Let consider some $N \in \bbn$. Then by total probability theorem and Assumption \ref{ass.KPY.1.1} we obtain for all $t\in\Lambda_T$
\beam \label{eq.KPY.3.19} \notag
&&\PP\left({\bf D}_r(t) \in x\,A \right) = \PP\left(\sum_{i=1}^{N(t)} \left[{\bf X}^{(i)}\,e^{-r\,\tau_i} + \sum_{j=1}^{M_i} {\bf Y}^{(i,j)}\,e^{-r\,(\tau_i + D_{ij})}\,{\bf 1}_{\{\tau_i + D_{ij} \leq t\}}  \right] \in x\,A \right)\\[2mm] \notag
&&=\left(\sum_{n=1}^{N} +\sum_{n=N+1}^{\infty}\right)\\[2mm] \notag
&& \PP\left(\sum_{i=1}^{n} {\bf X}^{(i)}\,e^{-r\,\tau_i} + \sum_{i=1}^{n} \sum_{j=1}^{M_i} {\bf Y}^{(i,j)}\,e^{-r\,(\tau_i + D_{ij})}\,{\bf 1}_{\{\tau_i + D_{ij} \leq t\}} \in x\,A\,,\;N(t)=n \right)\\[2mm]
&&=:I_1(x,\,N,\,r,\,t) + I_2(x,\,N,\,r,\,t)\,.
\eeam
For $I_2(x,\,N,\,r,\,t)$, again from Assumption \ref{ass.KPY.1.1} we find for all $t\in\Lambda_T$
\beam \label{eq.KPY.3.20} \notag
&&I_2(x,\,N,\,r,\,t)=\sum_{n=N+1}^{\infty} \sum_{m_1=0}^{\infty} \cdots \sum_{m_n=0}^{\infty} \\[2mm] \notag
&&\PP\left(\sum_{i=1}^{n} {\bf X}^{(i)}\,e^{-r\,\tau_i} + \sum_{i=1}^{n} \sum_{j=1}^{m_i} {\bf Y}^{(i,j)}\,e^{-r\,(\tau_i + D_{ij})}\,{\bf 1}_{\{\tau_i + D_{ij} \leq t\}} \in x\,A\,,\;N(t)=n \right)\\[2mm]
&&\times \PP\left(M_1=m_1,\,\ldots,\,M_n=m_n \right)\\[2mm] \notag
&&\leq \sum_{n=N+1}^{\infty} \sum_{m_1=0}^{\infty} \cdots \sum_{m_n=0}^{\infty} \PP\left(\sum_{i=1}^{n} {\bf X}^{(i)}\,e^{-r\,\tau_1} + \sum_{i=1}^{n} \sum_{j=1}^{m_i} {\bf Y}^{(i,j)}\,e^{-r\,\tau_1} \in x\,A\,,\;N(t)=n \right)\\[2mm] \notag
&&\times \PP\left(M_1=m_1,\,\ldots,\,M_n=m_n \right)\\[2mm] \notag
&&\leq \sum_{n=N+1}^{\infty} \sum_{m_1=0}^{\infty} \cdots \sum_{m_n=0}^{\infty} \PP\left(\sum_{i=1}^{n+m_1+\cdots +m_n} \left[{\bf X}^{(i)}+{\bf Y}^{(i)*}\right]\,e^{-r\,\tau_1} \in x\,A\,,\;N(t)=n \right)\\[2mm] \notag
&&\times \PP\left(M_1=m_1,\,\ldots,\,M_n=m_n \right) \,.
\eeam
where $\left\{{\bf Y}^{(i)*}\,,\;i \in \bbn \right\}$ are i.i.d.  random vectors with common distribution $G$, and are independent of the rest 
sources of randomness, so the sequence $\left\{{\bf X}^{(i)}+{\bf Y}^{(i)*}\,,\;i \in \bbn \right\}$ has i.i.d.  terms with common distribution $F*G$. 

We shall show that $F*G \in \mathcal{S}_A$. If, for all $i\in\mathbb{N}$, $Y_A^{(i)*}:= \sup \left\{u\;:\;{\bf Y}^{(i)*} \in u\,A \right\} \stackrel{d}{\sim} G_A$, then, 
because of $F,\,G \in \mathcal{S}_A \subsetneq \mathcal{L}_A$, due to \cite[Th. 3.4]{konstantinides:passalidis:2024g}, it is enough to show 
that $F_A*G_A \in \mathcal{S}$. Indeed, this is true, since $\bF_A(x) \asymp \bG_A(x)$, through \cite[Lem. 3.2]{tang:tsitsiashvili:2003}, 
see also \cite[Prop. 3.13(iv)]{leipus:siaulys:konstantinides:2023}. Hence, $F*G \in \mathcal{S}_A$, and further from 
\cite[Cor. 3.1]{konstantinides:passalidis:2024g} it holds
\beam \label{eq.KPY.3.21} 
F*G(x\,A)\sim F(x\,A)+ G(x\,A)\,.
\eeam
Next, since $G(x\,A) \asymp F(x\,A)$, there exists some finite $\omega >0$, such that 
\beao
\limsup \dfrac{G(x\,A)}{F(x\,A)}< \omega\,.
\eeao
Therefore, from relation \eqref{eq.KPY.3.20}, for arbitrary $\vep >0$, from Lemma \ref{lem.KPY.3.2}, for the ${\bf X}^{(i)}+{\bf Y}^{(i)*}$, 
we can find some $C_{\vep}>0$, such that for all $t\in\Lambda_T$ it holds
\beam \label{eq.KPY.3.22}  \notag
&&I_2(x,\,N,\,r,\,t) \\[2mm] \notag
&&\leq \sum_{n=N+1}^{\infty} \sum_{m_1=0}^{\infty} \cdots \sum_{m_n=0}^{\infty} \int_0^t \PP\left(\sum_{i=1}^{n+m_1+\cdots +m_n} \left[{\bf X}^{(i)}+{\bf Y}^{(i)*}\right] \in x\,e^{r\,s}\,A\right)\,\\[2mm] \notag
&&\times \PP\left(N(t-s)=n-1 \right)\,\PP(\tau_1 \in ds)\PP\left(M_1=m_1,\,\ldots,\,M_n=m_n \right) \\[2mm] \notag
&&\leq C_{\vep} \sum_{n=N+1}^{\infty} \sum_{m_1=0}^{\infty} \cdots \sum_{m_n=0}^{\infty} \int_0^t (1+\vep)^{n+m_1+\cdots +m_n} F*G\left(x\,e^{r\,s}\,A\right)\\[2mm]
&&\times \PP\left(N(t-s)=n-1 \right)\,\PP(\tau_1 \in ds)\,\PP\left(M_1=m_1,\,\ldots,\,M_n=m_n \right)\\[2mm] \notag
&&\lesssim C_{\vep} (1+\omega)\,\sum_{n=N+1}^{\infty} (1+\vep)^{n} \sum_{m_1=0}^{\infty} (1+\vep)^{m_1} \PP\left(M_1=m_1\right) \cdots \\[2mm] \notag
&&\sum_{m_n=0}^{\infty} (1+\vep)^{m_n} \PP\left(M_n=m_n\right)\,\int_0^t \PP\left({\bf X}\,e^{-r\,s}\in x\,A\right) \,\PP\left(N(t-s)=n-1 \right)\,\lambda(ds)\\[2mm] \notag
&&= C_{\vep} (1+\omega)\,\sum_{n=N+1}^{\infty} \E\left[\left(1+\vep\right)^M\right]^{n} (1+\vep)^{n} \\[2mm] \notag
&&\times \int_0^t \PP\left({\bf X}\,e^{-r\,s}\in x\,A \right)\, \PP\left(N(t-s)=n-1 \right)\,\lambda(ds)\leq C_{\vep} (1+\omega)\,\\[2mm] \notag
&&\times \E_{N(t)}\left[\left\{\left(1+\vep\right)\,\E_{M}\left[\left(1+\vep\right)^M\right]\right\}^{N(t)+1}\,{\bf 1}_{\{N(t) \geq N\}} \right] \int_0^t \PP\left({\bf X}\in x\,e^{r\,s}\,A \right)\,\lambda(ds) \,,
\eeam
where at the third step we used relation \eqref{eq.KPY.3.21} in combination with the fact that $1 \leq e^{r\,s} \leq e^{r\,t} < \infty$.

Because  $\{N(t)\,,\;t\geq 0\}$ represent a renewal process with finite renewal function, for any $t \in \Lambda_T$, there exists 
some $\delta^*>0$, such that $\E\left[e^{\delta^*\,N(t)} \right]< \infty$, see \cite{stein:1946}, hence for any $\delta'>0$, we 
can find some $\vep>0$, small enough, and some $N\in \bbn$, large enough, such that it holds
\beam \label{eq.KPY.3.23} 
\E_{N(t)}\left[\left\{\left(1+\vep\right)\,\E_{M}\left[\left(1+\vep\right)^M\right]\right\}^{N(t)+1}\,{\bf 1}_{\{N(t) \geq N\}} \right] < \delta'\,.
\eeam
So, from relations \eqref{eq.KPY.3.22} and \eqref{eq.KPY.3.23}, we obtain for all $t\in\Lambda_T$
\beam \label{eq.KPY.3.24} 
I_2(x,\,N,\,r,\,t) \leq \delta'\,C_{\vep} (1+\omega)\,h(x;\,t) \,.
\eeam  

Now we deal with $I_1(x,\,N,\,r,\,t)$. For sake of simplicity, for a sequence $\{a_i\,,\;i \in \mathcal{K}\}$, where $\mathcal{K} \subseteq \bbn_0$, we denote 
\beao
\check{a}_\mathcal{K}:=\bigvee_{i\in \mathcal{K}} a_i\,.
\eeao
So, we split $I_1(x,\,N,\,r,\,t)$ into $I_{11}(x,\,N,\,r,\,t)$ and $I_{12}(x,\,N,\,r,\,t)$, according to if for any $n=1,\,\ldots,\,N$ it holds
\beao
\check{M}_{\{1,\,\ldots,\,n\}} \leq N\,,
\eeao
or does not hold, respectively. For the $I_{11}(x,\,N,\,r,\,t)$ we obtain for all $t\in\Lambda_T$
\beam \label{eq.KPY.3.25}  \notag
&&I_{11}(x,\,N,\,r,\,t) =\sum_{n=1}^{N} \\[2mm] \notag
&&\PP\left(\sum_{i=1}^{n} {\bf X}^{(i)} e^{-r\,\tau_i} + \sum_{i=1}^{n} \sum_{j=1}^{M_i} {\bf Y}^{(i,j)} e^{-r (\tau_i + D_{ij})} {\bf 1}_{\{\tau_i + D_{ij} \leq t\}} \in x A ,\,N(t)=n ,\,\check{M}_{\{1, \ldots, n\}} \leq N \right) \\[2mm] \notag
&&=\sum_{n=1}^{N} \sum_{m_1=0}^{N} \cdots \sum_{m_n=0}^{N} \\[2mm] \notag
&&\PP\left(\sum_{i=1}^{n} {\bf X}^{(i)}\,e^{-r\,\tau_i} + \sum_{i=1}^{n} \sum_{j=1}^{m_i} {\bf Y}^{(i,j)}\,e^{-r\,(\tau_i + D_{ij})}\,{\bf 1}_{\{\tau_i + D_{ij} \leq t\}} \in x\,A\,,\;N(t)=n \right) \, \\[2mm] \notag
&&\times \PP\left(M_1=m_1,\,\ldots,\,M_n=m_n \right) =\sum_{n=1}^{N} \sum_{m_1=0}^{N} \cdots \sum_{m_n=0}^{N} \iint_{\{0<s_1\leq \ldots \leq s_n \leq t < s_{n+1}\}}\\[2mm] \notag
&&\iint_{\{u_{ij}>0,\,s_i+u_{ij} \leq t\}}\,\PP\left(\sum_{i=1}^{n} {\bf X}^{(i)}\,e^{-r\,s_i} + \sum_{i=1}^{n} \sum_{j=1}^{m_i} {\bf Y}^{(i,j)}\,e^{-r\,(s_i + u_{ij})} \in x\,A\right) \,{\bf H}(d{\bf u}) \\[2mm] \notag
&&\times \PP\left(\tau_1 \in ds_1,\,\ldots,\,\tau_{n+1} \in ds_{n+1} \right)\,\PP\left(M_1=m_1,\,\ldots,\,M_n=m_n \right) \\[2mm]
&&\sim \sum_{n=1}^{N} \sum_{m_1=0}^{N} \cdots \sum_{m_n=0}^{N} \iint_{\{0<s_1\leq \ldots \leq s_n \leq t < s_{n+1}\}}\,\iint_{\{u_{ij}>0,\,s_i+u_{ij} \leq t\}}\\[2mm] \notag
&&\left[\sum_{i=1}^{n}\PP\left( {\bf X}^{(i)}\,e^{-r\,s_i} \in x\,A\right) + \sum_{i=1}^{n} \sum_{j=1}^{m_i}\PP\left( {\bf Y}^{(i,j)}\,e^{-r\,(s_i + u_{ij})} \in x\,A\right) \right]\,{\bf H}(d{\bf u}) \\[2mm] \notag
&&\times \PP\left(\tau_1 \in ds_1,\,\ldots,\,\tau_{n+1} \in ds_{n+1} \right)\,\PP\left(M_1=m_1,\,\ldots,\,M_n=m_n \right)\\[2mm] \notag
&&= \sum_{n=1}^{N} \sum_{i=1}^{n} \PP\left( {\bf X}^{(i)}\,e^{-r\,\tau_i} \in x\,A\,,\;N(t) = n\right) \,\PP\left( \check{M}_{\{1,\,\ldots,\,n\}} \leq N\right)+ \sum_{n=1}^{N} \sum_{i=1}^{n} \\[2mm] \notag
&&\PP\left( {\bf Y}^{(i)}\,e^{-r\,(\tau_i + D)}\,{\bf 1}_{\{\tau_i + D \leq t\}} \in x\,A\,,\;N(t)=n\right) \sum_{m_i=0}^N \,\PP\left(M_i=m_i\right)\,\,\PP\left( \check{M}_{\{1,\,\ldots,\,n\}\setminus \{i\}} \leq N\right) \,,
\eeam
where by ${\bf H}$ we denote the joint distribution of the delays 
$D_{11},\,\ldots,\,D_{1,m_1},\,\ldots,\,D_{n_1},$ $\,\ldots,\,D_{n,m_n}$, at the fourth step we used Lemma \ref{lem.KPY.3.1}, with 
the $V_1,\,\ldots,\,V_n$ there, observing that for any $i \in \bbn$ it holds 
\beao
e^{-r\,s_i} \geq e^{-r\,(s_i + u_{ij})}\,,
\eeao 
with $e^{-r\,t} \leq e^{-r\,s_i}\leq 1$, and $0 \leq e^{-r\,(s_i + u_{ij})}\leq 1$, and at the last step, $D$ represent a general 
random variable from the $\{D_{ij}\,,\;i,\,j \in \bbn\}$.

From relation \eqref{eq.KPY.3.25}, for the upper bound of $I_{11}(x,\,N,\,r,\,t)$ we find that for all $t\in\Lambda_T$ it holds
\beam \label{eq.KPY.3.26} \notag
&&I_{11}(x,\,N,\,r,\,t) \lesssim \sum_{n=1}^{\infty} \sum_{i=1}^{n} \PP\left( {\bf X}^{(i)}\,e^{-r\,\tau_i} \in x\,A\,,\;N(t) = n\right) \\[2mm] 
&&+ \E[M]\,\sum_{n=1}^{\infty} \sum_{i=1}^{n} \PP\left( {\bf Y}^{(i)}\,e^{-r\,(\tau_i + D_{ij})}\,{\bf 1}_{\{\tau_i + D_{ij} \leq t\}} \in x\,A\,,\;N(t)=n\right) \\[2mm] \notag
&&=\sum_{i=1}^{\infty} \PP\left({\bf X}^{(i)}\,e^{-r\,\tau_i} \in x\,A\,,\;\tau_i \leq t\right)+ \E[M]\,\sum_{i=1}^{\infty} \PP\left( {\bf Y}^{(i)}\,e^{-r\,(\tau_i + D)}\,{\bf 1}_{\{\tau_i + D \leq t\}} \in x\,A\,,\;\tau_i \leq t\right) \\[2mm] \notag
&&=\int_0^t \PP\left({\bf X}\,e^{-r\,s} \in x\,A\right)\,\lambda(ds)\\[2mm] \notag
&&+ \E[M]\,\int_0^t \PP\left( {\bf Y}\,e^{-r\,(s + D)}\,{\bf 1}_{\{D \leq t -s\}} \in x\,A\right)\,\lambda(ds)=h(x;\,t)\,,
\eeam
where at the second step we change the order of summations. From the other hand side, since 
there exists a $\delta>0$, such that $\E\left[e^{\delta\,M} \right] < \infty$, we obtain $\E[M]< \infty$, hence for any $\delta' >0$, 
we can find large enough $N_1 \in \bbn$, such that for any $n=1,\,\ldots,\,N_1$ and $i=1,\,\ldots,\,n$ it holds
\beam \label{eq.KPY.3.27}
&&\PP\left( \check{M}_{\{1,\,\ldots,\,n\}\setminus \{i\}} > N_1 \right) \leq \PP\left( \check{M}_{\{1,\,\ldots,\,n\}} > N_1\right) \\[2mm] \notag
&&=\PP\left( \bigcup_{i=1}^n \{M_i > N_1\}\right) \leq N_1\,\PP\left( M> N_1 \right)\leq \sum_{m=N_1+1}^{\infty} m\,\PP(M=m) < \delta'\,.
\eeam
Also, from Lemma \ref{lem.KPY.3.3}, for the same $\delta'>0$, we can find some large enough $N_2 \in \bbn$, such that for 
all $t \in \Lambda_T$ and all $x>0$, it holds
\beam \label{eq.KPY.3.28}
\sum_{n=N_2+1}^{\infty} \sum_{i=1}^n \PP\left({\bf X}^{(i)}\,e^{-r\,\tau_i} \in x\,A\,,\;N(t) =n\right) \leq \delta' h^*(x;\,t) \leq \delta' h(x;\,t)\,.
\eeam
and 
\beam \label{eq.KPY.3.29} \notag
&&\sum_{n=N_2+1}^{\infty} \sum_{i=1}^n\PP\left( {\bf Y}^{(i)}\,e^{-r\,(\tau_i + D)}\,{\bf 1}_{\{\tau_i + D \leq t\}} \in x\,A\,,\;N(t) =n\right) \\[2mm] 
&&\leq \sum_{n=N_2+1}^{\infty} \sum_{i=1}^n \PP\left( {\bf Y}^{(i)}\,e^{-r\,\tau_i} \in x\,A\,,\;N(t) =n\right) \\[2mm] \notag
&&\leq \delta' \int_0^t \PP\left( {\bf Y} \in x\,e^{r\,s}\,A\right)\,\lambda(ds) \leq \delta'\,\omega\,\int_0^t \PP\left( {\bf X} \in x\,e^{r\,s}\,A\right)\,\lambda(ds)\leq \delta'\,\omega\, h(x;\,t)\,.
\eeam

From relations \eqref{eq.KPY.3.27} - \eqref{eq.KPY.3.29}, in combination with \eqref{eq.KPY.3.25}, and denoting $N = N_1\vee N_2$, we obtain
\beam \label{eq.KPY.3.30}\notag
&&I_{11}(x,\,N,\,r,\,t) \gtrsim \left[\sum_{n=1}^{\infty} \sum_{i=1}^{\infty}\PP\left( {\bf X}^{(i)}\,e^{-r\,\tau_i} \in x\,A\,,\;N(t) =n\right)- \delta'\,h(x;\,t) \right] (1-\delta') \\[2mm] 
&& + \left[\sum_{n=1}^{\infty} \sum_{i=1}^{n} \PP\left( {\bf Y}^{(i)}\,e^{-r\,(\tau_i + D)}\,{\bf 1}_{\{\tau_i + D \leq t\}} \in x\,A\,,\;N(t) =n\right) - \delta'\,\omega\,h(x;\,t)\right]\,\\[2mm] \notag
&&\times \left(\E[M]-\delta' \right)\,(1-\delta') \geq (1-\delta') \,[1-\delta'\,(1+\omega+\omega\,\E[M])]\,h(x;\,t) =:(1-\delta')\,(1-\delta'\,C_1)\,h(x;\,t)\,,
\eeam
for all $t \in \Lambda_T$, where 
\beao
C_1=1+\omega+\omega\,\E[M]>0\,,
\eeao 
is a constant. Hence, from relations \eqref{eq.KPY.3.26} and \eqref{eq.KPY.3.30} we obtain
\beam \label{eq.KPY.3.31}
(1-\delta')\,(1-\delta'\,C_1)\,h(x;\,t) \lesssim I_{11}(x,\,N,\,r,\,t) \lesssim h(x;\,t)\,,
\eeam
uniformly for $t \in \Lambda_T$. 

Now, we handle $I_{12}(x,\,N,\,r,\,t)$, for which we apply similar approach with that for $I_{2}(x,\,N,\,r,\,t)$. For any $t \in \Lambda_T$, we obtain
\beam \label{eq.KPY.3.32}\notag
&&I_{12}(x,\,N,\,r,\,t)= \sum_{n=1}^{N}\PP\Bigg(\sum_{i=1}^{n} {\bf X}^{(i)}\,e^{-r\,\tau_i}  + \sum_{i=1}^{n} \sum_{j=1}^{M_i} {\bf Y}^{(i,j)}\,e^{-r\,(\tau_i + D_{ij})}\,{\bf 1}_{\{\tau_i + D_{ij} \leq t\}} \in x\,A\,,  \\[2mm] \notag
&&\;N(t)=n,\;\check{M}_{\{1,\,\ldots,\,n\}} > N \Bigg) \leq \sum_{n=1}^{\infty} \sum_{k=1}^n \PP\Bigg(\sum_{i=1}^{n} {\bf X}^{(i)}\,e^{-r\,\tau_i}  + \sum_{i=1}^{n} \sum_{j=1}^{M_i} {\bf Y}^{(i,j)}\,e^{-r\,(\tau_i + D_{ij})}\,\\[2mm] \notag
&&\times {\bf 1}_{\{\tau_i + D_{ij} \leq t\}} \in x\,A\,,\;N(t)=n\,,\;M_k>N \Bigg) \leq \sum_{n=1}^{\infty} \sum_{k=1}^n \sum_{m_k=n+1}^{\infty} \sum_{m_1=0}^{\infty} \cdots \sum_{m_{k-1}=0}^{\infty}\\[2mm] \notag
&&  \sum_{m_{k+1}=0}^{\infty} \cdots \sum_{m_N=0}^{\infty} \int_0^t  \PP\left(\sum_{i=1}^{n+m_1+\cdots +m_n} \left[{\bf X}^{(i)}+{\bf Y}^{(i)*}\right] \in x\,e^{r\,s}\,A\right)\,\PP\left(N(t-s)=n-1 \right)\,\\[2mm] \notag
&&\times \PP(\tau_1 \in ds)\,\PP\left(M_1=m_1,\,\ldots,\,M_n=m_n \right) \\[2mm] \notag
&&\leq C_{\vep} \sum_{n=1}^{\infty} \sum_{k=1}^n   \sum_{m_k=n+1}^{\infty} \sum_{m_1=0}^{\infty} \cdots \sum_{m_{k-1}=0}^{\infty} \sum_{m_{k+1}=0}^{\infty} \cdots \sum_{m_N=0}^{\infty} (1+\vep)^{n+m_1+\cdots +m_n} \\[2mm] \notag
&&\times\int_0^t  F*G\left(x\,e^{r\,s}\,A\right)\,\PP\left(N(t-s)=n-1 \right)\,\lambda(ds)\,\PP\left(M_1=m_1,\,\ldots,\,M_n=m_n \right)\\[2mm]
&&\lesssim C_{\vep} (1+\omega)\,\E_M \left[(1+\vep)^M\,{\bf 1}_{\{M>N\}} \right] \sum_{n=1}^{\infty} n\, (1+\vep)^{n} \left(\E\left[(1+\vep)^M \right] \right)^{n-1} \\[2mm] \notag
&&\times \int_0^t \PP\left({\bf X}\in x\,e^{r\,s}\,A\right)\,\PP\left(N(t-s)=n-1\right)\,\lambda(ds)\leq C_{\vep} \,(1+\omega)\,\\[2mm] \notag
&& \times \E_M\left[(1+\vep)^M\,{\bf 1}_{\{M>N\}} \right]\,\E_{N(t)}\left[(N(t)+1)\,\left\{(1+\vep)\,\E_M\left[(1+\vep)^M \right] \right\}^{N(t)} \right] \,h(x;\,t)\,,
\eeam
Since the distribution of $M$ has light tail, for any $\delta'>0$, we can find some small enough $\vep>0$ and large enough 
$N_3 \in \bbn$, with $N_3 > N_1\bigvee N_2$, such that it holds
\beam \label{eq.KPY.3.33} 
\E_M\left[(1+\vep)^M\,{\bf 1}_{\{M>N\}} \right] < \delta'\,,
\eeam
for any $N > N_3$. Further, because $\{N(t)\,,\;t\geq 0\}$ is renewal process with finite renewal function, there exists 
some $\delta^*>0$, such that $\E\left[e^{\delta^*\,N(t)} \right]< \infty$, for any $t \in \Lambda_T$. Hence by the Caushy-Schwarz 
inequality we can find some small enough $\vep >0$, such that 
\beam \label{eq.KPY.3.34}
\E_{N(t)}\left[(N(t)+1)\,\left\{(1+\vep)\,\E_M\left[(1+\vep)^M \right] \right\}^{N(t)} \right]=:C_2 <\infty \,,
\eeam
for any $t \in \Lambda_T$. Therefore, from relations \eqref{eq.KPY.3.32} - \eqref{eq.KPY.3.34} it holds
\beam \label{eq.KPY.3.35}
I_{12}(x,\,N,\,r,\,t) \lesssim \delta'\,C_{\vep}\,(1+\omega)\,C_2\,h(x;\,t) \,,
\eeam
for any $t \in \Lambda_T$. From the \eqref{eq.KPY.3.24}, \eqref{eq.KPY.3.31}and \eqref{eq.KPY.3.35}, in combination with 
\eqref{eq.KPY.3.19} and the arbitrary choice of $\delta' >0$, we conclude that \eqref{eq.KPY.3.1} holds uniformly for $t \in \Lambda_T$.
\item[(ii)] 
In this case we have an important difference in relations \eqref{eq.KPY.3.22} and \eqref{eq.KPY.3.32}, where to show that 
$F*G \in \mathcal{S}_A$, we use Lemma \ref{lem.KPY.3.4}, and instead of \eqref{eq.KPY.3.21} we use relation  \eqref{eq.KPY.3.18}. 
Furthermore, in the fourth step of relation  \eqref{eq.KPY.3.25}, we use Lemma \ref{lem.KPY.3.1} in its complete form. 
\end{enumerate}
~\halmos

\section{Infinite time horizon} \label{sec.KPY.4}

In this section we present the second main result, that focuses in the case of infinite time horizon. For the proof, we need 
some closure properties of the distribution classes $\mathcal{A}^*$ and $\mathcal{A}_A^*$, that have their own merit.

\subsection{Main result}  \label{subsec.KPY.4.1}

In contrast to Theorem \ref{th.KPY.3.1}, here we restrict ourselves in class $\mathcal{A}_A^*$, instead of $\mathcal{S}_A$, 
however, this restriction is practically unnoticeable. Here, we also consider positive interest rate $r>0$, a necessary condition 
for non-defectiveness of randomly weighted sums of the form ${\bf D}_r(\infty)$. At the second case we need the extra condition 
$G \in \mathcal{S}_A$, in contrast to Theorem \ref{th.KPY.3.1}(ii), in which we did not assume  some distribution class for $G$. 
As we shall see below, in Remark \ref{rem.KPY.4.5}, this assumption can be relaxed in the form \eqref{eq.KPY.4.23}, or even in some of classes of $F$ can be removed, see Lemma \ref{lem.KPY.4.7}. 
However for sake of simplicity in the formulation of Theorem \ref{th.KPY.4.1} we keep the condition $G \in \mathcal{S}_A$. 

\bth \label{th.KPY.4.1}
Let $A \in \mathscr{R}$ some fixed set. We consider the discounted aggregate claims of relation \eqref{eq.KPY.1.1}, with $r>0$. 
We assume that Assumptions \ref{ass.KPY.1.1} and  \ref{ass.KPY.1.2} are satisfied, and there exists some $\delta>0$ such that 
$\E\left[e^{\delta\,M} \right] < \infty$.
\begin{enumerate}
\item[(i)]
If $F,\,G\,\in \mathcal{A}_A^*$ and $G(x\,A) \asymp F(x\,A)$, then it holds
\beam \label{eq.KPY.4.1}
&&\PP\left({\bf D}_r(\infty) \in x\,A \right)\\[2mm] \notag
&&\sim \int_0^{\infty} \PP\left({\bf X}\in x\,e^{r\,s}\,A \right)\,\lambda(ds) + \E[M]\,\int_0^{\infty} \int_0^{\infty} \PP\left({\bf Y}\in x\,e^{r\,(s+y)}\,A \right)\,H(dy)\,\lambda(ds)\,.
\eeam
\item[(ii)]
If $F\,\in \mathcal{A}_A^*$, $G\,\in \mathcal{S}_A$ and $G(x\,A) =o\left[ F(x\,A)\right]$, then it holds
\beam \label{eq.KPY.4.2} 
\PP\left({\bf D}_r(\infty) \in x\,A \right)\sim \int_0^{\infty} \PP\left({\bf X}\in x\,e^{r\,s}\,A \right)\,\lambda(ds)\,.
\eeam
\end{enumerate}
\ethe

\bre \label{rem.KPY.4.1}
As in Theorem \ref{th.KPY.3.1}, also in Theorem \ref{th.KPY.4.1}, if $F,\,G \in \mathcal{A}_{\mathscr{R}}^*$, and $G(x\,A) \asymp F(x\,A)$ 
for any $A \in \mathscr{R}$ (and $F \in \mathcal{A}_{\mathscr{R}}^*$, $G \in \mathcal{S}_{\mathscr{R}}$ and $G(x\,A) =o\left[F(x\,A)\right]$ 
for any $A \in \mathscr{R}$, respectively), then relation \eqref{eq.KPY.4.1} (relation \eqref{eq.KPY.4.2}, respectively), holds for 
any $A \in \mathscr{R}$. In the one-dimensional subcase, when $A=(1,\,\infty)$, Theorem \ref{th.KPY.4.1} still presents new result, 
since it extends \cite[Lem. 3.6]{li:2023b} (see also Lemma \ref{lem.KPY.4.7} for part (ii)). Theorem \ref{th.KPY.4.1}, can also gives immediately the infinite-time ruin probability in the case where the premiums have bounded densities. 
Hence our results in one-dimensional subcase can generalize \cite[Th. 2.2]{li:2023b}, under more relaxed conditions on premiums and more general distribution classes.
\ere

\bre \label{rem.KPY.4.2}
Using similar methodology, as in Remarks \ref{rem.KPY.3.3}, \ref{rem.KPY.3.4}, under the conditions of Examples \ref{exam.KPY.2.1} and 
\ref{exam.KPY.2.2} with marginals from class $\mathcal{A}^*$, which also have weak equivalent tails, then relations \eqref{eq.KPY.3.4}, \eqref{eq.KPY.3.5}, \eqref{eq.KPY.3.7}, 
\eqref{eq.KPY.3.8}, still hold with $t=\infty$. In order to use the dominated convergence theorem on above relations, it is enough to use the extension \cite[Lemma 5]{tang:yuan:2014}  
for regression dependent random variables. This can be proved by similar lines with the proof of their Lemma, under the using of \cite[Lemma 8]{wang:2011}.
\ere

The reduction of relations  \eqref{eq.KPY.4.1}, \eqref{eq.KPY.4.2}, to \eqref{eq.KPY.4.3}, \eqref{eq.KPY.4.4}, has 
the positive feature that the dependence structure of the components of ${\bf X}$ and ${\bf Y}$ are described completely by the Radon measures $\mu$ and $\nu$, respectively. 
Note also that in the second part we have remove the condition $G\,\in \mathcal{S}_A$, due to Lemma \ref{lem.KPY.4.7}.
 
\bco \label{cor.KPY.4.1}
Let $A \in \mathscr{R}$ be some fixed set, and consider the the discounted aggregate claims of \eqref{eq.KPY.1.1} with $r>0$.
\begin{enumerate}
\item[(i)]
Under the conditions of Theorem \ref{th.KPY.4.1}(i), with restrictions $F \in MRV(\alpha,\,B,\,\mu)$ and $G \in MRV(\alpha,\,Q,\,\nu)$, 
with $\alpha \in (0,\,\infty)$ and $\mu(A)\wedge\nu(A)>0$, it holds
 \beam \label{eq.KPY.4.3} 
\PP\left({\bf D}_r(\infty) \in x A \right)\sim \mu(A) \bB(x) \dfrac{\E\left[e^{-\alpha\,r\,\theta_1} \right]}{1-\E\left[e^{-\alpha r \theta_1} \right]} +\nu(A) \overline{Q}(x) \E[M] \dfrac{\E\left[e^{-\alpha r D} \right] \E\left[e^{-\alpha r \theta_1} \right]}{1-\E\left[e^{-\alpha r \theta_1}\right]}.
\eeam 
\item[(ii)]
Under the conditions of Theorem \ref{th.KPY.4.1}(ii), with $F \in MRV(\alpha,\,B,\,\mu)$, for some $\alpha \in (0,\,\infty)$, $\mu(A)>0$, and $G$ is an arbitrary distribution (such that $G(x\,A) =o\left[ F(x\,A)\right]$), then it holds
\beam \label{eq.KPY.4.4} 
\PP\left({\bf D}_r(\infty) \in x\,A \right)\sim \mu(A)\,\bB(x)\,\dfrac{\E\left[e^{-\alpha\,r\,\theta_1} \right]}{1-\E\left[e^{-\alpha\,r\,\theta_1} \right]}\,.
\eeam 
\end{enumerate}
\eco

\bre \label{rem.KPY.A*}
One of the quantities of interest on ruin theory is the distribution of the first entrance time, of the form
$\tau(x) = \inf\{ t>0\;:\;{\bf D}_r(t) \in x\,A\}$. Such quantities have been studied deeply in the 
random walk theory, but not in the case of risk models with interest rate. Usually we are focused 
on the study of conditional first entrance time $\{\tau(x)\;|\;\tau(x) < \infty \}$, since the $\tau(x)$ 
represents a defective random variable (note that under the conditions of Theorem \ref{th.KPY.4.1} we obtain 
$\PP(\tau(x) < \infty) = \PP\left({\bf D}_r(\infty) \in x\,A \right) <1$, let us recall that the 
$x\,A \in \mathscr{R}$, and consequently it is an increasing set, so the first entrance time of ${\bf D}_r(t)$ 
in $x\,A$ implies that the ${\bf D}_r(s)$ belongs to  $x\,A$, for any $t<s\leq \infty$). Therefore,  for any $t>0$ we find that 
\beam \label{eq.KPY.4.4.1} 
\PP\left(\tau(x)\leq t\;|\; \tau(x) < \infty \right) = \dfrac{\PP\left(\tau(x)\leq t\,,\;\tau(x) < \infty \right)}{\PP\left(\tau(x) < \infty \right)}= \dfrac{\PP\left({\bf D}_r(t) \in x\,A\right)}{\PP\left({\bf D}_r(\infty) \in x\,A \right)}\,.
\eeam 
In one-dimensional and multidimensional risk models, without delayed claims, and under the assumptions of Theorem \ref{th.KPY.4.1}, 
with restrictions of $F_A,\,G_A \in \mathcal{R}_{-\alpha}$, with $\alpha \in (0,\,\infty)$ and $\{N(t)\,,\;t \geq 0\}$ be a Poisson process 
with intensity $\lambda >0$, it is well-known that the $\{\tau(x)\;|\;\tau(x) < \infty \}$ is exponentially distributed 
with parameter $\alpha\,r$, see in \cite[Sec. 2]{hao:tang:2008} and \cite[Rem. 3.6]{chen:konstantinides:passalidis:2025} 
for the one-dimensional and multidimensional cases, respectively.
It is not difficult to see that the same is true here, in the case of Theorem \ref{th.KPY.4.1} (ii) (under the restrictions 
$F_A \in \mathcal{R}_{-\alpha}$, and $\{N(t)\,,\;t \geq 0\}$ be a Poisson process with intensity $\lambda >0$). However, in case of 
 Theorem \ref{th.KPY.4.1} (i) under the same restrictions, due to \eqref{eq.KPY.4.4.1}, \eqref{eq.KPY.3.1} and \eqref{eq.KPY.4.1}, 
through \eqref{eq.KPY.2.7}, for any fixed $t>0$ we obtain 
\beam \label{eq.KPY.4.4.2} \notag
&&\PP\left(\tau(x)\leq t\;|\;\tau(x) < \infty \right) \\[2mm]\notag
&&\sim  \dfrac{\lambda\,F(x\,A)\left[\dfrac{1-e^{-\alpha\,r\,t}}{\alpha\,r} \right]+ \lambda\,G(x\,A)\,\E[M]\int_0^t e^{-\alpha\,r\,s}\E\left[e^{-\alpha\,r\,D}{\bf 1}_{\{D \leq t-s\}}\right]ds}{\dfrac {\lambda}{\alpha\,r}\,F(x\,A)+\dfrac {\lambda}{\alpha\,r}\,G(x\,A)\E[M]\E[e^{-\alpha\,r\,D}]} \\[2mm]
&& \leq \max \left\{ 1-e^{-\alpha\,r\,t}\,,\; \dfrac{\alpha\,r \int_0^t e^{-\alpha\,r\,s}\E\left[e^{-\alpha\,r\,D}{\bf 1}_{\{D \leq t-s\}}\right]ds}{\E\left[e^{-\alpha\,r\,D}\right]}\right\} \leq 1- e^{-\alpha\,r\,t}\,.
\eeam 
From \eqref{eq.KPY.4.4.2} we see that if $Z^* \stackrel{d}{\sim} V^*$ with $V^*$ an exponential distribution with parameter 
$\alpha\,r > 0$, then $Z^* \leq_{st} \{\tau(x)\;|\;\tau(x) < \infty \}$ (as $x\to\infty$). Thus, in case where the delayed claims 
have asymptotically equivalent tails with the main claims, the first entrance time has not lighter tail than the exponential tail with 
parameter $\alpha\,r$, which is the expected outcome without delayed claims. This is interesting, since 
is in conflict with the intuition that 'when there exist important delayed claims the rare-event appears earlier'.
\ere

\bre \label{rem.KPY.4.B*}
The uniformity of the asymptotic estimations for all $t \in \Lambda$ (global uniformity) is important 
as from mathematical aspect, but also from practical point of view as well, see \cite[Sec. 3]{tang:2004} for more discussion on this topic. 
Here, under the conditions of Theorem \ref{th.KPY.4.1}(ii), relation \eqref{eq.KPY.3.2} can be proved to hold uniformly 
for all $t \in \Lambda$. Indeed, for this purpose, it is enough to take into account the local uniformity of \eqref{eq.KPY.3.2},  
that is implied by Theorem \ref{th.KPY.3.1} (ii), and Theorem \ref{th.KPY.4.1} (ii). By this fact, following the same treatment with 
\cite[Th. 3.2]{konstantinides:liu:passalidis:2025}, we can prove such a result. Correspondingly, under the assumptions of 
Corollary \ref{cor.KPY.4.1} (ii), relation \eqref{eq.KPY.3.c} holds uniformly for $t \in \Lambda$. However, under the 
conditions of Theorem \ref{th.KPY.4.1} (i), it is not a trivial task to get the uniformity of  \eqref{eq.KPY.3.1} for all $t \in \Lambda$. 
From the best of our knowledge, only few papers on delayed claims have tried to derive uniform estimations on the whole $\Lambda$.
Even more, we find serious mathematical drawbacks in their proofs, although the corresponging arguments for $t \in \Lambda_T$ and $t= \infty$ are correct.
\ere

\subsection{Closure properties of $\mathcal{A}^*$ and $\mathcal{A}_A^*$}  \label{subsec.KPY.4.2}

The closure properties of heavy-tailed distributions, play significant role in several branches of applied probability, see for example
\cite{leipus:siaulys:konstantinides:2023}. In this section we provide some lemmas on closure properties of the distribution 
classes $\mathcal{A}^*$ and $\mathcal{A}_A^*$, that are necessary for the proof of Theorem  \ref{th.KPY.4.1}. 

We begin with some properties of class  $\mathcal{A}^*$. In next result we establish sufficient and necessary conditions for the 
closure property of  $\mathcal{A}^*$ with respect to convolution. It is inspired by \cite[Th. 1.1]{leipus:siaulys:2020}, which 
studies such conditions in the class of convolution equivalence distributions. We recall that for two independent random variables 
$\Theta_1 \stackrel{d}{\sim} B_1$,  $\Theta_2 \stackrel{d}{\sim} B_2$, we denote with $\Theta_1\bigvee \Theta_2 \stackrel{d}{\sim} B_1\,B_2$, namely
$\PP\left(\Theta_1\vee \Theta_2 \leq x\right) =B_1\,B_2(x)$, for any $x \in \bbr$.

\bpr \label{prop.KPY.4.1}
Let $B_1,\,B_2$ be distributions on $\bbr$, with $B_1,\,B_2  \in \mathcal{T}^*$. Then the following statements are equivalent
\begin{enumerate}
\item[(i)]
$B_1*B_2 \in \mathcal{A}^*$.
\item[(ii)]
$B_1\,B_2 \in \mathcal{A}^*$.
\item[(iii)]
$p\,B_1 + (1-p)\,B_2 \in \mathcal{A}^*$, for any (or, equivalently, for some) $p\in (0,\,1)$.
\end{enumerate}
Further, each of these statements implies that 
\beam \label{eq.KPY.4.5}
\overline{B_1*B_2}(x) \sim \bB_1(x) + \bB_2(x)\,.
\eeam 
\epr

\pr~
At first, since $B_1,\,B_2  \in \mathcal{T}^* \subsetneq \mathcal{L}$ and $\mathcal{A}^*\subsetneq \mathcal{S}$, from \cite[Th. 1.1]{leipus:siaulys:2020}, 
each one of the statements $(i) - (iii)$ implies that \eqref{eq.KPY.4.5} holds. Hence, for any $v>1$ it holds,
\beao
&&\overline{(B_1*B_2)^*}(v) \\[2mm]
&&=\limsup \dfrac{\overline{B_1*B_2}(v\,x)}{\overline{B_1*B_2}(x)}\leq \limsup \dfrac{\bB_1(v\,x) + \bB_2(v\,x)}{\bB_1(x) + \bB_2(x)}\leq \limsup \left[ \dfrac{\overline{B_1}(v\,x)}{\overline{B_1}(x)}\bigvee \dfrac{\overline{B_2}(v\,x)}{\overline{B_2}(x)} \right] \\[2mm]
&& \leq \max\left\{ \limsup \dfrac{\overline{B_1}(v\,x)}{\overline{B_1}(x)}\,,\;\limsup \dfrac{\overline{B_2}(v\,x)}{\overline{B_2}(x)}\right\}\leq \max \{\overline{B_1^*}(v),\,\overline{B_2^*}(v)\}\,,
\eeao
 so we obtain
\beam \label{eq.KPY.4.6}
K_{B_1*B_2}^- =-\lim_{v\downarrow 1} \dfrac{\log\overline{(B_1*B_2)^*}(v)}{\log v} \geq -\lim_{v\downarrow 1} \dfrac{\log \left(\max \left\{\overline{B_1^*}(v),\,\overline{B_2^*}(v) \right\} \right)}{\log v} >0\,,
\eeam
where at the last step, we used that $B_1,\,B_2  \in \mathcal{T}^*$ and thus $K_{B_1}^-\wedge K_{B_2}^- >0$.

Next, we shall show the equivalence of assertions $(i) - (iii)$.

$(ii) \Rightarrow (i)$. We assume that $B_1\,B_2 \in \mathcal{A}^*$. Then from \cite[Th. 1.1]{leipus:siaulys:2020} follows $B_1*B_2 \in \mathcal{S}$, 
and further from \eqref{eq.KPY.4.6} we obtain that $B_1*B_2 \in \mathcal{A}^*$.

$(i) \Rightarrow (ii)$. We assume that $B_1*B_2 \in \mathcal{A}^*$. Then by \cite[Th. 1.1]{leipus:siaulys:2020} we find $B_1\,B_2 \in \mathcal{S}$. 
It remains to show that $K_{B_1\,B_2}^- >0$. At first we see that 
\beao
\overline{B_1\,B_2}(x) = \PP(\Theta_1 \vee \Theta_2 > x)=\bB_1(x) + \bB_2(x) - \bB_1(x)\,\bB_2(x) \sim \bB_1(x) + \bB_2(x)\,,
\eeao 
hence, it holds
\beao
&&\overline{(B_1\,B_2)^*}(v)  \\[2mm]
&&= \limsup \dfrac{\overline{B_1\,B_2}(v\,x)}{\overline{B_1\,B_2}(x)} \leq \limsup \dfrac{\bB_1(v\,x) + \bB_2(v\,x)}{\bB_1(x) + \bB_2(x)}\leq \limsup \left[ \dfrac{\overline{B_1}(v\,x)}{\overline{B_1}(x)}\bigvee \dfrac{\overline{B_2}(v\,x)}{\overline{B_2}(x)} \right] \\[2mm]
&& \leq \max\left\{ \limsup \dfrac{\overline{B_1}(v\,x)}{\overline{B_1}(x)}\,,\;\limsup \dfrac{\overline{B_2}(v\,x)}{\overline{B_2}(x)}\right\}\leq \max \{\overline{B_1^*}(v),\,\overline{B_2^*}(v)\}\,,
\eeao
for any $v>1$, and therefore we obtain
\beao
K_{B_1\,B_2}^- = -\lim_{v\downarrow 1} \dfrac{\log \overline{(B_1\,B_2)^*}(v)}{\log v} \geq -\lim_{v\downarrow 1} \dfrac{\log \left[\max\left\{\overline{B_1^*}(v)\,,\;\overline{B_2^*}(v) \right\} \right]}{\log v} > 0\,,
\eeao 
where at the last step we used the fact that $B_1, B_2 \in  \mathcal{T}^*$. Hence, we conclude that $B_1\,B_2 \in  \mathcal{A}^*$.

$(i) \Rightarrow (iii)$. Let $B_1*B_2 \in \mathcal{A}^*$. Then, by \cite[Th. 1.1]{leipus:siaulys:2020} holds $p\,B_1 + (1-p)\,B_2 \in \mathcal{S}$, 
for any (or, equivalently, for some) $p \in (0,\,1)$.

Further, it holds
\beao
&&\overline{\left(p\,B_1 + (1-p)\,B_2 \right)^*}(v)  \\[2mm]
&&= \limsup \dfrac{p\,\bB_1(v\,x) + (1-p)\,\bB_2(v\,x)}{p\,\bB_1(x) + (1-p)\,\bB_2(x)}\leq \limsup \left[ \dfrac{\overline{B_1}(v\,x)}{\overline{B_1}(x)}\bigvee \dfrac{\overline{B_2}(v\,x)}{\overline{B_2}(x)} \right] \\[2mm]
&& \leq \max \left\{\limsup \dfrac{\bB_1(v\,x) }{\bB_1(x)}\,,\; \limsup \dfrac{\bB_2(v\,x)}{\bB_2(x)}\right\} \leq \max \left\{\overline{B_1^*}(v)\,,\;\overline{B_2^*}(v) \right\}\,,
\eeao
for any (or, equivalently, for some) $p \in (0,\,1)$. Thus, it holds
\beao
&&K_{p\,B_1 + (1-p)\,B_2 }^-  \\[2mm]
&&= -\lim_{v\downarrow 1} \dfrac{\log \overline{\left(p\,B_1 + (1-p)\,B_2 \right)^*}(v)}{\log v}\geq -\lim_{v\downarrow 1} \dfrac{\log \left[\max \left\{\overline{B_1^*}(v)\,,\;\overline{B_2^*}(v) \right\} \right]}{\log v} > 0\,,
\eeao
for any (or, equivalently, for some) $p \in (0,\,1)$, which implies that $p\,B_1 + (1-p)\,B_2 \in \mathcal{A}^*$ for any (or, equivalently, for some) $p \in (0,\,1)$.

$(iii) \Rightarrow (i)$. Let for some $p \in (0,\,1)$, it holds $p\,B_1 + (1-p)\,B_2 \in \mathcal{A}^*$. Then for any $p \in (0,\,1)$, holds $p\,B_1 + (1-p)\,B_2 \in \mathcal{S}$, and further by \cite[Th. 1.1]{leipus:siaulys:2020} we obtain $B_1*B_2 \in \mathcal{S}$. Since relation \eqref{eq.KPY.4.5} is true, via \eqref{eq.KPY.4.6} we conclude that $B_1*B_2 \in \mathcal{A}^*$.   
~\halmos

If in Proposition \ref{prop.KPY.4.1} restrict the marginal distributions from class $\mathcal{T}^*$ to class $\mathcal{A}^*$, then 
relation \eqref{eq.KPY.4.5} becomes equivalent to any of assertions $(i) - (iii)$. This is described in the next result.

\ble \label{lem.KPY.4.1}
Let $B_1,\,B_2$ be distributions on $\bbr$, with $B_1,\,B_2  \in \mathcal{A}^*$. Then each of the statements $(i) - (iii)$ is equivalent 
to relation \eqref{eq.KPY.4.5}.
\ele

\pr~
Since  $B_1,\,B_2  \in \mathcal{A}^* \subsetneq \mathcal{T}^*$, by Proposition \ref{prop.KPY.4.1}, we obtain that the statements 
$(i) - (iii)$ are equivalent and imply  relation \eqref{eq.KPY.4.5}. So it remains to show that  relation \eqref{eq.KPY.4.5} implies $(i)$.

Let assume that relation \eqref{eq.KPY.4.5} is true. Then by \cite[Cor. 1.1]{leipus:siaulys:2020} it follows that $B_1*B_2 \in \mathcal{S}$. 
Further, from \eqref{eq.KPY.4.5}, we obtain \eqref{eq.KPY.4.6}, and consequently we conclude $B_1*B_2  \in \mathcal{A}^*$.       
~\halmos

In the following result, we provide some other sufficient conditions for the closure property of class $\mathcal{A}^*$ with respect to convolution.

\ble \label{lem.KPY.4.2}
Let $B_1,\,B_2$ be distributions on $\bbr$. 
\begin{enumerate}
\item[(i)]
If $B_1 \in \mathcal{A}^*$, $B_2 \in \mathcal{T}^*$ and $\bB_2(x)=O\left[\bB_1(x)\right]$, then \eqref{eq.KPY.4.5} holds, and $B_1*B_2 \in \mathcal{A}^*$.
\item[(ii)]
If $B_1 \in \mathcal{A}^*$ and $\bB_2(x)=o\left[\bB_1(x)\right]$, then 
\beam \label{eq.KPY.4.7}
\overline{B_1*B_2}(x) \sim \bB_1(x)\,,
\eeam
and $B_1*B_2 \in \mathcal{A}^*$.
\end{enumerate}
\ele

\pr~
\begin{enumerate}
\item[(i)]
Since $B_1 \in \mathcal{A}^* \subsetneq \mathcal{S}$, $B_2 \in \mathcal{T}^* \subsetneq \mathcal{L}$ and  $\bB_2(x)=O\left[\bB_1(x)\right]$, 
from \cite[Lem. 3.2]{tang:tsitsiashvili:2003} we obtain that \eqref{eq.KPY.4.5} holds and $B_1*B_2 \in \mathcal{S}$. From \eqref{eq.KPY.4.5}, 
via  \eqref{eq.KPY.4.6} we find that $K_{B_1*B_2}^- >0$. Hence $B_1*B_2 \in \mathcal{A}^*$.  
\item[(ii)]
From \cite[Cor. 3.18]{foss:korshunov:zachary:2013}, we find that relation \eqref{eq.KPY.4.7} is true, as also the $B_1*B_2 \in \mathcal{S}$. 
By relation \eqref{eq.KPY.4.7} follows immediately that $K_{B_1*B_2}^- = K_{B_1}^- >0$. ~\halmos
\end{enumerate}

Now, we proceed to closure properties of class $\mathcal{A}_A^*$ with respect to convolution. We notice that for two multidimensional 
distributions $V_1,\,V_2$, with support on the  non-negative orthant, the condition $V_1*V_2 \in \mathcal{B}_A$, with $A \in \mathscr{R}$ 
and some class $\mathcal{B}$, means that if ${\bf Z}^{(1)},\,{\bf Z}^{(2)}$ are independent random vectors with distributions  $V_1,\,V_2$, 
respectively, then the random variable $Z_A' = \sup\left\{u\;:\;{\bf Z}^{(1)}+{\bf Z}^{(2)} \in u\,A \right\}\stackrel{d}{\sim} V_A'
$ is such that it holds $V_A' \in \mathcal{B}$. See for more discussions \cite[Sec. 3.2]{konstantinides:passalidis:2024g}. 
The following statement says that when the marginal distributions $V_1,\,V_2$ belong to class $\mathcal{L}_A$, with the $V_A^{(1)},\,V_A^{(2)}$ to have positive lower Katamata index, the closure problem of  $\mathcal{A}_A^*$, with respect to convolution, is reduced to the corresponding one-dimensional problem.

\bpr \label{prop.KPY.4.2}
Let $A \in \mathscr{R}$ be some fixed set. If $V_1,\,V_2 \in \mathcal{L}_A$, with 
$K_{V_A^{(1)}}^- \wedge K_{V_A^{(2)}}^- >0$, then, $V_1*V_2 \in \mathcal{A}_A^*$, if and only if, $V_A^{(1)}*V_A^{(2)} \in \mathcal{A}^*$.
\epr

\pr~
$(\Rightarrow)$. Let $V_1*V_2 \in \mathcal{A}_A^*$. Then, by \cite[Th. 3.4]{konstantinides:passalidis:2024g} we find that 
$V_A^{(1)}*V_A^{(2)} \in \mathcal{S}$, and since  $V_A^{(1)},\,V_A^{(2)} \in \mathcal{T}^* \subsetneq \mathcal{L}$, from 
\cite[Th. 1.1]{leipus:siaulys:2020}, we obtain that it holds
\beam \label{eq.KPY.4.8}
\overline{V_A^{(1)}*V_A^{(2)}}(x) \sim \overline{V_A^{(1)}}(x)+\overline{V_A^{(2)}}(x)\,.
\eeam 
Similarly as for \eqref{eq.KPY.4.6}, here again via \eqref{eq.KPY.4.8}, we obtain that $K_{V_A^{(1)}*V_A^{(2)}}^- >0$. This, 
in combination with the previous assertions provides the relation $V_A^{(1)}*V_A^{(2)} \in \mathcal{A}^*$.  

$(\Leftarrow)$. Let $V_A^{(1)}*V_A^{(2)} \in \mathcal{A}^*$. From \cite[Th. 3.4]{konstantinides:passalidis:2024g} we obtain 
that $V_1*V_2 \in \mathcal{S}_A$, hence  $Z_A' \stackrel{d}{\sim} V_A' \in \mathcal{S}$. It remains to show $K_{V_A'}^- >0$. 

Let ${\bf Z}^{(1)},\,{\bf Z}^{(2)}$ non-negative random vectors, independent each other, with ${\bf Z}^{(1)} \stackrel{d}{\sim} V_1$ 
and ${\bf Z}^{(2)} \stackrel{d}{\sim} V_2$. Then for any $v>1$, it holds
\beao
&&\overline{(V_A')^*}(v) = \limsup \dfrac{\PP(Z_A' > v\,x)}{\PP(Z_A' > x)}\\[2mm]
&&= \limsup \dfrac{\PP({\bf Z}^{(1)}+{\bf Z}^{(2)} \in v\,x\,A)}{\PP({\bf Z}^{(1)}+{\bf Z}^{(2)} \in x\,A)} \leq \limsup \dfrac{\PP(Z_A^{(1)}+Z_A^{(2)} > v\,x)}{\PP(Z_A^{(1)}> x) + \PP(Z_A^{(2)} > x)} \\[2mm]
&&= \limsup \dfrac{\PP(Z_A^{(1)}> v\,x)+\PP(Z_A^{(2)}> v\,x)}{\PP(Z_A^{(1)}> x)+\PP(Z_A^{(2)}> x)} \leq\max\left\{\overline{\left(V_A^{(1)}\right)^*}(v)\,,\; \overline{\left(V_A^{(2)}\right)^*}(v) \right\}\,,
\eeao
where at the third step we used \cite[Prop. 2.4]{konstantinides:passalidis:2024g} in the numerator, and 
Bonferroni inequality in the denominator, at the fourth step we apply Proposition \ref{prop.KPY.4.1}, 
since $V_A^{(1)}*V_A^{(2)} \in \mathcal{A}^*$, $V_A^{(1)}\,,\;V_A^{(2)} \in \mathcal{T}^*$, hence relation \eqref{eq.KPY.4.8} 
is true. From the last formula we obtain
\beao
K_{V_A'}^- = -\lim_{v\downarrow 1} \dfrac{\log \overline{\left(V_A'\right)^*}(v)}{\log v} \geq -\lim_{v\downarrow 1} \dfrac{\log \left[\max\left\{\overline{\left(V_A^{(1)}\right)^*}(v)\,,\;\overline{\left(V_A^{(2)}\right)^*}(v) \right\} \right]}{\log v} >0\,,
\eeao
where at the last step we used that $V_A^{(1)}\,,\;V_A^{(2)} \in \mathcal{T}^*$. Therefore, we have $V_A' \in \mathcal{A}^*$, 
that implies $V_1*V_2 \in \mathcal{A}_A^*$.  
~\halmos

\bre \label{rem.KPY.4.3}
Direct consequence of Proposition \ref{prop.KPY.4.2}, through \cite[Cor. 3.1]{konstantinides:passalidis:2024g}, is the assertion: 
if $V_1,\,V_2 \in \mathcal{L}_A$ and $V_1*V_2 \in \mathcal{A}_A^*$ then it holds
\beam \label{eq.KPY.4.9}
V_1*V_2(x\,A) \sim V_1(x\,A)+V_2(x\,A)\,,
\eeam
without necessarily $K_{V_A^{(1)}}^- \vee K_{V_A^{(2)}}^- >0$. Additionally, if in Lemma \ref{lem.KPY.3.4}, instead of 
$V_1 \in \mathcal{S}_A$, we get restricted to  $V_1 \in \mathcal{A}_A^*$, then from relation \eqref{eq.KPY.3.18} we obtain 
that $V_1*V_2 \in \mathcal{A}_A^*$.
\ere

Finally, we provide a result on closure property of class $\mathcal{A}^*$ with respect to convolution product. 

\ble \label{lem.KPY.4.3}
Let $Z,\,W$ two independent non-negative random variables with $Z \stackrel{d}{\sim} B \in \mathcal{A}_A^*$ and $\PP(W=0)<1$, 
$\PP(W > c\,x) = o\left[\overline{Q}(x)\right]$, for any $c>0$, where $W\,Z \stackrel{d}{\sim} Q$. Then $Q \in \mathcal{A}_A^*$.
\ele

\pr~
At first, from \cite[Th. 1.3]{xu:cheng:wang:cheng:2018} we obtain that $Q \in \mathcal{S}$. It remains to show $K_{Q}^- >0$.

For any fixed $v>1$ and $\delta' >0$, we can find some $x_0=x_0(\delta') >0$, large enough, such that it holds
\beao
\bB(v\,x) \leq (1+ \delta')\,\overline{B^*}(v)\,\bB(x)\,,
\eeao
for any $x \geq x_0$. Furthermore, for the concrete $x_0>0$, and by assumption on the heaviness of tails of $W$ and $WZ$, we can find some $x_1> x_0$, such that $\PP\left( W > \dfrac x{x_0}\right) \leq \delta'\,\overline{Q}(x)$, for any $x_1 \geq x_0$. Thus, we obtain
\beao
&&\overline{Q}(v\,x) \leq \int_0^{x/x_0} \bB\left( \dfrac{v\,x}y\right)\,\PP(W \in dy) +\PP\left( W> \dfrac{x}{x_0}\right) \\[2mm]
&&\leq (1+\delta')\,\overline{B^*}(v)\,\int_0^{x/x_0} \bB\left( \dfrac{x}y\right)\,\PP(W \in dy) + \delta'\,\overline{Q}(x) \leq \left[(1+\delta')\,\overline{B^*}(v)+ \delta'\right]\,\overline{Q}(x)\,,
\eeao
for any $x > x_1$. From the last relation, via the arbitrary choice of $\delta'>0$, we find that 
\beam \label{eq.KPY.4.10}
\overline{Q^*}(v) \leq \overline{B^*}(v)\,,
\eeam
hence we have that
\beao
K_Q^- = -\lim_{v\downarrow 1} \dfrac{\log \overline{Q^*}(v)}{\log v} \geq -\lim_{v\downarrow 1} \dfrac{\log \overline{B^*}(v)}{\log v} = K_B^- >0\,,
\eeao
that gives $Q \in \mathcal{A}_A^*$.
~\halmos

\bre \label{rem.KPY.4.4}
We should mention that the tail condition of $W$ is quite common and handy, especially in case, with infinite upper Matuszewska 
index. We refer to \cite[Sec. 2]{tang:2006} for more discussions on this condition. We notice that if the distribution of $W$ has 
bounded from above support, then since $B$ has infinite right endpoint (that implies $Q$ has also infinite right endpoint), the heavy-tailed condition 
for $W$ in Lemma \ref{lem.KPY.4.3} holds automatically. In fact we shall use here Lemma \ref{lem.KPY.4.3} only in the case of bounded from above $W$.

Finally we note that if $A \in \mathscr{R}$, and ${\bf Z} \stackrel{d}{\sim} V \in \mathcal{A}_A^*$, then the question if the 
distribution of $W\,{\bf Z}$ belongs to $\mathcal{A}_A^*$ is reduced to the question if the distribution of $W\,Z_A$ belongs 
to $\mathcal{A}^*$, since it holds
\beao
\PP(W\,{\bf Z} \in x\,A) = \PP\left( \sup_{{\bf p} \in I_A} {\bf p}^{\top}\,W\,{\bf Z}>x\right) =\PP(W\,Z_A > x)\,,
\eeao
for any $x>0$.
\ere

\subsection{Section \ref{subsec.KPY.4.1} argumentation}  \label{subsec.KPY.4.3}

Now, we provide the proof of Theorem \ref{th.KPY.4.1} after some auxiliary lemmas. The first one is a multivariate reformulation 
of \cite[Lem. 5.2]{tang:yuan:2016}. The argument follows trivially with the help of Remark \ref{rem.KPY.4.4}.

\ble \label{lem.KPY.4.4}
Let $A \in \mathscr{R}$ be some fixed set. We consider a non-negative random vector 
${\bf Z} \stackrel{d}{\sim} V$ with $K_{V_A}^->0$, and a random variable $W$, independent of ${\bf Z}$, such that it holds 
$\PP(0\leq W \leq 1)=1$. Then, for any $\vep \in (0,\,K_{V_A}^- )$ and for any $c>1$, there exists some $x_0 =x_0(\vep,\,c)>0$, such that it holds
\beam \label{eq.KPY.4.11}
\dfrac{\PP(W\,{\bf Z} \in x\,A)}{\PP({\bf Z} \in x\,A)}  \leq c\,\E\left[W^{K_{V_A}^- -\vep} \right]\,,
\eeam
for all $x > x_0$.
\ele

The next lemma is a special case of \cite[Lem. 3.1]{konstantinides:passalidis:2025o}.

\ble \label{lem.KPY.4.5}
Let $A \in \mathscr{R}$ be some fixed set. We consider non-negative random vectors 
${\bf Z}^{(i)} \stackrel{d}{\sim} V_i \in \mathcal{S}_A$, for $i=1,\,2$, and a random variable $W$, independent of 
${\bf Z}^{(1)},\,{\bf Z}^{(2)}$, with $\PP(0 \leq W \leq b)=1$ for some $b \in (0,\,\infty)$, and $\PP(W =0)<1$. If, either 
$V_1(x\,A) = O\left[ V_2(x\,A) \right]$, or $V_2(x\,A) = O\left[V_1(x\,A) \right]$, then 
\beam \label{eq.KPY.4.12}
\PP(W\,[{\bf Z}^{(1)}+{\bf Z}^{(2)}] \in x\,A) \sim \PP(W\,{\bf Z}^{(1)} \in x\,A) + \PP(W\,{\bf Z}^{(2)} \in x\,A) \,.
\eeam
\ele

The following lemma represents a variation of Lemma \ref{lem.KPY.3.1}, and concretely it has an increased range of uniformity, 
namely instead of the interval $[a,\,b]$, with $0<a\leq b < \infty$, we have the $(0,\,b]$, with $b \in (0,\,\infty)$, however 
the distribution class get restricted from $\mathcal{S}_A$ to $\mathcal{A}_A$. 

\ble \label{lem.KPY.4.6}
Let $A \in \mathscr{R}$ be some fixed set. We consider $n+m$ non-negative, independent random vectors 
${\bf Z}^{(1)},\,\ldots,\,{\bf Z}^{(n+m)}$ with distributions $V_1,\,\ldots,\, V_{n+m}$, respectively. We assume that 
 $V_1,\,\ldots,\, V_{n} \in \mathcal{A}_A$, with $V_i(x\,A) \asymp V_1(x\,A)$, for any $i=2,\,\ldots,\,n$, and 
$V_j(x\,A) =o\left[ V_1(x\,A)\right]$, for any $j = n+1,\,\ldots,\,n+m$.
Let some constant $b\in (0,\,\infty)$, such that ${\bf c}_n \in (0,\,b]^n$ and the $(c_{n+1},\,\ldots,\,c_{n+m})^{\top}$ have values 
from $[0,\,b]^m$, and for any $j =1,\,\ldots,\,n+m$ there exists some $i_j \in \{1,\,\ldots,\,n\}$ 
with $c_j \leq c_{i_j}$. Then relation \eqref{eq.KPY.3.9} holds uniformly for ${\bf c}_n \in (0,\,b]^n$. 
\ele

\pr~
At first we show that it holds
\beam \label{eq.KPY.a}
\PP\left( \sum_{i=1}^n c_i\,{\bf Z}^{(i)} \in x\,A \right) \sim  \sum_{i=1}^n\PP\left( c_i\,{\bf Z}^{(i)} \in x\,A \right) 
\eeam
uniformly for ${\bf c}_n \in (0,\,b]^n$. From the one hand side we obtain
\beao
&&\PP\left( \sum_{i=1}^n c_i {\bf Z}^{(i)} \in x A \right) \leq \PP\left( \sum_{i=1}^n c_i Z_A^{(i)} > x \right) \sim \sum_{i=1}^n \PP\left( c_i Z_A^{(i)} > x \right)  \\[2mm]
&&=\sum_{i=1}^n\PP\left( c_i {\bf Z}^{(i)} \in x A \right)
\eeao 
uniformly for ${\bf c}_n \in (0,\,b]^n$, where at the first step we  used \cite[Prop. 2.4]{konstantinides:passalidis:2024g}, and at the 
second step we applied \cite[Lem. 5, Lem. 1]{tang:yuan:2014}. 

For the lower bound of \eqref{eq.KPY.a}, since the summands are non-negative, and the set $xA$ increasing, by Bonferroni inequality we find that it holds
\beao
&&\PP\left( \sum_{i=1}^n c_i\,{\bf Z}^{(i)} \in x\,A \right) \geq \PP\left( \bigcup_{i=1}^n \left\{ c_i {\bf Z}^{(i)} \in x A \right\}\right) \\[2mm]
&&=\sum_{i=1}^n \PP\left(c_i {\bf Z}^{(i)} \in x A \right) - \sum_{1\leq i < j \leq n} \PP\left(c_i {\bf Z}^{(i)} \in x A\,,\;c_j {\bf Z}^{(j)} \in x A \right) \sim \sum_{i=1}^n \PP\left(  c_i {\bf Z}^{(i)} \in x A \right),
\eeao
uniformly for ${\bf c}_n \in (0,\,b]^n$, where at the last step we took into account
the independence of ${\bf Z}^{(1)},\,\ldots,\,{\bf Z}^{(n)}$. Hence, relation  \eqref{eq.KPY.a} holds  uniformly for any ${\bf c}_n \in (0,\,b]^n$. 
Following similar path, as in the proof of Lemma \ref{lem.KPY.3.1}, but now using relation \eqref{eq.KPY.a}, we can establish  \eqref{eq.KPY.3.9}  
uniformly for any ${\bf c}_n \in (0,\,b]^n$. We  only mention that at the third step of  \eqref{eq.KPY.3.13} 
we now use the \cite[Lem. 5, Lem. 1]{tang:yuan:2014}.   
~\halmos

\noindent{\bf Proof of Theorem \ref{th.KPY.4.1}.}~
We examine the assertion $(i)$ in full details, and we provide only a sketch of the proof of the assertion $(ii)$, since it follows similar path.
\begin{enumerate}
\item[(i)]
From Assumption \ref{ass.KPY.1.1} we obtain
\beam \label{eq.KPY.4.13}
&&\PP\left({\bf D}_r(\infty)\in x\,A \right) = \PP\left(\sum_{i=1}^{\infty} \left[{\bf X}^{(i)}\,e^{- r\,\tau_i} + \sum_{j=1}^{M_i} {\bf Y}^{(i,j)}\,e^{- r\,(\tau_i + D_{ij})} \right] \in x\,A\right) \\[2mm] \notag
&&=\PP\left(\sum_{i=1}^{\infty} \left[{\bf X}^{(i)} + \sum_{j=1}^{M_i} {\bf Y}^{(i,j)} e^{- r\, D_{ij}} \right] e^{- r\,\tau_i} \in x A\right)=: \PP\left(\sum_{i=1}^{\infty} \left[{\bf X}^{(i)} + {\bf S}^{(i)} \right] e^{- r\,\tau_i} \in x A\right) ,
\eeam
where ${\bf S}^{(i)} = \sum_{j=1}^{M_i} {\bf Y}^{(i,j)}\,e^{- r\, D_{ij}}$ and thus from Assumptions \ref{ass.KPY.1.1} and \ref{ass.KPY.1.2}, 
the sequence $\left\{ {\bf S}^{(i)}\,,\;i\in \bbn \right\}$ contains i.i.d. random vectors with common distribution $G_{\bf S}$.

Since the distribution of ${\bf Y}$ is $G \in \mathcal{A}^*$ and $e^{- r\,D} \in [0,\,1]$, from Lemma \ref{lem.KPY.4.3}, see also Remark 
\ref{rem.KPY.4.4}, we find that the $\left\{  {\bf Y}^{(i,j)}\,e^{- r\, D_{ij}}\,,\;i,\,j \in \bbn \right\}$ are i.i.d. with common 
distribution from class $\mathcal{A}_A^*$.

From the fact that the distribution of $M$ has light tail, and  $\mathcal{A}_A^* \subsetneq \mathcal{S}_A$, therefore from 
\cite[Th. 4.2]{konstantinides:passalidis:2024g}, we obtain
\beam \label{eq.KPY.4.14} 
G_{\bf S}(x\,A) =\PP\left({\bf S} \in x\,A \right) = \PP\left(\sum_{i=1}^{M_i} {\bf Y}^{(i,j)}\,e^{- r\, D_{i,j}} \in x\,A\right)\sim \E[M]\,\PP\left({\bf Y}\,e^{- r\, D} \in x\,A\right)\,,
\eeam 
hence, if $\E[M]>0$, from closure property of class $\mathcal{A}_A^*$ with respect to strong tail equivalence, we have that $G_{\bf S} \in \mathcal{A}_A^*$. 
Further, since $G(x\,A) \asymp F(x\,A)$, via \eqref{eq.KPY.4.14} we obtain that $G_{\bf S}(x\,A) = O\left[ F(x\,A)\right]$. We should 
mention, that in the subcase where $\E[M]=0$, hence $M_i=0$ almost surely, the proof becomes immediate, without need to determine the 
distribution $G_{\bf S}$, since ${\bf S}={\bf 0}$, almost surely.

Considering the random variable $S_A := \sup \left\{u\;:\; {\bf S} \in u\,A \right\}\stackrel{d}{\sim} G_{{\bf S},A}$, then from 
$F_A,\,G_{{\bf S},A} \in \mathcal{A}^*$, and 
\beao
\bG_{{\bf S},A}(x) =O\left[\bF_A(x) \right]\,,
\eeao 
by Lemma \ref{lem.KPY.4.2}(i) we get $F_A * G_{{\bf S},A} \in \mathcal{A}^*$. Taking into account Proposition \ref{prop.KPY.4.2} this is 
equivalent to $F * G_{{\bf S}} \in \mathcal{A}_A^*$ and keeping in mind Remark \ref{rem.KPY.4.3} and \eqref{eq.KPY.4.9} we conclude 
\beam \label{eq.KPY.4.15}
F * G_{{\bf S}}(x\,A) \sim F(x\,A) + G_{\bf S}(x\,A)\,.
\eeam

Hence, the $\left\{ {\bf X}^{(i)} + {\bf S}^{(i)} \,,\;i \in \bbn \right\}$ is a sequence of i.i.d. random vectors, with common 
distribution  $F * G_{{\bf S}} \in \mathcal{A}_A^*$, that satisfies relation \eqref{eq.KPY.4.15}. But since $r >0$, we find 
$\PP\left(0 \leq e^{-r\,\theta_1} \leq 1\right)=1$, $\PP\left( e^{-r\,\theta_1}=1\right) <1$ and $\PP\left( e^{-r\,\theta_1}=0\right) <1$, 
where the last follows from the fact that the $\theta_1$ do not have mass at infinity, hence from \cite[Lem. 3.4]{konstantinides:passalidis:2025o}, 
see also Remark \ref{rem.KPY.3.3} of that paper, it holds
\beam \label{eq.KPY.4.16}
\PP\left(\sum_{i=1}^{n} \left[{\bf X}^{(i)} + {\bf S}^{(i)} \right]\,e^{- r\,\tau_i} \in x\,A\right) \sim \sum_{i=1}^{n} \PP\left(\left[{\bf X}^{(i)} + {\bf S}^{(i)} \right]\,e^{- r\,\tau_i} \in x\,A\right) \,.
\eeam
uniformly for $n \in \bbn$, and consequently also for $n = \infty$, in the sense that
the following limit is true
\beao
\lim \sup_{n \in\bbn} \left| \dfrac{\PP\left(\sum_{i=1}^{n} \left[{\bf X}^{(i)} + {\bf S}^{(i)} \right]\,e^{- r\,\tau_i} \in x\,A\right)}{\sum_{i=1}^{n} \PP\left(\left[{\bf X}^{(i)} + {\bf S}^{(i)} \right]\,e^{- r\,\tau_i} \in x\,A\right)} - 1 \right| =0\,.
\eeao

Next, we shall show that for arbitrary $\delta' >0$, we can find some large enough $N = N(\delta') \in \bbn$, such that it holds
\beam \label{eq.KPY.4.17} \notag
&&\sum_{i=N+1}^{n}\PP\left( \left[{\bf X}^{(i)} + {\bf S}^{(i)} \right]\,e^{- r\,\tau_i} \in x\,A\right) \\[2mm]
&&\lesssim \delta' \left[\PP\left({\bf X}\,e^{- r\,\theta_1} \in x\,A\right) + \E[M]\,\PP\left({\bf Y}\,e^{- r\,(\theta_1 + D)} \in x\,A\right) \right] \,.
\eeam
Relation \eqref{eq.KPY.4.17} represent a special case of the first relation of \cite[P. 24]{konstantinides:passalidis:2025o}. However, 
for sake of compactness, we provide below the full argument
\beao
&&\sum_{i=N+1}^{n}\PP\left( \left[{\bf X}^{(i)} + {\bf S}^{(i)} \right]\,e^{- r\,\tau_i} \in x\,A\right) \\[2mm]
&&= \sum_{i=N+1}^{n} \int_0^1 \int_0^1 \PP\left({\bf X}^{(i)} + {\bf S}^{(i)} \in  \dfrac x{y\,t}\,A\right)\,\PP(e^{- r\,\tau_{i-1}} \in dt)\,\PP(e^{- r\,\theta_1} \in dy) \\[2mm]
&&\sim \sum_{i=N+1}^{n} \int_0^1 \int_0^1 \left[ \PP\left({\bf X}^{(i)} \in  \dfrac x{y\,t}\,A\right) + \PP\left({\bf S}^{(i)} \in  \dfrac x{y\,t}\,A\right) \right]\,\PP(e^{- r\,\tau_{i-1}} \in dt)\,\PP(e^{- r\,\theta_1} \in dy) \\[2mm]
&&\sim \sum_{i=N+1}^{n} \int_0^1 \int_0^1 \left[ \PP\left({\bf X}^{(i)} \in  \dfrac x{y\,t}\,A\right) + \E[M]\,\PP\left({\bf Y}\,e^{-r\,D} \in  \dfrac x{y\,t}\,A\right) \right]\,\PP(e^{- r\,\tau_{i-1}} \in dt)\\[2mm]
&&\times \PP(e^{- r\,\theta_1} \in dy)\\[2mm]
&& = \sum_{i=N+1}^{n} \int_0^1 \left[ \PP\left({\bf X}^{(i)} e^{- r \tau_{i-1}} \in  \dfrac x{y} A\right) + \E[M] \PP\left({\bf Y}\,e^{- r \tau_{i-1}} e^{-r\,D} \in  \dfrac x{y} A\right) \right] \PP(e^{- r \theta_1} \in dy)\\[2mm]
&&\leq \sum_{i=N+1}^{n}c\,\left\{\E\left[ \left(e^{- r\,\tau_{i-1}}  \right)^{K_{F_A}^- -\vep} \right]\bigvee \E\left[ \left(e^{- r\,\tau_{i-1}}  \right)^{K_{G_A'}^- -\vep} \right] \right\}\, \\[2mm]
&&\times \int_0^1 \left[ \PP\left({\bf X} \in  \dfrac x{y}\,A\right) + \E[M]\,\PP\left({\bf Y}\,e^{-r\,D} \in  \dfrac x{y}\,A\right) \right]\,\PP(e^{- r\,\theta_1} \in dy) \\[2mm]
&&\leq c\,\left\{\dfrac{\left(\E\left[ \left(e^{- r\,\theta_1} \right)^{K_{F_A}^- -\vep} \right] \right)^N}{1-\E\left[ \left(e^{- r\,\theta_1} \right)^{K_{F_A}^- -\vep} \right]} + \dfrac{\left(\E\left[ \left(e^{- r\,\theta_1} \right)^{K_{G_A'}^- -\vep} \right] \right)^N}{1-\E\left[ \left(e^{- r\,\theta_1} \right)^{K_{G_A'}^- -\vep} \right]}  \right\}\, \\[2mm]
&&\times \left[ \PP\left({\bf X}\,e^{- r\,\theta_1} \in x\,A\right) + \E[M]\,\PP\left({\bf Y}\,e^{-r\,(\theta_1+D)} \in x\,A\right) \right]\,,
\eeao
where at the second step we used Lemma \ref{lem.KPY.4.5} and relation \eqref{eq.KPY.4.12}, at the third step we take advantage of 
\eqref{eq.KPY.4.14}, at the fifth step we apply twice Lemma \ref{lem.KPY.4.4} one for ${\bf X}$ and the other for 
${\bf Y}\,e^{-r\,D}  \in \mathcal{A}_A^*$, with $Y_A\,e^{-r\,D}\stackrel{d}{\sim} G'_{A}
$. In this application we choose $\vep \in \left(0,\,K_{F_A}^- \wedge K_{G_A'}^- \right)$. Hence, from the inequalities 
\beao
\E\left[ \left(e^{- r\,\theta_1} \right)^{K_{F_A}^- -\vep} \right]<1\,,\qquad \E\left[ \left(e^{- r\,\theta_1} \right)^{K_{G_A'}^- -\vep} \right]<1\,,
\eeao
for any $\delta' >0$ we can find large enough $N=N(\delta') \in \bbn$, such that relation \eqref{eq.KPY.4.17} holds.

Therefore, it holds
\beam \label{eq.KPY.4.18} \notag
&&\sum_{i=1}^{n}\PP\left( \left[{\bf X}^{(i)} + {\bf S}^{(i)} \right]\,e^{- r\,\tau_i} \in x\,A\right) \geq \sum_{i=1}^{N}\PP\left( \left[{\bf X}^{(i)} + {\bf S}^{(i)} \right]\,e^{- r\,\tau_i} \in x\,A\right)\\[2mm]
&&\sim \left(\sum_{i=1}^{n} - \sum_{i=N+1}^{n} \right) \left[\PP\left({\bf X}^{(i)}\,e^{- r\,\tau_i} \in x\,A\right) + \E[M]\,\PP\left({\bf Y}\,e^{- r\,(\tau_i + D)} \in x\,A\right) \right]\\[2mm] \notag
&&\geq (1-\delta')\left[\sum_{i=1}^{n} \PP\left({\bf X}^{(i)}\,e^{- r\,\tau_i} \in x\,A\right) + \E[M]\,\sum_{i=1}^{n} \PP\left({\bf Y}\,e^{- r\,(\tau_i + D)} \in x\,A\right) \right] \,.
\eeam
for any $n>N$, where at the second step used Lemma \ref{lem.KPY.4.5} and at the third step we take into consideration \eqref{eq.KPY.4.17} 
(similarly due to  Lemma \ref{lem.KPY.4.5} in combination with relation  \eqref{eq.KPY.4.17}).

From the other hand side, via similar arguments we obtain
\beam \label{eq.KPY.4.19} \notag
&&\sum_{i=1}^{n}\PP\left( \left[{\bf X}^{(i)} + {\bf S}^{(i)} \right]\,e^{- r\,\tau_i} \in x\,A\right) = \left(\sum_{i=1}^{N} + \sum_{i=N+1}^{n} \right)\PP\left( \left[{\bf X}^{(i)} + {\bf S}^{(i)} \right]\,e^{- r\,\tau_i} \in x\,A\right)\\[2mm] \notag
&&\sim \sum_{i=1}^{N}  \left[\PP\left({\bf X}^{(i)}\,e^{- r\,\tau_i} \in x\,A\right) + \E[M]\,\PP\left({\bf Y}\,e^{- r\,(\tau_i + D)} \in x\,A\right) \right] \\[2mm] 
&&+ \sum_{i=N+1}^{n} \PP\left( \left[{\bf X}^{(i)} + {\bf S}^{(i)} \right]\,e^{- r\,\tau_i} \in x\,A\right) \\[2mm] \notag
&&\lesssim (1-\delta')\left[\sum_{i=1}^{n} \PP\left({\bf X}^{(i)}\,e^{- r\,\tau_i} \in x\,A\right) + \E[M]\,\sum_{i=1}^{n} \PP\left({\bf Y}\,e^{- r\,(\tau_i + D)} \in x\,A\right) \right] \,.
\eeam

From \eqref{eq.KPY.4.18} and \eqref{eq.KPY.4.19} and the arbitrary choice of $\delta' > 0$, we find that it holds
\beam \label{eq.KPY.4.20}
&&\sum_{i=1}^{n}\PP\left( \left[{\bf X}^{(i)} + {\bf S}^{(i)} \right]\,e^{- r\,\tau_i} \in x\,A\right) \\[2mm]  \notag
&& \sim \sum_{i=1}^{n}  \PP\left({\bf X}^{(i)}\,e^{- r\,\tau_i} \in x\,A\right) + \E[M] \sum_{i=1}^{n} \,\PP\left({\bf Y}\,e^{- r\,(\tau_i + D)} \in x\,A\right) \,,
\eeam
uniformly for $n>N$. Further, via Lemma \ref{lem.KPY.4.5}, we obtain that relation  \eqref{eq.KPY.4.20} holds uniformly for $n \leq N$. 
Hence, \eqref{eq.KPY.4.20} holds uniformly for $n \in \bbn$. As result, we see that from \eqref{eq.KPY.4.16} and \eqref{eq.KPY.4.20} is 
implied the relation
\beam \label{eq.KPY.4.21} 
&&\PP\left( \sum_{i=1}^{n}\left[{\bf X}^{(i)} + {\bf S}^{(i)} \right]\,e^{- r\,\tau_i} \in x\,A\right)  \\[2mm] \notag
&&\sim \sum_{i=1}^{n}  \PP\left({\bf X}^{(i)}\,e^{- r\,\tau_i} \in x\,A\right) + \E[M] \sum_{i=1}^{n} \,\PP\left({\bf Y}\,e^{- r\,(\tau_i + D)} \in x\,A\right) \,,
\eeam
uniformly for $n \in \bbn$. So, applying \eqref{eq.KPY.4.21} for $n=\infty$ on relation \eqref{eq.KPY.4.13} we obtain
\beao
&&\PP\left({\bf D}_r(\infty)\in x\,A \right) \sim \sum_{i=1}^{\infty}\PP\left( {\bf X}^{(i)}\,e^{- r\,\tau_i} \in x\,A\right)+ \E[M] \sum_{i=1}^{n} \,\PP\left({\bf Y}\,e^{- r\,(\tau_i + D)} \in x\,A \right) \\[2mm] 
&&=\int_0^{\infty} \PP\left({\bf X}^{(i)}  \in x\,e^{r\,s}\,A\right)\,\lambda(ds)+\E[M] \int_0^{\infty} \int_0^{\infty} \PP\left({\bf Y} \in x\,e^{r\,(s + y)}\,A\right) \,H(dy)\,\lambda(ds)\,,
\eeao
which gives relation \eqref{eq.KPY.4.1}.

\item[(ii)] 
Since $G \in \mathcal{S}_A$, the distribution of $Y\,e^{-r\,D}$ belongs to class $\mathcal{S}_A$, see \cite[Th. 2.1]{cline:samorodnitsky:1994} 
in combination with Remark \ref{rem.KPY.4.4}, or \cite[Th. 3.2(iv)]{konstantinides:passalidis:2024g}. 

Hence relation \eqref{eq.KPY.4.14} holds with $G_{\bf S} \in \mathcal{S}_A$. Due to asymptotic relation $G(x\,A) =o\left[F(x\,A) \right]$ it 
follows that $G_{\bf S}(x\,A) =o\left[F(x\,A) \right]$, so by Remark \ref{rem.KPY.4.3} and by \eqref{eq.KPY.3.18} we find the relation
\beam \label{eq.KPY.4.22} 
F*G_{\bf S}(x\,A) \sim F(x\,A)\,,
\eeam
and thus the $\left\{ {\bf X}^{(i)} + {\bf S}^{(i)} \,,\;i \in \bbn \right\}$ is a sequence of i.i.d. random vectors, with common distribution 
$F * G_{{\bf S}} \in \mathcal{A}_A^*$. Next, the only difference, exists in derivation of relation \eqref{eq.KPY.4.17}, where instead of Lemma 
\ref{lem.KPY.4.5}, we should use relation \eqref{eq.KPY.3.15}, that is derived by Lemma \ref{lem.KPY.4.6}, see also Remark \ref{rem.KPY.3.5}, 
for expression of the form
\beao
\PP\left( \left[{\bf X}^{(i)} + {\bf S}^{(i)} \right]\,e^{- r\,\tau_i} \in x\,A\right) \sim \PP\left( {\bf X}^{(i)}\,e^{- r\,\tau_i} \in x\,A\right) \,,
\eeao
since here the inequality $0<e^{- r\,\tau_i} \leq_{st} e^{- r\,\tau_i} \leq 1$ is satisfied as equality. Further, we follow the same path of 
argumentation as before.
\end{enumerate}
 ~\halmos

\bre \label{rem.KPY.4.5}
In fact, as we see from the first assertion of Theorem \ref{th.KPY.4.1}(ii), the assumption $G\,\in \mathcal{S}_A$, can be relaxed into
\beam \label{eq.KPY.4.23}
G_{\bf S}(x\,A) =o\left[F(x\,A) \right]\,,
\eeam
that includes several cases of delayed claims, as for example the case when the distribution $G$ has light tail. However, we keep the condition 
$G\,\in \mathcal{S}_A$ on the Theorem, for the sake of simplicity of the assumption, but \eqref{eq.KPY.4.23} contains many other cases.
\ere

The following lemma provides a sufficient condition for the validity of relation \eqref{eq.KPY.4.23}, and is inspired from \cite[Lem. 3.3]{li:2023b}. 
Namely, this lemma permit the complete removal of the conditions $G \in \mathcal{S}_A$, in the case of $F \in \mathcal{C}_A$, with $K_{F_A}^- > 0$, where 
we say $F \in \mathcal{C}_A$, if $F_A \in \mathcal{C}$, in the sense of
\beao
\lim_{b\uparrow 1} \limsup \dfrac{\bF_A(b\,x)}{\bF_A(x)} = 1\,.
\eeao 
We recall that we say $F \in (\mathcal{C} \cap \mathcal{P_D})_A$, if $F_A \in (\mathcal{C} \cap \mathcal{P_D})$, and $\mathcal{R}_{-\alpha} \subsetneq \mathcal{C} \cap \mathcal{P_D}$.
In the following Lemma we need also the upper Matuszewska index of a distribution $B$, which defined as:
\beao
J_B^+ = -\lim_{v\to \infty} \dfrac{\log \overline{B_{*}}(v)}{\log v}\,,
\eeao
where
\beao
\overline{B_{*}}(v)= \liminf \dfrac {\bB(v\,x)}{\bB(x)}\,,
\eeao
Recall also that if $F_A\in\mathcal{C}$ then $J_{F_A}^+<\infty$. 
Hence, the following lemma is helpful in the proof of Corollary \ref{cor.KPY.4.1} (ii).

\ble \label{lem.KPY.4.7}
Let $A \in \mathscr{R}$ be some fixed set and $\left\{{\bf Z}^{(i)}\,,\;i \in \bbn \right\}$ be a sequence of i.i.d. non-negative random vectors 
with distribution $V$. We also assume that $F$ is a distribution with support on the non-negative orthant, such that $F \in (\mathcal{C} \cap \mathcal{P_D})_A $, 
and $V(x\,A) = o[F(x\,A)]$. Further, we suppose that $\left\{ \Theta_{i}\,,\;i \in \bbn \right\}$ represents a sequence of non-negative random 
weights, independent of $\left\{{\bf Z}^{(i)}\,,\;i \in \bbn \right\}$, and such that for some 
$0<p_1 < J_{F_A}^- \leq J_{F_A}^+ < p_2 < \infty $ it holds
\beam \label{eq.KPY.4.27}
\sum_{i=1}^{\infty} \left( \E\left[\Theta_i^{p_1} \right] \bigvee \E\left[\Theta_i^{p_2} \right]\right)^{1/p} < \infty\,,
\eeam
where $p = {\bf 1}_{\{0 < J_{F_A}^+ < 1\}} +p_2\,{\bf 1}_{\{1 \leq J_{F_A}^+ < \infty \}}$. Then we obtain
\beam \label{eq.KPY.4.28}
\PP\left( \sum_{i=1}^{\infty} \Theta_i\,{\bf Z}^{(i)} \in x\,A \right)= o[F(x\,A)]\,.
\eeam
\ele

\pr~
Firstly, for any $x>0$ it holds
\beao
&&\PP\left( \sum_{i=1}^{\infty} \Theta_i\,{\bf Z}^{(i)} \in x\,A \right)=\PP\left( \sup_{{\bf p}\in I_A} {\bf p}^{\top} \left[\sum_{i=1}^{\infty} \Theta_i\,{\bf Z}^{(i)}\right] > x \right)\\[2mm]
&&\leq \PP\left( \sum_{i=1}^{\infty} \Theta_i\, Z_A^{(i)} > x \right)\,,
\eeao
Further, following the line of the proof of \cite[Lem. 3.3]{li:2023b}, but now using \cite[Th. 2]{yi:chen:su:2011} (due to \eqref{eq.KPY.4.27}), 
instead of Lemma 3.2 in that paper, we can find that
\beao
\PP\left( \sum_{i=1}^{\infty} \Theta_i\, Z_A^{(i)} > x \right)= o[\bF_A(x)]\,.
\eeao
Hence, by these two last relations we conclude \eqref{eq.KPY.4.28}.
~\halmos

Now we can prove Corollary \ref{cor.KPY.4.1}.

\noindent{\bf Proof of Corollary 4.1}~
(i) The proof follows immediately from relation \eqref{eq.KPY.4.1} by application of relation \eqref{eq.KPY.2.6a}, 
since $F_A,\,G_A \in \mathcal{R}_{-\alpha}$ (remind that $\mu(A)\wedge\nu(A)>0$). The application of the dominated convergence theorem, is possible 
via the Potter inequalities for the class $\mathcal{R}_{-\alpha}$. 
\\(ii) Firstly, since $F \in MRV(\alpha,\,B,\,\,\mu)$ and $\mu(A)>0$, we obtain $F_A \in \mathcal{R}_{-\alpha}$, and thus we have
$F_A \in \mathcal{C}$, with $K_{F_A}^-=J_{F_A}^- =\alpha >0$; therefore it follows 
that $F_A \in \mathcal{C} \cap \mathcal{P_D}$. Now, the result is implied by relation \eqref{eq.KPY.4.2} (with the same arguments with part (i)), 
because of Lemma \ref{lem.KPY.4.7} (keeping in mind that $r>0$ and $\{D_{ij}\,,\;i,\,j \in \bbn\}$ are identically distributed).
~\halmos

\vspace{1em}

\noindent\textbf{Funding Statement}: This work is supported by the Basic Scientific Research Project of the Department of Education of Liaoning Province (No. LJ212510173010) and Doctoral Research Start-up Program of the Department of Science and Technology of Liaoning Province (No. 2026-BS-0753).


\begin{thebibliography}{99}

%\selectlanguage{american}

%AAAAAAAAAAAAAAAAAAAAAAAAAAAAAAAAAAAAAAAAAAAAAAAAAAAAAAAA 

%\bibitem{albrecher:asmussen:2006}
%{\sc Albrecher, H., Asmussen, S.}\ (2006) 
%Ruin probabilities and aggregate claims distributions for shot noise Cox processes.
%{\em  Scand. Actuar. J.} \textbf{2}, 86--110.

%\bibitem{asimit:furman:tang:vernic:2011}
%{\sc Asimit, A.V., Furman, E., Tang, Q., Vernic, R.}\ (2011) 
%Asymptotics for risk capital allocations based on Conditional Tail Expectation.
%{\em  Insur. Math. Econom.} \textbf{49}, 310--324.

%\bibitem{asmussen:schmidli:schmidt:1999}
%{\sc Asmussen, S., Schmidli, H., Schmidt, V.}\ (1999) 
%Tail probabilities for non-standard risk and queueing processes with subexponential jumps.
%{\em  Adv. Appl. Probab.} \textbf{31}, 442--447.

%BBBBBBBBBBBBBBBBBBBBBBBBBBBBBBBBBBBBBBBBBBBBBBBBBBBBBBBB

%\bibitem{bardoutsos:konstantinides:2011}
%{\sc Bardoutsos, A.G., Konstantinides, D.G.}\ (2011)
%Characterization of tails through hazard rate and convolution closure properties. 
%{\em J. Appl. Probab.}, \textbf{48A}, 123--132.

%\bibitem{basrak:davis:mikosch:2002}
%{\sc Basrak, B., Davis, R.A., Mikosch, T.}\ (2002) 
%A characterization of multivariate regular variation.
%{\em  Ann. Appl. Probab.} \textbf{12}, 908--920.

%\bibitem{basrak:davis:mikosch:2002b}
%{\sc Basrak, B., Davis, R.A., Mikosch, T.}\ (2002) 
%Regular variation of GARCH processes.
%{\em  Stoch. Process. Appl.} \textbf{99}, no. 1, 95--115.

%\bibitem{basrak:segers:2009}
%{\sc Basrak, B., Segers, J.}\ (2009) 
%Reularly varying multivariate time series.
%{\em  Stoch. Process. Appl.} \textbf{119}, 1055--1080.

%\bibitem{beirlant:goegebeur:segers:teugels:2004}
%{\sc Beirlant, J., Goegebeur, Y., Segers, J., Teugels, J.L.}\ (2004)
%{\em Statistics of extremes: theory and applications} 
%Wiley, Vol. 558, Chichester.

\bibitem{bingham:goldie:teugels:1987} 
{\sc Bingham. N.H., Goldie, C.M., Teugels, J.L.} \ (1987)
{\em Regular Variation.}
Cambridge University Press, Cambridge.

\bibitem{buraczewski:damek:mikosch:2016} 
{\sc Buraczewski, D., Damek, E., Mikosch, T.}\ (2016) 
{\em Stochastic Models with Power-Law Tails}.
Springer, New York.

%CCCCCCCCCCCCCCCCCCCCCCCCCCCCCCCCCCCCCCCCCCCCCCCCCCCCCCC

%\bibitem{cai:tang:2004}
%{\sc Cai, J., Tang, Q.}\  (2004)
%On max-sum equivalence and convolution closure of heavy-tailed distributions and their applications. 
%{\em J. Appl. Probab.} \textbf{41}, 117--130.

%\bibitem{chen:blanchet:rhee:zwart:2019}
%{\sc Chen, B., Blanchet, J., Rhee, C.H., Zwart, B.}\  (2019)
%Efficient Rare-Event Simulation for Multiple Jump Events in Regularly Varying Random Walks and Compound Poisson Processes. 
%{\em Math. Oper. Resear.} \textbf{44}, no. 3, 919--942.

%\bibitem{chen:xu:cheng:2019}
%{\sc Chen, J., Xu, H., Cheng, F.}\  (2019)
%The product of dependent random variables with applications to a discrete-time risk model. 
%{\em Commun. Stat. Theory Methods} \textbf{48}, 3325--3340.

%\bibitem{chen:2011}
%{\sc Chen, Y.}\  (2011)
%The finite-time ruin probabilities with dependent insurance and financial risks. 
%{\em J. Appl. Probab.} \textbf{48}, 1035--1048.

%\bibitem{chen:2017}
%{\sc Chen, Y.}\  (2017)
%Interplay of subexponential and dependent insurance and financial risks. 
%{\em Insur. Math. Econom. } \textbf{77}, 78--83.

%\bibitem{chen:cui:wang:2024}
%{\sc Chen, Y., Cui, Z., Wang, Y.}\ (2024)
%Precise large deviations of some objectives related t the net loss process in two nonstandard risk modes. 
%{\em Preprint, arXiv: 2305.00475}.

%\bibitem{chen:cheng:2024b}
%{\sc Chen, Z., Cheng, D.}\ (2024)
%Precise large deviations for non-centralized sums of partial sums and random sums of heavy-tailedd END random variables. 
%{\em Stat. Probab. Lett.}, \textbf{211}, 110134.

%\bibitem{chen:liu:2022}
%{\sc Chen, Y., Liu, J.}\  (2022)
%An asymptotic study of systemic expected shortfall and marginal expected shortfall. 
%{\em Insur. Math. Econom.} \textbf{105}, 238--251.

%\bibitem{chen:yuen:2009}
%{\sc Chen, Y., Yuen, K.C.}\ (2009)
%Sums of pairwise quasi-asymptotic independent random variables with consistent variation.
%{\em Stochastic Models}, \textbf{25}, 76--89.

%\bibitem{chen:wang:shang:2021}
%{\sc Chen, Y., Wang, J., Zhang, W.}\ (2021)
%Tail distortion risk measure for prortfolio with multivariate regularly variation. 
%{\em Commun. Math. Statist.}, \textbf{10}, 263--285.

%\bibitem{chen:ng:tang:2005}
%{\sc Chen, Y., Ng, K.W., Tang, Q.}\ (2005)
%Weighted sums of subexponential random variables and their maxima. 
%{\em Adv. Appl. Probab.}, \textbf{37}, 510--522.

%\bibitem{chen:cheng:zheng:2025}
%{\sc Chen, Z., Cheng, D., Zheng, H.}\ (2025)
%On the joint tail behavior of randomly weighted sums of dependent random variables with applications to risk theory. 
%{\em Scand. Actuar. J.},  1-20

\bibitem{chen:konstantinides:passalidis:2025}
{\sc Chen, Z., Konstantinides, D.G., Passalidis, C.D.}\ (2025)
Asymptotics for aggregated interdependet multivariate subexponential claims with general investment returns. 
{\em Preprint, arXiv:2507.23713}.


\bibitem{chen:passalidis:2028}
{\sc Chen, Z., Passalidis, C.D.}\ (2026)
Asymptotics for a continuous-time risk model with multivariate stronly regularly varying claims and stochastic investment returns. 
{\em In preparation}.

%\bibitem{chen:li:cheng:2023}
%{\sc Chen, Z., Li, M., Cheng, D.}\ (2023)
%Asymptotics for sum-ruin probabilities of a bidimensional risk model with heavy-tailed claims and stochastic returns. 
%{\em Stochastics}, \textbf{96}, no.2, 947-967.

\bibitem{chen:wang:wang:2013}
{\sc Chen, Y., Wang, L., Wang, Y.}\ (2013)
Uniform asymptotics for the finite-time ruin probabilities of two kinds of nonstandard bidimensional risk modes. 
{\em J. Math. Anal. Appl.}, \textbf{401}, no. 1, 114--129.

%\bibitem{chen:yang:2019}
%{\sc Chen, Y., Yang, Y.}\  (2019)
%Bivariate regular variation among randomly weighted sums in general insurance. 
%{\em Eur. Actuar. J.} \textbf{9}, 301--322.

%\bibitem{chen:yuan:2017}
%{\sc Chen, Y., Yuan, Z.}\  (2017)
%A rivisit to ruin probabilities in the presence of heavy-tailed insurance and financial risks. 
%{\em Insur. Math. Econom.} \textbf{73}, 75--81.

%\bibitem{cheng:2014}
%{\sc Cheng, D.}\  (2014)
%Randomly weighted sums of dependent random variables with dominated variation. 
%{\em J. Math. Anal. Appl.} \textbf{420}, no. 3, 1617--1633.

%\bibitem{cheng:2021}
%{\sc Cheng, D.}\  (2021)
%Uniform asymptotics for the finite-time ruin probability of a generalized bidimensional risk model with Brownian perturbations.
%{\em Stochastics} \textbf{93}, Vol 1. 56--71.

%\bibitem{cheng:cheng:2018}
%{\sc Cheng, F., Cheng, D.}\ (2018)
%Randomly weighted sums of dependent subexponential random variables with application to risk theory. 
%{\em Scand. Actuar. J.}, \textbf{3}, 191-202.

%\bibitem{cheng:cheng:chen:2021}
%{\sc Cheng, F., Cheng, D., Chen, Z.}\ (2021)
%Asymptotic behavior for finite-time ruin probabilities in a generalized bidimensional risk model with subexponential claims. 
%{\em Jap. J. Industr. Appl. Math.}, \textbf{38}, 947-963.

%\bibitem{cheng:konstantinides:wang:2022}
%{\sc Cheng, M., Konstantinides, D.G., Wang, D}\ (2022)
%Uniform asymptotic estimates in a time-dependent risk model with general investment returns and multivariate regularly varying claims. 
%{\em Appl. Math. and Comput.}, \textbf{434}, 127436.

\bibitem{cheng:konstantinides:wang:2024}
{ \sc Cheng, M., Konstantinides, D.G., Wang, D.}\ (2024)
Multivariate regular varying insurance and financial risks in $d$-dimensional risk model. 
{ \em J. Appl. Probab.}, \textbf{61}, no. 4, 1319 -- 1342.

%\bibitem{cheng:yu:2019}
%{\sc Cheng, D., Yu, C.}\  (2019)
%Uniform asymptotics for the ruin probabilities in a bidimensional renewal risk model with strongly subexponential claims.
%{\em Stochastics} \textbf{91}, Vol 1. 643--656.

%\bibitem{cheng:xu:2020}
%{\sc Cheng, F., Xu, H.}\  (2020)
%The finite-time ruin probability of the nonhomogeneous Poisson risk model with conditionally independent subexponential claims.
%{\em Comm. Stat. Theor. Meth.} \textbf{51}, Vol. 12, 4119--4132.

%\bibitem{cline:resnick:1992} 
%{\sc Cline, D.B.H., Resnick, S.}\ (1992)
%Multivariate subexponential distributions.
%{\em Stoch. Process. Appl.}, \textbf{42}, no.1, 49--72.

\bibitem{cline:samorodnitsky:1994} 
{\sc Cline, D.B.H., Samorodnitsky, G.}\ (1994)
Subexponentiality of the product of independent random variables.
{\em Stoch. Process. Appl.}, \textbf{49}, 75--98.

%DDDDDDDDDDDDDDDDDDDDDDDDDDDDDDDDDDDDDDDDDDDDDDDDDDDDDDD

%\bibitem{dirma:nakiliuda:siaulys:2023} 
%{\sc Dirma, M., Nakliuda, N., {\v{S}}iaulys, J.}\ (2023)
%Generalized moments of sums with heavy-tailed random summands.
%{\em  Lith. Math. J.}, \textbf{63}, no. 3, 254--271.

%EEEEEEEEEEEEEEEEEEEEEEEEEEEEEEEEEEEEEEEEEEEEEEEEEEEEEE

%\bibitem{embrechts:klueppelberg:mikosch:1997}
%{\sc Embrechts, P., Kl\"{u}pellberg, C. and Mikosch, T.}\ (1997) 
%{\em Modelling Extremal Events for Insurance and Finance.} 
%Springer, New York.

%FFFFFFFFFFFFFFFFFFFFFFFFFFFFFFFFFFFFFFFFFFFFFFFF

\bibitem{foss:korshunov:zachary:2013} 
{\sc Foss, S., Korshunov, D., Zachary, S.} \ (2013)
{\em An Introduction to Heavy-Tailed and Subexponential Distributions.}
Springer, New York, 2nd ed. 

%\bibitem{foss:richards:2010} 
%{\sc Foss, S., Richards, A.} \ (2010)
%On sums of conditionally independent subexponential random variables.
%{\em Math. Oper. Resear.}, \textbf{35}, 102--119.

%GGGGGGGGGGGGGGGGGGGGGGGGGGGGGGGGGGGGGGGGGGGGGGGGGGGGGGGGGGGGGGGGGG

%\bibitem{gao:wang:2010}
%{\sc Gao, Q., Wang, Y.}\ (2010)
%Randomly weighted sums with dominatedly varying-tailed increments and application to risk theory. 
%{\em J. Korean Stat. Soc.}, \textbf{39}, 305--314.

\bibitem{gao:yang:2014}
{\sc Gao, Q., Yang, X.}\ (2014)
Asymptotic ruin probabilities in a generalized bidimensional risk model perturbed by diffusion with constant force of interest. 
{\em J. Math. Anal. Appl.}, \textbf{419}, no. 2, 1193--1213.

\bibitem{gao:zhuang:huang:2019}
{\sc Gao, Q., Zhuang, J., Huang, Z.}\ (2019)
Asymptotics for a delay-claim risk model with diffusion, dependence structures and constant force of interest. 
{\em J. Comput. Appl. Math.}, \textbf{353}, 219--231.

%\bibitem{geluk:ng:2006}
%{\sc Geluk, J., Ng, K.W.}\ (2006)
%Tail behavior of negatively associated heavy-tailed sums. 
%{\em J. Appl. Probab.}, \textbf{43}, no. 2, 587--593.

\bibitem{geluk:tang:2009}
{\sc Geluk, J., Tang, Q.}\ (2009)
Asymptotic tail probabilities of sums of dependent subexponential random variables. 
{\em J. Theor. Probab.}, \textbf{22}, 871--882.

%\bibitem{geng:liu:wang:2023}
%{\sc Geng, B., Liu, Z., Wang, S.}\  (2023)
%A Kesten-type inequality for randomly  weighted sums of dependent subexponential random variables with application to risk theory. 
%{\em Lith. Math. J.} \textbf{63}, 81--91.

%\bibitem{goldie:1978}
%{\sc Goldie, C.M.}\ (1978)
%Subexponential distributions and dominated variation tails 
%{\em J. Appl. Probab.}, \textbf{15}, 440--442.

%\bibitem{guo:2022}
%{\sc Guo, F.}\ (2022)
%Ruin probability of a continuous-time model with dependence between insurance  and financial risks caused by systematic risk factors.
%{\em Appl. Math. Comput.}, \textbf{413}, 126634.

%\bibitem{guo:wang:2013}
%{\sc Guo,F., Wang, D.}\ (2013)
%Finite-and infinite-time ruin probabilities with general stochastic investment return processes and bivariate upper tail independent and heavy tailed claims.
%{\em Adv. Appl. Probab.}, \textbf{12}, no.4, 241--273.

%HHHHHHHHHHHHHHHHHHHHHHHHHHHHHHHHHHHHHHHHHHHHHHHHHHHHHHH

%\bibitem{haan:ferreira:2006}
%{\sc Haan, L. de, Ferreira, A.}\ (2006) 
%{\em Extreme Value Theory: An Introduction.} 
%Springer Verlag, New York.

%\bibitem{haan:resnick:1984} 
%{\sc Haan, L. de, Resnick, S.}\ (1984)
%Stochastic compactness and point processes.
%{\em J. Aust. Math. Soc. Ser. A}, \textbf{37}, 307--316.

%\bibitem{hazra:maulik:2012}
%{\sc Hazra, R.S., Maulik, K.}\ (2012)
%Tail behavior of randomly weighted sums. 
%{\em Adv. Appl. Probab.}, \textbf{44}, 794--814.

%\bibitem{haegele:lehtomaa:2021} 
%{\sc Haegele, M., Lehtomaa, J.}\ (2021)
%Large deviations for a class of multivariate heavy-tailed risk processes used in insurance and finance.
%{\em  J. Risk Fin. Manag.}, \textbf{14}, 202.

\bibitem{hao:tang:2008} 
{\sc Hao, X., Tang, Q.}\ (2008)
A uniform asymptotic estimate for discounted aggregate claims with subexponential tails.
{\em  Insur. Math. Econom.}, \textbf{43}, 116--120.

%\bibitem{hao:tang:2012} 
%{\sc Hao, X., Tang, Q.}\ (2012)
%Asymptotic ruin probabilities for a bivariate L\'{e}vy-Driven risk model with heavy-tailed claims and risky investments.
%{\em  J. Appl. Probab.}, \textbf{49}, 939--953.

%\bibitem{horn:steutel:1978}
%{\sc Horn, R.A., Steutel, F.W.}\ (1978)
%On multivariate infinitely divisible distributions.
%{\em Stoch. Process. Appl.}, \textbf{6}, 139--151.

%\bibitem{hu:jiang:2013}
%{\sc Hu, Z., Jiang, B.}\  (2013)
%On joint ruin probabilities of a two-dimensional risk model with constant interest rate. 
%{\em J. Appl. Probab.}, \textbf{50}, no. 2, 309--322.

%\bibitem{hult:lindskog:2006}
%{\sc Hult, H., Lindskog, F.}\ (2006)
%Heavy-tailed insurance protfolios: buffer capital and ruin probabilities.
%{\em Technical report}.

%\bibitem{hult:samorodnitsky:2008}
%{\sc Hult, H., Samorodnitsky, G.}\ (2008)
%Tail probabilities for infinite series of regularly varying random vectors.
%{\em Bernoulli}, \textbf{14}, no. 3, 838--864.
%JJJJJJJJJJJJJJJJJJJJJJJJJJJJJJJJJJJJJJJJJJJJJJJJJJJJJJJJJJJJJ

%\bibitem{jaune:ragulina:siaulys:2018}
%{\sc Jaune, E., Ragulina, O., {\v S}iaulys, J.}\ (2018)
%Expectation of the truncated randomly weighted sums with dominatedly varying summands.
%{\em Lith. Math. J.}, \textbf{58}, no. 4, 421--440.

%\bibitem{jaune:siaulys:2022}
%{\sc Jaune, E.,  {\v S}iaulys, J.}\ (2022)
%Asymptotic risk decomposition for regularly varying distributions with tail dependence.
%{\em Appl. Math. Comput.}, \textbf{427}, no. 127164.

\bibitem{jia:chen:cheng:2025}
{\sc Jia, Y., Chen, Z., Cheng, D.}\  (2025)
Asymptotics for ruin probabilities of a bidimensional risk model with a random number of delayed claims. 
{\em Commun. Stat. Th. Meth.}, \textbf{54}, no. 10, 2990--3007.

%\bibitem{jiang:tang:2011}
%{\sc Jiang, T., Tang, Q.}\  (2011)
%The product of two dependent random variables with rebularly varying or rapidly varyng tails. 
%{\em Stat. Probab. Lett.}, \textbf{81}, 957--961.

%\bibitem{jiang:gao:wang:2014}
%{\sc Jiang, T., Gao, Q., Wang, Y.}\  (2014)
%Max-sum equivalence of conditionally dependent random variables. 
%{\em Stat. Probab. Lett.}, \textbf{84}, 60--66.

\bibitem{jiang:wang:chen:xu:2015}
{\sc Jiang, T., Wang, Y., Chen, Y., Xu, H.}\  (2015)
Uniform asymptotic estimate for finite-time ruin probabilities of a time-dependent bidimensional renewal model. 
{\em Insur. Math. Econom.}, \textbf{64}, 45--53.

%\bibitem{ji:wang:yan:cheng:2023}
%5{\sc Ji, X., Wang, B., Yan, J., Cheng, D.}\  (2023)
%Asymptotic estimates for finite-time ruin probabilities in a generalized dependent bidimensional risk model with CMC simulations. 
%{\em J. Industr. Manag. Optim.}, \textbf{19}, 2140--2155.
%KKKKKKKKKKKKKKKKKKKKKKKKKKKKKKKKKKKKKKKKKKKKKKKK

%\bibitem{kluppelberg:resnick:2008}
%{\sc Kl{\"{u}}ppelberg, C., Resnick, S.}\ (2008)
%The pareto copula, aggregation of risks and the emperor's socks.
%{\em J. Appl. Probab.}, \textbf{45}, 67--84.

\bibitem{ko:tang:2008} 
{\sc Ko, B.W., Tang Q.H.} \ (2008)
Sums of dependent non-negative random variables with subexponential tails.
{\em J. Appl. Probab.}, \textbf{45}, 85--94.

%\bibitem{konstantinides:2018} 
%{\sc Konstantinides, D.G.} \ (2018)
%{\em Risk Theory. A Heavy Tail Approach.}
%World Scientific, New Jersey.

%\bibitem{konstantinides:2025} 
%{\sc Konstantinides, D.G.} \ (2025)
%Infinite-time ruin probability of a multivariate renewal risk model with Brownian perturbations.
%{\em Preprint, arXiv:2510.11229}.

%\bibitem{konstantinides:li:2016}
%{\sc Konstantinides, D.G., Li, J.}\ (2016)
%Asymptotic ruin probabilities for a multidimensional renewal risk model with multivariate regularly varying claims.
%{\em Insur. Math. Econom.}, \textbf{69}, 38--44.

\bibitem{konstantinides:liu:passalidis:2025} 
{\sc Konstantinides, D.G., Liu, J., Passalidis, C.D.} \ (2026)
Uniform asymptotics for a multidimensional renewal risk model with multivariate subexponential claims.
{\em Scand. Act. J.}, 2026(6), 609--629.

%\bibitem{konstantinides:passalidis:2025a} 
%{\sc Konstantinides, D.G., Passalidis, C.D.} \ (2025)
%Background risk model in presence of heavy tails under dependence.
%{\em Non. Anal. Mod. Contr.}, \textbf{30}.

%\bibitem{konstantinides:passalidis:2024b} 
%{\sc Konstantinides, D.G., Passalidis, C.D.} \ (2024)
%Closure properties and heavy tails: random vectors in the presence of dependence
%{\em Preprint, arXiv:2402.09041}.



\bibitem{konstantinides:passalidis:2025h} 
{\sc Konstantinides, D.G., Passalidis, C.D.} \ (2025)
Heavy-tailed random vectors: theory and applications.
{\em Preprint, arXiv:2503.12842}.

\bibitem{konstantinides:passalidis:2025c} 
{\sc Konstantinides, D.G., Passalidis, C.D.} \ (2025)
A new approach in two-dimensional heavy-tailed distributions.
{\em Ann. Actuar. Scien.}, \textbf{19}, no.2, 317 -- 349.

%\bibitem{konstantinides:passalidis:2024d} 
%{\sc Konstantinides, D.G., Passalidis, C.D.} \ (2025)
%Positively decreasing and related distributions under dependence.
%{\em Theory Probab. Appl.}, \textbf{70}, no.3, 461--486.

\bibitem{konstantinides:passalidis:2024j} 
{\sc Konstantinides, D.G., Passalidis, C.D.} \ (2025)
Uniform asymptotic estimates for ruin probabilities of a multidimensional risk model with c\'{a}dl\'{a}g returns and multivariate heavy-tailed claims.
{\em  Insur. Math. Econom.}, \text{125}, 103148.

\bibitem{konstantinides:passalidis:2025o} 
{\sc Konstantinides, D.G., Passalidis, C.D.} \ (2025)
Multivariate subexponentiality and interplay of insurance and financial risks in a renewal risk model.
{\em Preprint, arXiv:2510.17377}.

\bibitem{konstantinides:passalidis:2024g} 
{\sc Konstantinides, D.G., Passalidis, C.D.} \ (2026)
Random vectors in the presence of a single big jump.
{\em J. Theor. Probab. to appear}.

\bibitem{konstantinides:tang:tsitsiashvili:2002}
{\sc Konstantinides, D., Tang, Q., Tsitsiashvili, G.}\ (2002)
Estimates for the ruin probability in the classical risk model with
constant interest force in the presence of heavy tails.
{\em Insur. Math. Econom.}, \textbf{31}, 447--460.

%LLLLLLLLLLLLLLLLLLLLLLLLLLLLLLLLLLLLLLLLLLLLLLLL

\bibitem{lehmann:1966}
{\sc Lehmann, E.L.}\ (1966)
Some concepts of dependence. 
{\em Ann. Math. Stat.}, \textbf{37}, 1137--1153.

\bibitem{leipus:siaulys:2020}
{\sc Leipus, R.,  {\v S}iaulys, J.}\ (2020)
On a closure property of convolution equivalent class of distributions.
{\em J. Math. Anal. Appl.}, \textbf{490}, no. 124226.

%\bibitem{leipus:paukstys:siaulys:2021}
%{\sc Leipus, R., Pauk\v s tys, S., {\v S}iaulys, J.}\ (2021)
%Tails of higher-order moments of sums with heavy-tailed  increments and application to the Haezendonck-Goovaerts risk measure.
%{\em Stat. Probab. Lett.}, \textbf{170}, no. 108998.

%\bibitem{leipus:siaulys:dirma:zove:2023}
%{\sc Leipus, R., {\v S}iaulys, J., Dirma, M., Zov{\.e}, R.}\ (2023)
%On the distribution-tail behaviour of the product of normal random variables.
%{\em J. Inequal. Appl.}, \textbf{32}, no. 1, DOI:10.1186/s13660-023-02941-1.

\bibitem{leipus:siaulys:konstantinides:2023}
{\sc Leipus, R., \v{S}iaulys, J., Konstantinides, D.G.}\ (2023)
{\em Closure Properties for Heavy-Tailed and Related Distributions: An Overview.}
Springer Nature, Cham Switzerland.

%\bibitem{leipus:siaulys:vareikaite:2019}
%{\sc Leipus, R., {\v S}iaulys, J., Vareikaite, I.}\ (2019)
%Tails of higher-order moments with dominatedly varying summands. 
%{\em Lith. Math. J.}, \textbf{59}, no. 3, 389--407.

%\bibitem{leipus:surgailis:2007}
%{\sc Leipus, R., Surgailis, D.}\ (2007)
%On long-range dependence in regenerative processes based on a general ON/OFF scheme. 
%{\em J. Appl. Probab.}, \textbf{44}, 379--392.

%\bibitem{leslie:1989}
%{\sc Leslie, J.R.}\ (1989)
%On the non-closure under convolution of the subexponential family. 
%{\em J. Appl. Probab.}, \textbf{26}, 58--66.  

\bibitem{li:2013}
{\sc Li, J.}\ (2013)
On pairwise quasi-asymptotically independent random variables and  their applications. 
{\em Stat. Prob. Lett.}, \textbf{83}, 2081--2087.

\bibitem{li:2016}
{\sc Li, J.}\ (2016)
Uniform asymptotics for a multidimensional time-dependent risk model with multivariate regularly varying claims and stochastic return. 
{\em Insur. Math. Econom.}, \textbf{71}, 195--204.

%\bibitem{li:2017}
%{\sc Li, J.}\ (2017)
%A note on the finite-time ruin probabilities of a renewal risk model with Brownian perturbation. 
%{\em Stat. Prob. Lett.}, \textbf{127}, 49--55.

%\bibitem{li:2018a}
%{\sc Li, J.}\ (2018)
%A revisit to asymptotic ruin probabilities of a two-dimesional renewal risk model. 
%{\em Stat. Probab. Lett.}, \textbf{140}, 23--32.

%\bibitem{li:2018b}
%{\sc Li, J.}\ (2018b)
%On the joint tail behavior of randomly weighted sums of heavy-tailed random variables.
%{\em J. Mutlivar. Anal.}, \textbf{164}, 40--53.

%\bibitem{li:2022a}
%{\sc Li, J.}\ (2022)
%Asymptotic results on marginal expected shortfalls for dependent risks.
%{\em Insur. Math. Econom.}, \textbf{102}, 146--168.

%\bibitem{li:2023a}
%{\sc Li, J.}\ (2023)
%Asymptotic ruin probabilities for a two-dimensional risk model with dependent claims and stochastic return. 
%{\em Comm. Stat. Th. Meth.}, \textbf{53}, no.16, 5773--5784.

\bibitem{li:2023b}
{\sc Li, J.}\ (2023)
Asymptotic ruin probabilities for a renewal risk model with a random number of delayed claims. 
{\em J. Indust. Manag. Optim.}, \textbf{19}, no. 6, 3840--3853.

%\bibitem{li:sun:2009}
%{\sc Li, H., Sun, Y.}\ (2009)
%Tail dependence for heavy-tailed scale mixtures of multivariate distributions.
%{\em J. Appl. Probab.}, \textbf{46}, 925--937.

%\bibitem{li:yang:2015}
%{\sc Li, J., Yang, H.}\ (2015)
%Asymptotic ruin probabilities for a bidimensional renewal risk model with constant interest rate and dependent claims.
%{\em J. Math. Anal. Appl.}, \textbf{426}, no. 1, 247--266.

%\bibitem{li:tang:2015}
%{\sc Li, J., Tang, Q.} \ (2015)
%Interplay of insurance and financial risks in a discrete time model with strongly regular variation.
%{\em Bernoulli} 21, no.3, 1800--1823

%\bibitem{li:tang:wu:2010} 
%{\sc Li, J., Tang, Q., Wu, R.}\ (2010)
%Subexponential tails of discounted aggregate claims in a time-dependent renewal risk model. 
%{\em Adv. Appl. Probab.}, \textbf{42}, no. 4, 1126-1146.

%\bibitem{liu:gao:dong:2024}
%{\sc Liu, X., Gao, Q., Dong, Z.}\ (2024)
%Asymptotics for a diffusion-perturbed risk model with dependence structures, 
%constant interest force, and a random number of delayed claims. 
%{\em Stochastic Models}, \textbf{40}, no. 1, 97--122.

\bibitem{liu:gao:chen:2024}
{\sc Liu, X., Gao, Q., Chen, Y.}\ (2024)
Uniform asymptotics for a risk model with constant force of interest and a random number of delayed claims. 
{\em Stochastics}, \textbf{97}, no. 7, 1817--1839.

%\bibitem{liu:yang:2021}
%{\sc Liu, J., Yang, Y.}\ (2021)
%Asymptotics for systemic risk with dependent heavy-tailed losses. 
%{\em ASTIN bull.}, \textbf{51}, no. 2, 571--605.

%\bibitem{lu:2012} 
%{\sc Lu, D.}\ (2012)
%Lower bounds of large deviation for sums of long-tailed claims in a multi-risk model. 
%{\em Stat. Probab. Lett.}, \textbf{82}, no. 7, 1242-1250.

\bibitem{lu:qin:yuan:2025}
{\sc Lu, D., Qin, X., Yuan, M.}\ (2026)
Uniform asymptotics for a time-dependent bidimensional delay-claim risk model  with stochastic return and dependent subexponential claims. 
{\em Adv.  Appl. Probab.}, \textbf{58}, no. 2, 745--781.

%\bibitem{lu:yuan:2022}
%{\sc Lu, D., Yuan, M.}\ (2022)
%Asymptotic finite-time ruin probabilities for a bidimensional delay-claim risk model with subexponential claims. 
%{\em Meth. Comp. Appl. Probab.}, \textbf{24}, 2265-2286.

%ΜΜΜΜΜΜΜΜΜΜΜΜΜΜΜΜΜΜΜΜΜΜΜΜΜΜΜΜΜΜΜΜΜΜΜΜΜΜΜΜΜΜΜΜΜΜ

%\bibitem{mcneil:frey:embrechts:2005} 
%{\sc McNeil, A.J., Frey, R., Embrechts, P.}\ (2005)
%{\em Quantitative Risk Management: Concepts, Techniques and Tools}, 
%Princeton: Princeton Univ. Press.

%\bibitem{matuszewska:1964}
%{\sc Matuszewska, W.}\ (1964)
%On generalization of regularly increasing functions.
%{\em Studia Mathematica}, 24, 271--279.

%\bibitem{maulik:resnick:2004} 
%{\sc Maulik, K., Resnick, S.I.}\ (2004)
%Characterizations and examples of hidden regular variation.
%{\em Extremes}, \textbf{7}, 31--67.

%\bibitem{mikosch:nagaev:1998}
%{\sc Mikosch, T., Nagaev, A.V.}\ (1998)
%Large deviations of heavy-tailed sums with applications in insurance. 
%{\em Extremes}, \textbf{1}, no. 1, 81--110.

%\bibitem{mikosch:samorodnitsky:2000b} 
%{\sc Mikosch, T., Samorodnitsky, G.}\ (2000) 
%Ruin probability with claims modeled by a stationary ergodic stable process. 
%{\em Ann. Probab.}, \textbf{28}, no. 4, 1814--1851.

%\bibitem{mikosch:wintenberger:2024} 
%{\sc Mikosch, T., Wintenberger, O.}\ (2024)
%Extreme value theory for time series: Models with power-law tails.
%{\em Springer Nature}, Cham Switzerland.

%NNNNNNNNNNNNNNNNNNNNNNNNNNNNNNNNNNNNNNNNNNNNNNNNNNNNNNNNN

%\bibitem{ng:tang:yan:yang:2003}
%{\sc Ng, K.W., Tang, Q., Yan, J., Yang, H.}\ (2003)
%Precise large deviation for the prospective-loss process.
%{\em J. Appl. Probab.}, \textbf{40}, no. 2, 391--400.

%\bibitem{ng:tang:yang:2002}
%{\sc Ng, K.W., Tang, Q., Yang, H.}\ (2002)
%Maxima of sums of heavy-tailed random variables.
%{\em ASTIN bull.}, \textbf{32}, no.1, 43--55.

%OOOOOOOOOOOOOOOOOOOOOOOOOOOOOOOOOOOOOOOOOOOOOOOO

%PPPPPPPPPPPPPPPPPPPPPPPPPPPPPPPPPPPPPPPPPPPPPPPPP
%\bibitem{pakes:2004}
%{\sc Pakes, A.G.} \ (2004)
%Convolution equivalence and infinite divisibility.
%{\em J. Appl. Probab.}, \textbf{41}, 407--424.

%\bibitem{pakes:2007}
%{\sc Pakes, A.G.} \ (2007)
%Convolution equivalence and infinite divisibility: Corrections and corollaries.
%{\em J. Appl. Probab.}, \textbf{44}, 295--305.

%\bibitem{palmowski:pojer:thonhauser:2025}
%{\sc Palmowski, Z., Pojer, S., Thonhauser, S.} \ (2025)
%Exact asymptotics of ruin probabilities with linear Hawkes arrivals.
%{\em Stoch. Proc. Appl.}, \textbf{182}, 104571.

%\bibitem{passalidis:2025}
%{\sc Passalidis, C.D.} \ (2025)
%Multivariate strong subexponential distributions: Properties and Applications.
%{\em Preprint, arXiv:2503.22267}.

%RRRRRRRRRRRRRRRRRRRRRRRRRRRRRRRRRRRRRRRRRRRRRRRR
\bibitem{resnick:2007}
{\sc Resnick, S.}\ (2007) 
{\em Heavy-Tail Phenomena. Probabilistic and Statistical Modeling.} 
Springer, New York.

%\bibitem{resnick:samorodnitsky:2015}
%{\sc Resnick, S., Samorodnitsky, G.}\ (2015) 
%Tauberian theory for multivariate regularly varying  distributions with application to preferential attachment networks. 
%{\em Extremes}, \textbf{18}, 349--367.

%\bibitem{resnick:willekens:1991} 
%{\sc Resnick, S. I., Willekens, E.}\   (1991)
%Moving averages with random coefficients and random coefficient autoregressive models. \emph{ Stochastic Models},  7(4), 511--525.

%\bibitem{rolski:schmidli:schmidt:teugels:1999}
%{\sc Rolski, T., Schmidli, H., Schmidt, V., Teugels, J.L.}\ (1999) 
%{\em Stochastic Processes for Insurance and Finance} 
%Wiley, Chichester.

%\bibitem{ross:1983} 
%{\sc Ross, S.M.} (1983)
%{\em Stochastic Processes.}
%Wiley, New York.

%SSSSSSSSSSSSSSSSSSSSSSSSSSSSSSSSSSSSSSSSSSSSSSSSSSSSSSSSSSSS

\bibitem{samorodnitsky:2016} 
{\sc Samorodnitsky, G.}\ (2016) 
{\em Stochastic Processes and Long Range Dependence.}
Springer, New York.

%\bibitem{samorodnitsky:resnick:towsley:davis:willis:wan:2016}
%{\sc Samorodnitsky, G., Resnick, S., Towsley, D., Davis, R., Willis, A., Wan, P.}\ (2016) 
%Nonstantard regular variation of in-degree and out-degreee in the preferential attachment model. 
%{\em J. Appl. Probab.}, \textbf{53}, 146--161.

\bibitem{samorodnitsky:sun:2016} 
{\sc Samorodnitsky, G., Sun, J.}\ (2016) 
Multivariate subexponential distributions and their applications. 
{\em Extremes}, \textbf{19}, no. 2, 171--196.

%\bibitem{samorodnitsky:taqqu:1994} 
%{\sc Samorodnitsky, G., Taqqu, M.}\ (1994) 
%{\em Stable Non-Gaussian Random Processes: Stochastic Models with Infinite Variance}
%Chapman and Hall, New York.

%\bibitem{schmidli:2017}
%{\sc Schmidli, H.}\ (2017) 
%{\em Risk Theory} 
%Springer, Cham.

%\bibitem{sgibniev:1988}
%{\sc Sgibnev, M. S.}\ (1988) 
%Banach algebras of measures of class $\mathcal{S}(\gamma)$.
%{\em Siberian Math. J.}, \textbf{29}, no. 4, 647--655.

%\bibitem{shen:du:2023}
%{\sc Shen, X., Du, K.}\ (2023) 
%Uniform approximation for the tail behavior of bidimensional randomly weighted sums.
%{\em Methodol. Comput. Appl. Probab.}, \textbf{25}, no. 26.

%\bibitem{shen:ge:fu:2020}
%{\sc Shen, X., Ge, M., Fu, K.A.}\ (2020) 
%Approximation of the tail probabilities for bidimensional randomly weighted sums with dependent components.
%{\em Probab. Eng. Inf. Sci.}, \textbf{34}, no. 1, 112--130.

\bibitem{shen:yuan:lu:2022}
{\sc Shen, X., Yuan, M., Lu, D.}\ (2022) 
Uniform asymptotics for ruin probabilities of multidimensional risk models with stochastic returns and regular variation claims.
{\em Commun. Stat. Th. Method.}, \textbf{52}, no. 19, 6878--6895.

%\bibitem{shimura:watanabe:2005}
%{\sc Shimura, T., Watanabe, T.}\ (2005)
%Infinite divisibility and generalized subexponentiality. 
%{\em Bernoulli}, \textbf{11}, 445--469.

\bibitem{stein:1946}
{\sc Stein, C.}\  (1946) 
A note on comulative sums.
{\em Ann. Math. Stat.}, \textbf{4}, 498--499.

%TTTTTTTTTTTTTTTTTTTTTTTTTTTTTTTTTTTTTTTTTTTTTTTT

\bibitem{tang:2004}
{\sc Tang, Q.}\ (2004)
Asymptotics for the finite-time ruin probability in the renewal risk model with consistent variation.
{\em Stoch. Models}, \textbf{20}, no. 3, 281--297.

\bibitem{tang:2006}
{\sc Tang, Q.}\ (2006)
The subexponentiality of products revisited.
{\em Extremes}, \textbf{9}, 231--241.

%\bibitem{tang:2007}
%{\sc Tang, Q.}\ (2007)
%Heavy tails of discounted aggregate claims in the continuous-time renewal model. 
%{\em J. Appl. Probab.}, \textbf{44}, 285--294.

%\bibitem{tang:2008}
%{\sc Tang, Q.}\ (2008)
%Insensitivity to negative dependence of asymptotic tail probabilities of sums and maxima of sums.
%{\em Stoch. Anal. Appl.}, \textbf{26}, 435--450.

%\bibitem{tang:2008b}
%{\sc Tang, Q.}\ (2008) 
%From light tails to heavy tails through multiplier.
%{\em Extremes} \textbf{11}, 379--391.

\bibitem{tang:tsitsiashvili:2003}
{\sc Tang, Q., Tsitsiashvili, G.}\ (2003)
Precise estimates for the ruin probability in finite horizon in a discrete-time model with heavy-tailed insurance and financial risks.
{\em Stoch. Process. Appl.}, \textbf{108}, 299--325.

%\bibitem{tang:yang:2019}
%{\sc Tang, Q., Yang, Y.}\ (2019)
%Interplay of insurance and financial risks in a stochastic environment. 
%{\em Scand. Actuar. J.}, no. 5, 432--451.

\bibitem{tang:yuan:2014}
{\sc Tang, Q., Yuan, Z.}\ (2014) 
Randomly weighted sums of subexponential random variables with application to capital allocation.
{\em Extremes}, \textbf{17}, 467--493.

\bibitem{tang:yuan:2016}
{\sc Tang, Q., Yuan, Z.}\ (2016) 
Random difference equations with subexponential innovations.
{\em Scien. China Math.}, \textbf{59}, 2411--2426.

%\bibitem{teugels:1975}
%Teugels, J.L., (1975)
%The class of subexponential distributions.
%\emph{ Ann. Probab.}, {\bf 3}, 1000--1011.

%WWWWWWWWWWWWWWWWWWWWWWWWWWWWWWWWWWWWWWWWWWWWWWWWWW

%\bibitem{wang:tang:2006}
%{\sc Wang, D., Tang, Q.}\ (2006)
%Tail probabilities of randomly weighted sums of random variables with dominated variation. 
%{\em Stoch. Models}, \textbf{22}, 253--272.

\bibitem{wang:2011}
{\sc Wang, K.}\ (2011)
Randomly weighted sums of dependent subexponential random variables.
{\em Lith. Math. J.}, \textbf{51}, no. 4, 573--586.

%\bibitem{wang:su:yang:2024}
%{\sc Wang, H., Su, Q., Yang, Y.}\ (2024)
%Asymptotics for ruin probabilities in a bidimensional discrete-time risk model with dependence and consistently varying net losses. 
%{\em Stochastics}, \textbf{96}, no. 1, 667--695.

\bibitem{waters:papatriandafylou:1985}
{\sc Waters, H.R., Papatriandafylou, A.}\ (1985)
Ruin probabilities allowing for delay in claims settlement. 
{\em Insur. Math. Econom.}, \textbf{4}, no. 2, 113--122.

%XXXXXXXXXXXXXXXXXXXXXXXXXXXXXXXXXXXXXXXXXXXXXXXXXXXXXXXXXX

\bibitem{xu:cheng:wang:cheng:2018}
{\sc Xu, H., Cheng, F., Wang, Y., Cheng, D.}\ (2018) 
A necessary and sufficient condition for the subexponentiality of product convolution. 
{\em Adv. Appl. Probab.} \textbf{50}, no.1,  57--73.

%\bibitem{xu:shen:wang:2025}
%{\sc Xu, C., Shen, X., Wang, K.}\ (2025) 
%The finite-time ruin probabilities of a dependent bidimensional risk model with subexponential claims and Brownian perturbations. 
%{\em Non. Anal. Mod. Contr.}, \textbf{30}.

%YYYYYYYYYYYYYYYYYYYYYYYYYYYYYYYYYYYYYYYYYYYYYYYYYYYYYYYY
%\bibitem{yang:chen:yuen:2024}
%{\sc Yang, H., Chen, S., Yuen, C.}\ (2024)
%Asymptotics for the joint tail probability of bidimensional randomly weighted sums with applications to insurance.
%{\em Sci. China Math.}, \textbf{67}, 163--186.

%\bibitem{yang:cheng:zhang:2025} 
%{\sc Yang, Y., Cheng, B., Zhang, Z.} \ (2025)
%A revisit to tail risk measures in the presence of bivariate regularly varying tailed insurance and financial risks.
%{\em Non. Anal. Mod. Contr.} \textbf{30}, no.5, 892-912.

%\bibitem{yang:fan:yuen:2023}
%{\sc Yang, Y., Fan, Y., Yuen, K.C.}\ (2023)
%Ruin in a continuous-time risk model with arbitrarily dependent insurance and financial risks triggered by systematic factors.
%{\em Scand. Actuar. J.}, 2024(4), 361--382.

%\bibitem{yang:gao:li:2016}
%{\sc Yang, H., Gao, W., Li, J.}\ (2016)
%Asymptotic ruin probabilities for a discrete-time risk model with dependent insurance and  financial risks.
%{\em Scand. Actuar. J.} no. 1, 1--17.

%\bibitem{yang:li:2014}
%{\sc Yang, H., Li, J.}\ (2014)
%Asymptotic finite-time ruin probabilities for a bidimensional renewal risk model with constant interest force and dependent subexponential claims.
%{\em Insur. Math. Econom.}, \textbf{58}, 185--192.

\bibitem{yang:li:2019}
{\sc Yang, H., Li, J.}\ (2019)
On asymptotic finite-time ruin probabilities for a renewal risk model with subexponential main claims and delayed claims.
{\em Stat. Probab. Lett.}, \textbf{58}, 185--192. 

%\bibitem{yang:zhang:jiang:cheng:2015}
%{\sc Yang, Y., Zhang, Z., Jiang, T., Cheng, D.}\ (2015)
%Uniformly asymptotic behavior of ruin probabilities in a time-dependent renewal risk model with stochastic return.
%{\em J. Comput. Appl. Math.} {\bf 287}, 32--43.

%\bibitem{yang:konstantinides:2015}
%{\sc Yang, Y., Konstantinides, K.}\ (2015)
%Asymptotic for ruin probabilities in a discrete-time risk model with dependent financial and insurance risks.
%{\em Scand. Actuar. J.}, \textbf{8}, 641--659.

%\bibitem{yang:leipus:siaulys:2012}
%{\sc Yang, Y., Leipus, R., \v{S}iaulys, J.}\ (2012)
%Tail probability of randomly weighted sums of subexponential random variables under a dependence structure. 
%{\em Stat. Probab. Lett.}, \textbf{82}, 1727--1736.

\bibitem{yang:su:2023}
{\sc Yang, Y., Su, Q.}\ (2023)
Asymptotic behavior of ruin probabilities in a multidimensional risk model with investment and multivariate regularly varying claims.
{\em J. Math. Anal. Appl.}, \textbf{525}, 127319.

%\bibitem{yang:wang:2010}
%{\sc Yang, Y., Wang, Y.}\ (2010)
%Asymptotics for ruin probability of some negatively dependent risk models with a constant interest rate and dominatedly-varying tailed claims.
%{\em Stat. Probab. Lett.}, \textbf{80}, no.3-4, 143--154.

%\bibitem{yang:wang:2013}
%{\sc Yang, Y., Wang, Y.}\ (2013)
%Tail behavior of the product of two dependent random variables with applications to risk theory.
%{\em Extremes}, \textbf{16}, 55--74.

%\bibitem{yang:wang:konstantinides:2014}
%{\sc Yang, Y., Wang, K., Konstantinides, D.G.}\ (2014)
%Uniform asymptotics for discounted aggregate claims in dependent risk models.
%{\em J. Appl. Probab.}, \textbf{51}, no. 3, 669--684.

%\bibitem{yang:wang:leipus:siaulys:2011}
%{\sc Yang, Y., Wang, K., Leipus, R., \v{S}iaulys, J.}\ (2011)
%Tail behavior of sums and maxima of sums of dependent subexponential random variables.
%{\em Acta Appl. Math.}, \textbf{114}, 219--231.

%\bibitem{yang:wang:leipus:siaulys:2013}
%{\sc Yang, Y., Wang, K., Leipus, R., \v{S}iaulys, J.}\ (2013)
%A note on the max-sum equivalence of randomly weighted sums of heavy-tailed random variables.
%{\em Nonlin. Anal. Mod. Contr.}, \textbf{18}, no.4, 519--525.

\bibitem{yang:wang:zhang:2021}
{\sc Yang, Y., Wang, X., Zhang, Z.}\ (2021)
Finite-time ruin probability of a perturbed risk model with dependent main and delayed claims.
{\em Nonlin. Anal. Mod. Contr.}, \textbf{26}, no. 5, 801--820.

%\bibitem{yang:yuen:liu:2018}
%{\sc Yang, Y., Yuen, K.C., Liu, J.-F.}\ (2018) 
%Asymptotics for ruin probabilities in L\'{e}vy-driven risk models with heavy-tailed claims. 
%{\em J. Ind. Manag. Optim.}, \textbf{14}, 231--247.

%\bibitem{yang:wang:liu:zhang:2019}
%{\sc Yang, Y., Wang, K., Liu, J., Zhang, Z.}\ (2019) 
%Asymptotics for a bidimensional risk model with two geometric L\'{e}vy price processes. 
%{\em J. Ind. Manag. Optim.}, \textbf{15}, 481--505.

\bibitem{yi:chen:su:2011}
{\sc Yi, L., Chen Y., Su, C.}\ (2011) 
Approximation of the tail probability of randomly weighted sums of dependent random variables with dominated variation. 
\emph{ J. Math. Anal. Appl.}, {\bf 376}, 365--372.

%\bibitem{yu:wang:cheng:2015}
%{\sc Yu, C., Wang, Y., Cheng, D.}\ (2015)
%Tail behavior of the sums of dependent and heavy-tailed random variables.
%{\em J. Korean Stat. Soc.}, \textbf{44}, 12--27.

\bibitem{yuan:lu:2023}
{\sc Yuan, M., Lu, D.}\ (2023)
Asymptotics for a time-dependent by-claim model with dependent subexponential claims.
{\em Insur. Math. Econom.}, \textbf{112}, 120--141.

\bibitem{yuan:lu:fu:2025}
{\sc Yuan, M., Lu, D., Fu, Y.}\ (2026)
Asymptotics for a multidimensional risk model with a random number of delayed claims and multivariate 
regularly varying distribution.
{\em Adv. Appl. Probab.}, \textbf{58}, no. 2, 782--811.

\bibitem{yuen:guo:2001}
{\sc Yuen, K.C., Guo, J.Y.}\ (2001)
Ruin probabilities for time-correlated claims in the compound binomial model.
{\em Insur. Math. Econom.}, \textbf{29}, no. 1, 47--57.

\bibitem{yuen:guo:ng:2005}
{\sc Yuen, K.C., Guo, J.Y., Ng, K.W.}\ (2005)
On ultimate ruin in a delayed-claims risk model.
{\em J. Appl. Probab.}, \textbf{42}, no. 1, 143--174.

%ZZZZZZZZZZZZZZZZZZZZZZZZZZZZZZZZZZZZZZZZZZZZZZZZZZZZZZZZZZZz

%\bibitem{zhou:wang:wang:2012}
%{\sc Zhou, M., Wang, K., Wang, Y.}\ (2012) 
%Estimates for the finite-time ruin probability with insurnace and financial risks.
%{\em Acta Math. Appl. Sin.}, \textbf{28}, no. 4, 795--806.

%\bibitem{zhu:li:2012}
%{\sc Zhu, L., Li, H.}\ (2012) 
%Tail distortion risk and its asymptotic analysis.
%{\em Insur. Math. Econom.}, \textbf{51}, no. 1, 115--132.


\end{thebibliography}
\end{document}